\documentclass[12pt,a4paper]{article}
\usepackage{amssymb,amsmath,amsthm,euscript}
\usepackage{leftidx}
\usepackage[matrix,arrow,curve]{xy}

\usepackage[utf8]{inputenc}
\usepackage[english]{babel}

\newcommand{\chara}{\mathop{\mathrm{char}}\nolimits}
\newcommand{\End}{\mathop{\mathrm{End}}\nolimits}

\newcommand{\Hom}{\mathop{\mathrm{Hom}}\nolimits}
\newcommand{\bfhom}{\mathop{\mathrm{\bf hom}}\nolimits}
\newcommand{\bfend}{\mathop{\mathrm{\bf end}}\nolimits}
\newcommand{\bfcohom}{\mathop{\mathrm{\bf cohom}}\nolimits}
\newcommand{\bfcoend}{\mathop{\mathrm{\bf coend}}\nolimits}
\newcommand{\Mat}{\mathop{\mathrm{Mat}}\nolimits}

\newcommand{\id}{\mathop{\mathrm{id}}\nolimits}

\newcommand{\Spec}{\mathop{\mathrm{Spec}}\nolimits}

\newcommand{\DiS}{\mathop{\mathrm{DiS}}\nolimits}
\newcommand{\CoP}{\mathop{\mathrm{CoP}}\nolimits}
\newcommand{\TeP}{\mathop{\mathrm{TeP}}\nolimits}

\newcommand{\ootimes}{\mathbin{\underline{\otimes}}}
\newcommand{\isoright}{\xrightarrow{\smash{\raisebox{-0.65ex}{\ensuremath{\sim}}}}}
\newcommand{\dotsqcup}{\mathbin{\dot\sqcup}}
\newcommand{\dototimes}{\mathbin{\dot\otimes}}
\newcommand{\dotcirc}{\mathbin{\dot\circ}}
\newcommand{\dotbullet}{\mathbin{\dot\bullet}}

\numberwithin{equation}{section}
\renewcommand{\Im}{\mathop{\mathrm{Im}}\nolimits}

\newcommand{\op}{{\mathrm{op}}}
\newcommand{\cop}{{\mathrm{cop}}}

\renewcommand{\le}{\leqslant}
\renewcommand{\ge}{\geqslant}

\newcommand{\la}{\langle}
\newcommand{\ra}{\rangle}

\newcommand{\A}{{\mathcal A}}
\newcommand{\B}{{\mathcal B}}
\newcommand{\C}{{\mathcal C}}
\newcommand{\D}{{\mathcal D}}

\newcommand{\M}{{\mathcal M}}
\newcommand{\N}{{\mathcal N}}
\newcommand{\K}{{\mathcal K}}

\newcommand{\gR}{{\mathfrak R}}
\newcommand{\gS}{{\mathfrak S}}
\newcommand{\gA}{{\mathfrak A}}
\newcommand{\gX}{{\mathfrak X}}
\newcommand{\gU}{{\mathfrak U}}
\newcommand{\gC}{{\mathfrak C}}
\newcommand{\gB}{{\mathfrak B}}

\newcommand{\wt}{\widetilde}
\newcommand{\wh}{\widehat}

\newtheorem{Th}{Theorem}[section]
\newtheorem{Lem}[Th]{Lemma}
\newtheorem{Prop}[Th]{Proposition}
\newtheorem{Cor}[Th]{Corollary}

\newtheorem{Def}[Th]{Definition}
\newtheorem{Rem}[Th]{Remark}
\newcommand{\BB}{{\mathbb{B}}}
\newcommand{\GG}{{\mathbb{G}}}

\newcommand{\ZZ}{{\mathbb{Z}}}
\newcommand{\MM}{{\mathbb{M}}}
\newcommand{\OO}{{\mathbb{O}}}
\newcommand{\NN}{{\mathbb{N}}}
\newcommand{\KK}{{\mathbb{K}}}

\renewcommand{\AA}{{\mathbb{A}}}
\newcommand{\SSS}{{\mathbb{S}}}

\newcommand{\Set}{{\mathbf{Set}}}
\newcommand{\AlgSet}{{\mathbf{AlgSet}}}
\newcommand{\SLAlgSet}{{\mathbf{SLAlgSet}}}
\newcommand{\AffSch}{{\mathbf{AffSch}}}
\newcommand{\SLAffSch}{{\mathbf{SLAffSch}}}
\newcommand{\Vect}{{\mathbf{Vect}}}
\newcommand{\FVect}{{\mathbf{FVect}}}
\newcommand{\GrVect}{{\mathbf{GrVect}}}
\newcommand{\Alg}{{\mathbf{Alg}}}
\newcommand{\CommAlg}{{\mathbf{CommAlg}}}
\newcommand{\FAlg}{{\mathbf{FAlg}}}
\newcommand{\GrAlg}{{\mathbf{GrAlg}}}
\newcommand{\QA}{{\mathbf{QA}}}
\newcommand{\QAsc}{{\mathbf{QA}}_{\mathrm sc}}
\newcommand{\FTA}{{\mathbf{FTA}}}
\newcommand{\FQA}{{\mathbf{FQA}}}
\newcommand{\FQAsc}{{\mathbf{FQA}}_{\mathrm sc}}
\newcommand{\CommQSAe}{{\mathbf{CommQSA}}_{\mathrm even}}
\newcommand{\CommQSAo}{{\mathbf{CommQSA}}_{\mathrm odd}}
\newcommand{\FCommQSAe}{{\mathbf{FCommQSA}}_{\mathrm even}}
\newcommand{\FCommQSAo}{{\mathbf{FCommQSA}}_{\mathrm odd}}

\newcommand{\bfC}{{\mathbf{C}}}
\newcommand{\bfD}{{\mathbf{D}}}
\newcommand{\bfP}{{\mathbf{P}}}
\newcommand{\bfQ}{{\mathbf{Q}}}

\newcommand{\Mon}{\mathop{\mathbf{Mon}}\nolimits}
\newcommand{\cMon}{\mathop{\mathbf{\leftidx{_c}Mon}}\nolimits}
\newcommand{\Comon}{\mathop{\mathbf{Comon}}\nolimits}
\newcommand{\cocComon}{\mathop{\mathbf{\leftidx{_{coc}}Comon}}\nolimits}
\newcommand{\Bimon}{\mathop{\mathbf{Bimon}}\nolimits}
\newcommand{\Lact}{\mathop{\mathbf{Lact}}\nolimits}
\newcommand{\Rep}{\mathop{\mathbf{Rep}}\nolimits}
\newcommand{\Corep}{\mathop{\mathbf{Corep}}\nolimits}
\newcommand{\Lcoact}{\mathop{\mathbf{Lcoact}}\nolimits}
\newcommand{\eval}{\mathop{\mathrm{ev}}\nolimits}
\newcommand{\coev}{\mathop{\mathrm{coev}}\nolimits}

\usepackage{hyperref}

\newcounter{bbcount}[subsection]
\renewcommand{\thebbcount}{{\thesection}.\arabic{subsection}.\arabic{bbcount}}
\newcommand{\bbo}[1]{\noindent\refstepcounter{bbcount}{\bf\thebbcount.}{ \bfseries #1}}
\newcommand{\bb}[1]{\vspace{2mm}\noindent\refstepcounter{bbcount}{\bf\thebbcount.}{ \bfseries #1}}

\title{Quantum Representation Theory and Manin matrices I: finite-dimensional case}
\author{Alexey Silantyev\thanks{aleksejsilantjev@gmail.com}}
\date{}

\begin{document}


\maketitle

\vspace{-5mm}
\begin{center}
{\it Bogoliubov Laboratory of Theoretical Physics, Joint Institute for Nuclear Research, 141980~Dubna, Moscow region, Russia} \\
{\it State University ``Dubna''{}, Universitetskaya st. 19, 141980~Dubna, Moscow region, Russia} \\
\end{center}
\vspace{5mm}

\begin{abstract}
 We construct Quantum Representation Theory which describes quantum analogue of representations in frame of `non-commutative linear geometry' developed by Manin in~\cite{Manin88}. To do it we generalise the internal hom-functor to the case of adjunction with a parameter and construct a general approach to representations of a monoid in a symmetric monoidal category with a parameter subcategory. Quantum Representation Theory is obtained by application of this approach to a monoidal category of some class of graded algebras with Manin product, where the parameter subcategory consists of connected finitely generated quadratic algebras. We formulate this theory in the language of Manin matrices and obtain quantum analogues of direct sum and tensor product of representations. Finally, we give some examples of quantum representations.
\end{abstract}

{\bf Keywords:} quadratic algebras; Manin matrices; quantum groups; non-commu\-ta\-tive spaces; Representation Theory; monoidal categories.


\tableofcontents

\section{Introduction}

In the end of 80s Yuri Manin generalised the notion of finite-dimensional vector (linear) space to the quantum
\footnote{We use the word `quantum' instead of `non-commutative' in the sense of non-commutative spaces by following Manin~\cite{Manin88} and Drinfeld~\cite{Dr}, because sometimes the word `non-commutative' is confusing. For example, the term `non-commutative group' does not mean a non-commutative space with a group structure (i.e. a quantum group), but rather a non-abelian group. Moreover, we want to consider the usual `commutative' spaces as a particular case of the non-commutative ones.
}
case~\cite{Manin88}. He introduced the term `quantum linear spaces' for (finitely generated connected) quadratic algebras considered as objects of the opposite category. The usual finite-dimensional vector spaces correspond to the polynomial algebras. Manin defined four binary operations with quadratic algebras, they where denoted by $\circ$, $\bullet$, $\otimes$, $\ootimes$. He also introduced a duality of quadratic algebras, which coincides with the Koszul duality for the quadratic Koszul algebras (see~\cite{Manin87,Manin88,ManinBook91}).

The most important binary operation for our purposes is `$\circ$', which we call `Manin white product' or just `Manin product'. As it was shown in~\cite{Manin88} the category of quantum linear spaces with these product is a closed monoidal category and that the internal hom can be written via the Koszul duality and the binary operation `$\bullet$', which is called `black Manin product'. This gives a construction of the internal hom-object as a quadratic algebra generated by entries of some matrix.

It turned out that the description of the internal hom-object for the polynomial algebras gives some non-trivial commutation relations for the entries of this matrix. The matrices (over a non-commutative ring in general) satisfying these commutation relations were called `Manin matrices' in~\cite{CF}, where the properties of these matrices and their determinants were applied to Talalaev's formula~\cite{T} and many other questions. A lot of interesting properties were studied in details in~\cite{CFR}. The theory of Manin matrices was applied to Sugawara operators in~\cite{CM} and to some elliptic integrable systems in~\cite{RST}.

In the work~\cite{qManin} the properties of Manin matrices were generalised for the $q$-case, where $q$ is a parameter of the deformation of the polynomial algebras. The MacMahon Master Theorem and Cayley--Hamilton--Newton identities for the $q$-Manin matrices were obtained in~\cite{GLZ,FH07,FH08} and \cite{IO}. Some properties and applications of the super-version of Manin matrices were derived in~\cite{MR}. Manin matrices for the series $B,C,D$ appeared in the book~\cite{MolevSO}. Note also that the super-case with multi-parametric deformation was considered by Manin~\cite{Manin89,ManinBook91,Manin92}.

The notion of Manin matrices was generalised in the work~\cite{S} in terms of idempotent operators (in~\cite{IO} these matrices were called `half-quantum matrices'). It was shown that any quadratic algebra can be given by a matrix idempotent $A$ acting in a tensor product of a vector space. This algebra was denoted by $\gX_A(\KK)$, where $\KK$ is a basic field. Then it was defined a Manin matrix for a pair of idempotents $A$ and $B$ and showed how this notion can be associated to the pair of quadratic algebras $\gX_A(\KK)$ and $\gX_B(\KK)$.

\vspace{3mm}
The main aim of this paper is to generalise Representation Theory to the case of a quantum representation space. The objects which are to be represented (such as algebras, groups, algebraic groups) are generalised up to `quantum algebras' (these include quantum groups).

To do it we formulate a general representation theory in a monoidal category with some generalisation of internal hom-functor. For such category the representations of a monoid can be functorially identified with left actions of this monoid on some objects. In a particular case this approach gives Quantum Representation Theory: we need to take the monoidal category of quadratic algebras with the Manin product `$\circ$' together with some extension of this category, where the category of quadratic algebras plays the role of parameter subcategory. This extension is needed to include the representations of quantum groups. As a `classical' representation can be given by a square matrix in a basis of the representation space, a quantum representation is given by a Manin matrix for the pair of idempotents $A$ and $B$, where $A=B$. 

In the non-commutative geometry there is a widely known quantum analogue of actions of affine algebraic groups on algebraic sets -- comodule algebras (see e.g.~\cite[Def.~III.7.1]{Kass}). The quantum representations can be considered as their `linear' version, i.e. as a particular case of comodule algebras equipped with an additional structure. The notion of linear representation was generalised for quantum groups in~\cite[Def.~20]{FRT89}. It is a special case of quantum representation, when the quantum representation space is given by the tensor algebra.

This work is mainly devoted to the finite-dimensional case. In the quantum level this means that the algebras $\A$ we consider are finitely generated and, as consequence, all the graded components $\A_k$ of these algebras are finite-dimensional vector spaces. For wider generality we sometimes include the infinite-dimensional (infinitely generated) case into consideration, but here we do not define quantum representations on infinitely generated quadratic algebras.

\vspace{3mm}
The paper is organised as follows. In Section~\ref{sec2} we give basic notations and definitions including theory of monoids and their actions. In Section~\ref{sec3} we describe the general approach to construct representation theory in a monoidal category by using the generalised internal hom. Since it is convenient to work in an opposite category we also define the dual notions: generalised internal cohom, corepresentations etc. Section~\ref{sec4} prepares the necessary notions and facts for the quantum representation theory: the Manin's concept of quantum linear spaces, operations with quadratic algebras, semi-linear spaces and its quantum version, the internal cohom for quadratic algebras and its generalisation. Section~\ref{sec5} is devoted to the description of the quantum representations, to the consideration of classical representations as quantum ones and to binary operations with the quantum representations. In Section~\ref{sec6} we give some examples of quantum representations.

\vspace{3mm}
{\it Acknowledgements}. 
The author is grateful to A.~Isaev, Y.~Manin for useful references and to E.~Patrin for discussions.  The author particularly appreciates advice and guidance from V. Rubtsov.

\section{Preliminaries}
\label{sec2}

\subsection{Vector spaces and algebras}

Fist of all, let us fix some terms and notations that we use in this work.

\bb{Basic field.} Let $\KK$ be an infinite field such that $\chara\KK\ne2$. All the vector spaces will be over $\KK$. For two vector spaces $V$ and $W$ we denote by $V\otimes W=V\otimes_\KK W$ their tensor product over $\KK$. For briefness we say `{\it algebra}' for a unital associative algebra over $\KK$ and we suppose that algebra homomorphisms map unity to unity. The unity of an algebra $\A$ we denote by $1$ or $1_\A$. 

\bb{Notations for morphisms in a category.} We denote the composition of morphisms $f\colon Y\to Z$ and $g\colon X\to Y$ in a category $\bfC$ by $f\cdot g$. We do not use the symbol `$\circ$' for composition, because it will be used for the Manin product.

 For a natural transformation $\theta$ between functors $F$ and $G$ we sometimes omit the subscript of its components $\theta_X\colon FX\to GX$. In particular, for a morphism $f\colon Y\to Z$ we denote by $f_*$ and $f^*$ the maps
\begin{align}
 &f_*\colon\Hom(X,Y)\to\Hom(X,Z), &&f_*(g)=f\cdot g\colon X\to Z, &&g\colon X\to Y, \\
 &f^*\colon\Hom(Z,X)\to\Hom(Y,X), &&f^*(h)=h\cdot f\colon Y\to X, &&h\colon Z\to X. \label{fh}
\end{align}
They are the components of the natural transformations $\Hom(-,Y)\to\Hom(-,Z)$ and $\Hom(Z,-)\to\Hom(Y,-)$ respectively.

\bb{Basis.} A {\it basis} of a vector space $V$ is a map $e\colon I\to V$ from a set $I$ (which may be infinite) such that for any vector space $W$ and any map $f\colon I\to W$ there exists a unique linear map $h\colon V\to W$ satisfying $h\cdot e=f$. The basis $e\colon I\to V$ is denoted by $(e_i)_{i\in I}$ or $(e_i)$, where $e_i$ is the value of $e$ on $i\in I$. In the case of the set $I=\{1,\ldots,n\}$ we also use the notation $(e_i)_{i=1}^n$. One can check that $(e_i)_{i\in I}$ is a basis of $V$ iff the elements $e_i\in V$ are pairwise different, the set $\{e_i\mid i\in I\}=e(I)$ is linear independent and any vector $v\in V$ is a finite linear combination of the form $\sum_i\alpha_i e_i$ where $\alpha_i\in\KK$. Notice that the basis $(e_i)$ is not just a set, because the permutation of the elements $e_i$ gives another basis $(e'_i)$ with the same set $\{e'_i\}=\{e_i\}$.

\bb{Bimodules.} Let $\gR$ and $\gS$ be algebras. By an $(\gR,\gS)$-bimodule we understand a vector space $M$ that has a structures of left $\gR$-module and right $\gS$-module agreed with the structure of the vector space. In other words, $M$ is an $(\gR,\gS)$-bimodule in the ring-theoretic sense satisfying the condition $\alpha m=m\alpha$ $\;\forall\,\alpha\in\KK,m\in M$. We also use the term $\gR$-bimodule or two-sided $\gR$-module for $(\gR,\gR)$-bimodule. For any vector space $V$ the tensor product $\gR\otimes V$ has a structure of $\gR$-bimodule: $r_1(r_2\otimes v)=(r_1\otimes m)r_2=(r_1r_2)\otimes v$, $r_1,r_2\in\gR$, $v\in V$.

\bb{Grading.} \label{bbGr}
 Let $\NN_0$ be the monoid of natural numbers (including $0$). We call an $\NN_0$-graded vector space simply {\it graded vector space}. The $k$-th component of a graded vector space $V$ is denoted by $V_k$, so that $V=\bigoplus\limits_{k\in\NN_0}V_k$. By a {\it graded algebra} we mean an algebra $\A$ with a structure of graded vector space such that $\A_k\A_l\subset\A_{k+l}$ for any $k,l\in\NN_0$. The component $\A_0$ is an algebra, so $\A$ has a structure of an $\A_0$-bimodule. Each component $\A_k$ is its $\A_0$-submodule. Note that the tensor product $\gR\otimes\A$ of an algebra $\gR$ with a graded algebra $\A$ is a graded algebra with the components $(\gR\otimes\A)_k=\gR\otimes\A_k$.

\bb{Connected affinely generated graded algebras.} A graded algebra $\A$ is called {\it connected} if $\A_0=\KK$. The most interesting graded algebras in the projective (non-com\-mutative) geometry are the connected algebras generated by the subspace $\A_1$. Let us call a graded algebra $\A$ {\it affinely generated} if it is generated by the subspace $\A_0\oplus\A_1$, so the connected graded algebras $\A$ generated by $\A_1$ are exactly the connected affinely generated graded algebras. Let us call an affinely generated graded algebra $\A$ {\it semi-connected} if it is isomorphic to the graded algebra $\gR\otimes\B$ for an algebra $\gR$ and a connected affinely generated graded algebra $\B$ (since $\B_0=\KK$ we have $\gR\cong\A_0$, so that $\A\cong\A_0\otimes\B$).

\bb{Quadratic algebras.} By a {\it quadratic algebra} over an algebra $\gR$ we mean an algebra isomorphic to $T_\gR M/\mathcal I$ where
\begin{align}
 T_\gR M:=
\gR\oplus M\oplus\big(M\otimes_{\gR}M\big)\oplus\bigoplus\limits_{k\ge3}\big(\underbrace{M\otimes_\gR\cdots\otimes_\gR M}_{k}\big)
\end{align}
is the tensor algebra of an $\gR$-bimodule $M$ and $\mathcal I$ is a (two-sided) ideal of $T_\gR M$ generated by a subset of $M\otimes_{\gR}M$. Any quadratic algebra $\A=T_\gR M/\mathcal I$ has a structure of graded algebra with the components $\A_0=\gR$, $\A_1=M$, $\A_k=M^{\otimes k}/({\mathcal I}\cap M^{\otimes k})$, where $M^{\otimes k}=M\otimes_\gR\cdots\otimes_\gR M$. Hence any quadratic algebra is affinely generated graded algebra.

In a narrow sense the notion `quadratic algebra' is used for the quadratic algebras over $\KK$. In this case the two-sided module $M$ is a vector space $V$ and $\alpha v=v\alpha$ for any $\alpha\in\KK$ and $v\in V$. Such a quadratic algebra is isomorphic to $TV/I$, where $TV=T_\KK V=\bigoplus\limits_{k\in\NN_0}V^{\otimes k}$ is the tensor algebra of the vector space $V$ and $I$ is the ideal of $TV$ generated by a subspace $R\subset V\otimes V$. These quadratic algebras are exactly the connected quadratic algebras.

Any semi-connected quadratic algebra over $\gR$ is isomorphic to $\gR\otimes(TV/I)\cong T_\gR M/\mathcal I$ where $V$ is a vector space, $\M=\gR\otimes V$, $I$ and $\mathcal I$ are the ideals of $TV$ and $T_\gR M$ generated by a subspace $R\subset V^{\otimes2}$ (over $\KK$ and over $\gR$ respectively).

\bb{Notations for some categories.}
\begin{align*}
 &\Set &&\text{the category of sets,} \\
 &\Vect &&\text{the category of vector spaces ($\KK$-modules),} \\
 &\FVect &&\text{the category of finite-dimensional vector spaces,} \\
 &\GrVect &&\text{the category of graded vector spaces,} \\
 &\Alg &&\text{the category of algebras,} \\
 &\FAlg &&\text{the category of finite-dinetional algebras,} \\
 &\CommAlg &&\text{the category of commutative algebras,} \\
 &\GrAlg &&\text{the category of graded algebras,} \\
 &\QA &&\text{the category of quadratic algebras over $\KK$,} \\
 &\FQA &&\text{its full subcategory of $\A\in\QA$ such that $\A_1\in\FVect$,} \\
 &\QAsc &&\text{the category of quadratic algebras $\A\cong\A_0\otimes\B$ where $\B\in\QA$,} \\ 
 &\FQAsc &&\text{the category of quadratic algebras $\A\cong\A_0\otimes\B$ where $\B\in\FQA$.}
\end{align*}

The categories $\QA$ and $\QAsc$ consist of the connected and semi-connected quadratic algebras respectively. The subcategory $\FQA\subset\QA$ consists of finitely generated quadratic algebras $\A\in\QA$, while $\FQAsc\subset\QAsc$ consists of algebras $\A\in\QAsc$ finitely generated over $\A_0$. Note that a quadratic algebra $\A=T V/(R)\in\QA$, where $(R)$ is the ideal generated by a subspace $R\subset V^{\otimes2}$, belongs to $\FQA$ iff  $V\in\FVect$. We obtain the following picture:
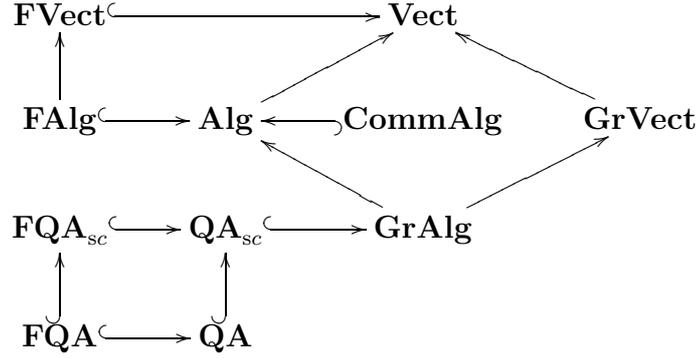
\begin{align*}
 \xymatrix{\FVect\ar@{^{(}->}[rr]      & &\Vect  \\
           \FAlg\ar@{^{(}->}[r]\ar[u]      &\Alg\ar@{->}[ur] & \CommAlg\ar@{^{(}->}[l]&\GrVect\ar[ul]   \\
            \FQAsc\ar@{^{(}->}[r] & \QAsc\ar@{^{(}->}[r]&\GrAlg\ar@{->}[ul]\ar[ur] \\
           \FQA\ar@{^{(}->}[r]\ar@{^{(}->}[u] & \QA\ar@{^{(}->}[u]  }
\end{align*}
where `$\hookrightarrow$' are inclusions of full subcategories (fully faithful functors) and `$\to$' are `forgetful' functors (they are faithful, but not full).

\subsection{Algebraic sets and affine schemes}

Let us introduce some notions from algebraic geometry~\cite{Perrin}, \cite{Fulton}, \cite{ManinAlgGeom}, \cite{Har}.

\bb{Algebraic sets.} Define {\it (affine) algebraic set} (over $\KK$) as the set of all the solutions of a system of algebraic equations with coefficients in $\KK$, that is a set of all $n$-tuples $(x^1,\ldots,x^n)\in\KK^n$ satisfying
\begin{align}
 F_\alpha(x^1,\ldots,x^n)=0, \label{Fixx}
\end{align}
where $F_\alpha\colon\KK^n\to\KK$ are polynomial functions of $n$ variables (with coefficients in $\KK$) and $\alpha$ runs over some set of indices. If the system~\eqref{Fixx} is empty then the corresponding algebraic set consists of all the $n$-tuples; it is denoted by $\AA^n=\AA^n_\KK$ and called {\it affine space}. As a set the affine space $\AA^n$ coincides with the vector space $\KK^n$, but we consider it as an object of another category, where there are more morphisms. Namely, a morphism of affine spaces $\Phi\colon\AA^n\to\AA^m$ is a map $(x^1,\ldots,x^n)\mapsto(P_1(x^1,\ldots,x^n),\ldots,P_m\big(x^1,\ldots,x^n)\big)$, where $P_i$ are polynomial functions.

\bb{Category of algebraic sets.} Any algebraic set is a subset of $\AA^n$ for a certain~$n$. A morphism between algebraic sets $X\subset\AA^n$ and $Y\subset\AA^m$ is a map $\varphi\colon X\to Y$ induced by a morphism $\Phi\colon\AA^n\to\AA^m$ in the sense that the diagram
\begin{align} \label{phiXY}
\xymatrix{X\ar[r]^{\varphi}\ar@{^{(}->}[d] & Y\ar@{^{(}->}[d] \\
\AA^n\ar[r]^{\Phi} & \AA^m}
\end{align}
commutes. Denote the category of algebraic sets over $\KK$ with these morphisms by $\AlgSet$.

\bb{Algebra of regular functions.} A {\it regular function} on an algebraic set $X$ is a morphism $X\to\AA^1$. Since $\AA^1=\KK$, the regular functions on $X$ form an algebra, denote it by $A(X)$. The algebra of regular functions on $\AA^n$ is the algebra of polynomials $\KK[x^1,\ldots,x^n]$. By setting $Y=\AA^m=\AA^1=\KK$ in the diagram~\eqref{phiXY} we see that a function $f\colon X\to\KK$ is regular iff there exists a polynomial $F\in\KK[x^1,\ldots,x^n]$ such that $f(x^1,\ldots,x^n)=F(x^1,\ldots,x^n)$ for any $(x^1,\ldots,x^n)\in X$. We obtain an algebra epimorphism $\KK[x^1,\ldots,x^n]\twoheadrightarrow A(X)$, so that $A(X)\cong\KK[x^1,\ldots,x^n]/I(X)$, where $I(X)\subset\KK[x^1,\ldots,x^n]$ is the ideal of polynomials $F\in\KK[x^1,\ldots,x^n]$ such that $F(x^1,\ldots,x^n)=0$ for any $(x^1,\ldots,x^n)\in X$.

\bb{Affine schemes.} For any morphism $\varphi\colon X\to Y$ its pull-back $\varphi^*\colon A(Y)\to A(X)$, $\varphi^*(f)=f\cdot\varphi$, is an algebra homomorphism. It is known that by mapping $X\mapsto A(X)$, $\varphi\mapsto\varphi^*$ we obtain a contravariant fully faithful functor $\AlgSet\to\CommAlg$. This means that the category $\AlgSet$ is embedded into $\CommAlg^{\op}$. The notion of algebraic set over~$\KK$ is generalised to the notion of {\it affine scheme over $\KK$} such that these schemes form the category $\AffSch=\CommAlg^{op}$. The affine scheme corresponding to an algebra $\gR\in\CommAlg$ is denoted by $\Spec\gR$. (see e.g.~\cite{ManinAlgGeom} for details). 


\bb{Sheaves of modules and vector bundles.} Any quasi-coherent sheaf of modules on an affine scheme $\Spec\gR$ is uniquely determined by the $\gR$-module of its global sections. Moreover, the category of all quasi-coherent sheaves of modules is equivalent to the category of $\gR$-modules. A vector bundle on an affine scheme is a particular case of a (quasi-)coherent sheaf. A trivial vector bundle (free sheaf) on $\Spec\gR$ with a fibre $V\in\FVect$ corresponds to the $\gR$-module $\gR\otimes V$. If $V=\KK^n$ and $\gR=A(X)$ for $X\in\AlgSet$, then a global section is a morphism $X\to X\times\AA^n$ such that $x\mapsto\big(x,v(x)\big)$ for a morphism $v\colon X\to\AA^n$. These sections form an $A(X)$-module isomorphic to $A(X)\otimes\KK^n$, an element $a\otimes v\in A(X)\otimes\KK^n$ gives a section $x\mapsto\big(x,a(x)v\big)$.

\bb{Quantum affine schemes.} It is natural to define the quantum or `non-commutative' version of the category of affine schemes as $\Alg^\op$ (see e.g.~\cite{Dr}). Thus, quantum affine schemes are algebras $\gR\in\Alg$, their morphisms are reversed algebra homomorphisms. As an analogue of the category of quasi-coherent sheaves on $\gR\in\Alg^\op$ one should regard the category of left, right or two-sided $\gR$-modules.

\subsection{Monoidal categories and functors}

Here we remind the basic concepts from the theory of monoidal categories~\cite{Mcl}, \cite{Bor2}.

\bb{Monoidal categories.} A {\it monoidal category} is a category $\bfC$ equipped with a bifunctor $\otimes\colon\bfC\times\bfC\to\bfC$ (called {\it a monoidal product}), {\it a unit object} $I\in\bfC$ and natural isomorphisms $(X\otimes Y)\otimes Z\cong X\otimes(Y\otimes Z)$, $I\otimes X\cong X$ and $X\otimes I\cong X$ satisfying some conditions (see~\cite{Mcl} for details). It is called {\it strict monoidal} if these isomorphisms are identities. Due to~\cite[11.3,Th.~1]{Mcl} any monoidal category is monoidally equivalent to a strict monoidal category, therefore we can suppose without loss of generality that all considered monoidal categories are strict monoidal (if it is not, one can always insert necessary isomorphisms into diagrams and formulae). We write $\bfC=(\bfC,\otimes,I)$ or simply $\bfC=(\bfC,\otimes)$, since the unit object is unique up to an isomorphism. In the second case we sometimes write $I_\bfC$ for the unit object. Note that if $\bfC=(\bfC,\otimes)$ is a monoidal category, then its opposite $\bfC^\op=(\bfC^\op,\otimes)$ is also a monoidal category with the same unit object $I_{\bfC^\op}=I_\bfC$. A subcategory $\bfC'$ of a monoidal category $(\bfC,\otimes,I)$ is called {\it monoidal subcategory} if $I\in\bfC'$ and $f\otimes g\in\bfC'$ for any morphisms $f,g$ in $\bfC'$ (in particular, $X\otimes Y\in\bfC'$ $\;\forall\,X,Y\in\bfC'$).

\bb{Examples of monoidal categories.} \label{bbExMonCat}
\begin{itemize}
 \item If there exists a terminal object and a product of any pair of objects in a category $\bfC$, then $(\bfC,\times)$ is a monoidal category and its terminal object is its unit object (in this case $\bfC$ is called {\it a category with finite products}, see~p.~\ref{bbCatFinPr}). In particular, $(\Vect,\oplus)$ and its subcategory $(\FVect,\oplus)$ are monoidal categories with unit object $0$, where $V\oplus W$ is the direct sum of the vector spaces $V$ and $W$.
 \item If there exists an initial object and a coproduct of any pair of objects in a category $\bfC$, then $(\bfC,\amalg)$ is a monoidal category. For instance, the category $(\CommAlg,\otimes)$ is monoidal, since the tensor product $\gR\otimes\gS$ of commutative algebras $\gR$ and $\gS$ is their coproduct. The opposite monoidal category is $(\AffSch,\times)$, where $\Spec\gR\times\Spec\gS=\Spec(\gR\otimes\gS)$ is the categorical product.
 \item Also, the category $\Vect$ and its subcategory $\FVect$ are monoidal with respect to the usual tensor product $\otimes$. The unit object is the one-dimensional vector space~$\KK=\KK^1$.
 \item The tensor product $\otimes$ gives a structure of monoidal category on $\GrVect$. The monoidal product of two graded vector spaces $V$ and $W$ is the usual tensor product $V\otimes W$ with the components $(V\otimes W)_k=\bigoplus\limits_{l=0}^k V_l\otimes W_{k-l}$. The unit object is the vector space $\KK$ with vanishing components of non-zero order: $(\KK)_k=\delta_{k0}\KK$. In this way we obtain the monoidal category $(\GrVect,\otimes,\KK)$.
 \item The tensor product of two algebras $\A,\B\in\Alg$ is the vector space $\A\otimes\B$ with the multiplication $(a\otimes b)(a'\otimes b')=(aa')\otimes(bb')$. This gives the monoidal category $(\Alg,\otimes,\KK)$.
 \item In the same way we obtain the monoidal category $(\GrAlg,\otimes,\KK)$. One can check that all the full subcategories of $\GrAlg$ defined above are monoidal subcategories.
 \item Another useful example of a monoidal product on $\GrVect$ is {\it Manin product}. For two objects $V,W\in\GrVect$ it is defined as the graded vector space $V\circ W$ with the components
\begin{align}
 (V\circ W)_k=V_k\otimes W_k. \label{MWP}
\end{align}
The monoidal category $(\GrVect,\circ)$ is equivalent to the direct product $\prod\limits_{k\in\NN_0}(\Vect,\otimes)$. Its unit object is a graded vector space with one-dimensional components, that is the polynomial algebra $\KK[u]$ with the standard grading. Note also that we have the (non-graded) inclusion of vector spaces $V\circ W\subset V\otimes W$.
 \item One can show that for any $\A,\B\in\GrAlg$ their Manin product $\A\circ\B$ is a subalgebra of $\A\otimes\B$. Thus we obtain the monoidal category $(\GrAlg,\circ,\KK[u])$ and its monoidal subcategories (the product `$\circ$' was originally defined by Manin for the category $\FQA$).
\end{itemize}

\bbo{Monoidal functors.} Let $\bfC=(\bfC,\otimes)$ and $\bfD=(\bfD,\odot)$ be two monoidal categories. The functor $F\colon\bfC\to\bfD$ equipped with a morphism $\varphi\colon I_\bfD\to FI_\bfC$ and a natural transformation $\phi_{X,Y}\colon FX\odot FY\to F(X\otimes Y)$ is called {\it lax monoidal functor} if for any $X,Y,Z\in\bfC$ the diagrams
\begin{gather}
\xymatrix{FX\odot FY\odot FZ\ar[rr]^{\id\odot\phi_{Y,Z}}\ar[d]^{\phi_{X,Y}\odot\id}& & FX\odot F(Y\otimes Z)\ar[d]^{\phi_{X,Y\otimes Z}} \\
F(X\otimes Y)\odot FZ\ar[rr]^{\phi_{X\otimes Y,Z}}& & F(X\otimes Y\otimes Z)
} \label{laxFd1} \\
\xymatrix{I_\bfD\odot FX\ar@{=}[r]\ar[d]^{\varphi\odot\id}&  FX\ar@{=}[d] \\
F(I_\bfC)\odot FX\ar[r]^{\phi_{I_\bfC,X}}&  F(I_\bfC\otimes X)
} \qquad\qquad
\xymatrix{FX\odot I_\bfD\ar@{=}[r]\ar[d]^{\id\odot\varphi}&  FX\ar@{=}[d] \\
FX\odot F(I_\bfC)\ar[r]^{\phi_{X,I_\bfC}}&  F(X\otimes I_\bfC)
} \label{laxFd2}
\end{gather}
commute. The functor $F\colon\bfC\to\bfD$ is called {\it colax monoidal} if its opposite $F^\op\colon\bfC^\op\to\bfD^\op$ is lax monoidal. In other words, this is a triple $(F,\varphi,\phi)$, where $\varphi\colon FI_\bfC\to I_\bfD$ and $\phi$ is a natural transformation with components $\phi_{X,Y}\colon F(X\otimes Y)\to FX\odot FY$ such that the reversed diagrams~\eqref{laxFd1}, \eqref{laxFd2} commute, i.e. $(\phi_{X,Y}\odot\id_{FZ})\cdot\phi_{X\otimes Y,Z}=(\id_{FX}\odot\phi_{Y,Z})\cdot\phi_{X,Y\otimes Z}$, $(\varphi\odot\id_{FX})\cdot\phi_{I_\bfC,X}=\id_{FX}=(\id_{FX}\odot\varphi)\cdot\phi_{X,I_\bfC}$.

A {\it strong monoidal functor} is a lax (or colax) monoidal functor $(F,\varphi,\phi)$ such that $\varphi$ and $\phi_{X,Y}$ are isomorphisms. A strong monoidal functor is called {\it strict monoidal functor} if they are identities. For instance, the inclusion of a monoidal subcategory $\bfC'\subset\bfC$ into $\bfC$ is a faithful strict monoidal functor. Note that a composition of two lax/colax/strong monoidal functors is a lax/colax/strong monoidal functor.

\bb{Examples of monoidal functors.} \label{bbExMonFunc}
\begin{itemize}
 \item Consider the dualisation functor $(-)^*\colon\Vect^\op\to\Vect$ which sends a vector space $V\in\Vect$ to its dual $V^*$, it is defined on linear maps $f\colon V\to W$ by the formula~\eqref{fh}. For any vector spaces $V,W\in\Vect$ the tensor product $V^*\otimes W^*$ is naturally embedded into $(V\otimes W)^*$, so the functor $(-)^*\colon\Vect^\op\to\Vect$ has a structure of lax monoidal functor $(\Vect^\op,\otimes)\to(\Vect,\otimes)$. The opposite functor $(\Vect,\otimes)\to(\Vect^\op,\otimes)$ is colax monoidal.
 \item The same functor considered on finite-dimensional vector spaces is a strong monoidal functor $(-)^*\colon(\FVect,\otimes)\to(\FVect^\op,\otimes)$, since the embedding gives the natural isomorphisms $V^*\otimes W^*\cong (V\otimes W)^*$ for any $V,W\in\FVect$. 
 \item Note that $(V\oplus W)^*\cong V^*\oplus W^*$ for any $V,W\in\Vect$, so we obtain strong monoidal contravariant functors $(-)^*\colon(\Vect,\oplus)\to(\Vect,\oplus)$ and $(-)^*\colon(\FVect,\oplus)\to(\FVect,\oplus)$.
 \item The functor $T\colon\Vect\to\GrAlg$, which gives the tensor algebra $T V$ of a vector space $V\in\Vect$, is a strong monoidal functor $T\colon(\Vect,\otimes)\to(\GrAlg,\circ)$. The isomorphisms $\phi_{V,W}\colon T V\circ T W\isoright T(V\otimes W)$ have the graded components
\begin{align}
 &(\phi_{V,W})_k\colon V^{\otimes k}\otimes W^{\otimes k}\isoright(V\otimes W)^{\otimes k}, \label{phiVWT} \\
 &v_1\otimes\cdots\otimes v_k\otimes w_1\otimes\cdots\otimes w_k\mapsto(v_1\otimes w_1)\otimes\cdots\otimes(v_k\otimes w_k). \notag
\end{align}
It gives also the strong monoidal functors
\begin{align}
 &T\colon(\Vect,\otimes)\to(\QA,\circ), &
 &T\colon(\FVect,\otimes)\to(\FQA,\circ).
\end{align}
\end{itemize}

\bbo{Symmetric monoidal category.} A monoidal category $\bfC=(\bfC,\otimes)$ equipped with a natural isomorphism $\sigma_{X,Y}\colon X\otimes Y\isoright Y\otimes X$ is called {\it symmetric} if for any $X,Y,Z\in\bfC$ we have $\sigma_{Y,X}\sigma_{X,Y}=\id_{X\otimes Y}$ and the diagrams
\begin{gather} \label{sigmaDiag}
\xymatrix{X\otimes Y\otimes Z\ar[rr]^{\sigma_{X,Y}\otimes\id}\ar[drr]_{\sigma_{X,Y\otimes Z}} && Y\otimes X\otimes Z\ar[d]^{\id\otimes\sigma_{X,Z}} \\
   && Y\otimes Z\otimes X
} \qquad\qquad
\xymatrix{X\otimes I_\bfC\ar@{=}[d]\ar[dr]^{\sigma_{X,I_\bfC}}  \\
X\ar@{=}[r] & I_\bfC\otimes X
}
\end{gather}
commute. The natural transformation $\sigma$ gives an additional structure on the monoidal category $\bfC$.

All the monoidal categories considered in p.~\ref{bbExMonCat} are symmetric. For example, the symmetric structure of $(\Vect,\otimes)$ is given by the isomorphisms
\begin{align} \label{sigmaVect}
 &\sigma_{V,W}\colon V\otimes W\isoright W\otimes V, &&\sigma_{V,W}(v\otimes w)=w\otimes v.
\end{align}
One can check that it induces the symmetric structure of the monoidal categories $(\Alg,\otimes)$, $(\GrVect,\otimes)$, $(\GrAlg,\otimes)$, $(\GrVect,\circ)$, $(\GrAlg,\circ)$ and of their monoidal subcategories considered above. In fact, we already used the maps~\eqref{sigmaVect}: the natural isomorphisms~\eqref{phiVWT} are compositions of some operators of the form $\id^{\otimes l}\otimes\sigma_{V,W}\otimes\id^{\otimes(2k-l-2)}$, $l=1,\ldots,2k-3$, $k\ge2$.

If a monoidal category $\bfC=(\bfC,\otimes)$ is symmetric, then $\bfC^\op=(\bfC^\op,\otimes)$ is also symmetric: the role of the structure isomorphisms $V\otimes W\isoright W\otimes V$ in $\bfC^\op$ is played by the structure isomorphisms $\sigma_{W,V}\colon W\otimes V\isoright V\otimes W$ of the category $\bfC$.

\bb{Symmetric monoidal functors.} Let $(\bfC,\otimes)$ and $(\bfD,\odot)$ be symmetric monoidal categories. A lax (or strong) monoidal functor $F\colon(\bfC,\otimes)\to(\bfD,\odot)$ with the structure natural transformation $\phi_{X,Y}\colon FX\odot FY\to F(X\otimes Y)$ is called {\it symmetric} if it preserves the symmetric structure in the sense that all the diagrams
\begin{gather} \label{diagSymMonFunc}
\xymatrix{FX\odot FY\ar[rr]^{\sigma_{FX,FY}}\ar[d]^{\phi_{X,Y}} && FY\odot FX\ar[d]^{\phi_{Y,X}} \\
 F(X\otimes Y)\ar[rr]^{F\sigma_{X,Y}} && F(Y\otimes X)
}
\end{gather}
commute. A colax monoidal functor $F\colon(\bfC,\otimes)\to(\bfD,\odot)$ with the structure natural transformation $\phi_{X,Y}\colon F(X\otimes Y)\to FX\odot FY$ is called {\it symmetric} if the diagrams~\eqref{diagSymMonFunc} with the reversed vertical arrows commute (i.e. if the lax monoidal functor $F^\op\colon(\bfC^\op,\otimes)\to(\bfD^\op,\odot)$ is symmetric). A composition of symmetric monoidal functors is symmetric. All the monoidal functors considered in p.~\ref{bbExMonFunc} are symmetric.

\subsection{Monoids and their actions}

We define a notion of monoid, bimonoid and Hopf monoid for a general (symmetric) monoidal category~\cite{Mcl}, \cite{Porst}. They generalise the notions of algebra, bialgebra and Hopf algebra (for this case see~\cite{Kass}).

\bb{Monoids and comonoids.} A {\it monoid} in a monoidal category $(\bfC,\otimes)$ is a triple $\MM=(X,\mu_X,\eta_X)$ of an object $X\in\bfC$ with structure morphisms $\mu_X\colon X\otimes X\to X$ and $\eta_X\colon I_\bfC\to X$ such that the diagrams
\begin{gather} \label{diagMonoid}
\xymatrix{X\otimes X\otimes X\ar[rr]^{\id\otimes\mu_X}\ar[d]^{\mu_X\otimes\id} && X\otimes X\ar[d]^{\mu_X} \\
 X\otimes X\ar[rr]^{\mu_X} && X
}  \qquad\qquad
\xymatrix{I\otimes X\ar[r]^{\eta_X\otimes\id}\ar@{=}[dr] & X\otimes X\ar[d]^{\mu_X} & X\otimes I\ar[l]_{\id\otimes\eta_X}\ar@{=}[dl] \\
 & X
} 
\end{gather}
commute. Morphisms $\mu_X$ and $\eta_X$, which give a structure of monoid on $X$, are often called {\it multiplication} and {\it unit}.

A {\it comonoid} in $(\bfC,\otimes)$ is a triple $\OO=(X,\Delta_X,\varepsilon_X)$ of an object $X\in\bfC$ and morphisms $\Delta_X\colon X\to X\otimes X$ and $\varepsilon_X\colon X\to I_\bfC$ in $\bfC$ such that the reversed diagrams~\eqref{diagMonoid} commute, i.e. $(\id_X\otimes\Delta_X)\cdot\Delta_X=(\Delta_X\otimes\id_X)\cdot\Delta_X$ and $(\varepsilon_X\otimes\id_X)\cdot\Delta_X=\id_X=(\id_X\otimes\varepsilon_X)\cdot\Delta_X$. For $\Delta_X$ and $\varepsilon_X$ one uses the terms `{\it comultiplication}' and `{\it counit}'.

\bb{Morphisms of monoids and comonoids.} \label{bbMorMonComon}
 A~morphism of monoids $(X,\mu_X,\eta_X)$ and $(Y,\mu_Y,\eta_Y)$ or comonoids $(X,\Delta_X,\varepsilon_X)$ and $(Y,\Delta_Y,\varepsilon_Y)$ is a morphism $f\colon X\to Y$ in $\bfC$ preserving the structure morphisms:
\begin{gather} \label{diagMonoidMor}
\xymatrix{X\otimes X\ar[r]^{\mu_X}\ar[d]_{f\otimes f} & X\ar[d]^{f} \\
 Y\otimes Y\ar[r]^{\mu_Y} & Y
}  \qquad
\xymatrix{I_\bfC\ar[r]^{\eta_X}\ar[rd]_{\eta_Y} & X\ar[d]^{f} \\
 & Y
}  \qquad\text{or}\qquad
\xymatrix{X\ar[r]^{\Delta_X}\ar[d]_{f} & X\otimes X\ar[d]^{f\otimes f} \\
 Y\ar[r]^{\Delta_Y} & Y\otimes Y
}  \qquad
\xymatrix{X\ar[d]_{f}\ar[r]^{\varepsilon_X} & I_\bfC \\
  Y\ar[ru]_{\varepsilon_Y}
}  \qquad
\end{gather}
Since the composition of two morphisms of (co)monoids is a morphisms of (co)monoids, the monoids and comonoids in a monoidal category $(\bfC,\otimes)$ form categories denoted by $\Mon(\bfC,\otimes)$ and $\Comon(\bfC,\otimes)$ respectively. Note that $\Comon(\bfC,\otimes)=\big(\Mon(\bfC^\op,\otimes)\big)^\op$. 

\begin{Prop} \label{PropMonIso}
 Consider two monoids $(X,\mu_X,\eta_X),(Y,\mu_Y,\eta_Y)\in\Mon(\bfC)$ or comonoids $(X,\Delta_X,\varepsilon_X),(Y,\Delta_Y,\varepsilon_Y)\in\Comon(\bfC)$ in a monoidal category $\bfC=(\bfC,\otimes)$.  Let $f\colon X\isoright Y$ be an isomorphism in the category~$\bfC$. If $f\colon X\to Y$ preserves structure of these monoids or comonoids, then $f^{-1}\colon Y\to X$ also preserves structure of these monoids or comonoids. In particular, this implies an isomorphism of the corresponding monoids or comonoids in $\Mon(\bfC)$ or $\Comon(\bfC)$ respectively.
\end{Prop}

\noindent{\bf Proof.} It is enough to consider the case of monoids. By using the commutativity of the left diagrams~\eqref{diagMonoidMor} one yields $f^{-1}\cdot\eta_Y=f^{-1}\cdot f\cdot\eta_X=\eta_X$ and $f^{-1}\cdot\mu_Y=f^{-1}\cdot\mu_Y\cdot(f\otimes f)\cdot(f^{-1}\otimes f^{-1})=f^{-1}\cdot f\cdot\mu_X\cdot(f^{-1}\otimes f^{-1})=\mu_X\cdot(f^{-1}\otimes f^{-1})$. \qed

\bb{Monoidal products of monoids and comonoids.} \label{bbTPMon}
If $\bfC$ is symmetric then $\Mon(\bfC)$ and $\Comon(\bfC)$ are monoidal categories and they are also symmetric: the monoidal product of monoids  $(X,\mu_X,\eta_X)$ and $(Y,\mu_Y,\eta_Y)$ or comonoids $(X,\Delta_X,\varepsilon_X)$ and $(Y,\Delta_Y,\varepsilon_Y)$ is the object $X\otimes Y$ with the structure defined by the diagrams
\begin{align*}
\xymatrix{X\otimes Y\otimes X\otimes Y\ar[d]_{\id\otimes\sigma_{Y,X}\otimes\id}\ar[rrd]^{\mu_{X\otimes Y}}   \\
 X\otimes X\otimes Y\otimes Y\ar[rr]^{\mu_X\otimes\mu_Y} && X\otimes Y
}  \qquad
\xymatrix{I_\bfC\ar@{=}[d]\ar[rrd]^{\eta_{X\otimes Y}} \\
 I_\bfC\otimes I_\bfC\ar[rr]^{\eta_X\otimes\eta_Y}&& X\otimes Y
}  \qquad\text{or} \\
\xymatrix{X\otimes Y\ar[d]_{\Delta_X\otimes\Delta_Y}\ar[rrd]^{\Delta_{X\otimes Y}}   \\
 X\otimes X\otimes Y\otimes Y\ar[rr]^{\id\otimes\sigma_{X,Y}\otimes\id} && X\otimes Y\otimes X\otimes Y
}  \qquad
\xymatrix{X\otimes Y\ar[d]_{\varepsilon_X\otimes\varepsilon_Y}\ar[rrd]^{\varepsilon_{X\otimes Y}} \\
 I_\bfC\otimes I_\bfC\ar@{=}[rr]&& I_\bfC
}
\end{align*}

A usual monoid is a monoid in the monoidal category $\Set=(\Set,\times)$. A monoid in $(\Vect,\otimes)$ is an algebra: $\Mon(\Vect,\otimes)=(\Alg,\otimes)$. In its subcategory of finite-dimensional vector spaces we have $\Mon(\FVect,\otimes)=(\FAlg,\otimes)$. Comonoids in $(\Vect,\otimes)$ are called {\it coalgebras}. Analogously, we have $\Mon(\GrVect,\otimes)=(\GrAlg,\otimes)$.

\bb{Eckmann--Hilton Principle.} \label{bbEHP}
 A monoid $(X,\mu_X,\eta_X)$ in a symmetric monoidal category $(\bfC,\otimes)$ is called {\it commutative} if $\mu_X\cdot\sigma_{X,X}=\mu_X$. Analogously, a comonoid $(X,\Delta_X,\varepsilon_X)$ is called {\it cocommutative} if $\sigma_{X,X}\cdot\Delta_X=\Delta_X$. Thus we obtain full monoidal subcategories $\cMon(\bfC,\otimes)\subset\Mon(\bfC,\otimes)$ and $\cocComon(\bfC,\otimes)\subset\Mon(\bfC,\otimes)$. For example, for the category $\bfC=\Vect$ with the tensor product $\otimes$ we have $\cMon(\Vect,\otimes)=(\CommAlg,\otimes)$.

The Eckmann--Hilton Principle was formulated in~\cite[Th~1.12]{EH} for the groups: if a group homomorphism $\mu\colon G\times G\to G$ has a neutral element as a binary operation, then the group $G$ is abelian and $\mu$ coincides with the group multiplication. This fact is generalised for any symmetric monoidal category $(\bfC,\otimes)$ (see e.g.~\cite[\S~4]{Porst}).

\begin{Prop} \label{PropMonMon}
 The category $\Mon\big(\Mon(\bfC,\otimes)\big)$ is equivalent to $\cMon(\bfC,\otimes)$. Namely, the monoid structure $(\mu_\MM,\eta_\MM)$ on an object $\MM=(X,\mu_X,\eta_X)\in\Mon(\bfC,\otimes)$ exists iff the monoid $\MM$ is commutative. This structure is unique: $\mu_\MM=\mu_X$, $\eta_\MM=\eta_X$.
 The category $\Comon\big(\Comon(\bfC,\otimes)\big)$ is equivalent to $\cocComon(\bfC,\otimes)$ in the same way.
\end{Prop}

\bbo{The functors $\Mon(F)$ and $\Comon(F)$.} \label{bbMonF}
 Any lax monoidal functor $F\colon\bfC\to\bfD$ between monoidal categories $\bfC=(\bfC,\otimes)$ and $\bfD=(\bfD,\odot)$ induces the functor
\begin{align} \label{MonF}
 \Mon(F)\colon\Mon(\bfC)\to\Mon(\bfD).
\end{align}
Namely, if $\MM=(X,\mu_X,\eta_X)$ is a monoid in $\bfC$, then $\Mon(F)\MM=(FX,\mu_{FX},\eta_{FX})$ is a monoid in $\bfD$, where $\mu_{FX}\colon FX\odot FX\to FX$ and $\eta_{FX}\colon I_\bfD\to FX$ are the compositions
\begin{align} \label{muFX}
 &FX\odot FX\xrightarrow{\phi_{X,X}}F(X\otimes X)\xrightarrow{F\mu_X} FX &&\text{and} && I_\bfD\xrightarrow{\varphi}FI_\bfC\xrightarrow{F\eta_X} FX
\end{align}
respectively. Dually, a colax monoidal functor $F\colon\bfC\to\bfD$ induces the functor
\begin{align} \label{ComonF}
 \Comon(F)\colon\Comon(\bfC)\to\Comon(\bfD).
\end{align}
Note that the functors~\eqref{MonF} and \eqref{ComonF} depend on the monoidal structure on $F$, that is $\Mon(F)=\Mon(F,\varphi,\phi)$ and $\Comon(F)=\Comon(F,\varphi,\phi)$.

A faithful lax or colax monoidal functor $F$ induces a faithful functor~\eqref{MonF} or \eqref{ComonF} respectively. A strong monoidal functor $F$ gives strong monoidal functors~\eqref{MonF} and \eqref{ComonF}. If $F$ is strong monoidal and fully faithful, then~\eqref{MonF} and \eqref{ComonF} are also fully faithful (the condition to be strong monoidal can be weakened: it is enough to require $\varphi$ and all $\phi_{X,X}$ to be epimorphisms/monomorphisms for lax/colax case).

Suppose the categories $\bfC$ and $\bfD$ are symmetric. If $F$ is a symmetric lax/strong monoidal functor, then $\Mon(F)$ is a symmetric lax/strong monoidal functor as well. If $F$ is a symmetric colax/strong monoidal functor, then $\Comon(F)$ is also a symmetric colax/strong monoidal functor.

\bb{Bimonoids.} \label{bbBimon}
Let $\bfC=(\bfC,\otimes)$ be a symmetric monoidal category. A {\it bimonoid} in $\bfC$ is $\BB=(X,\mu_X,\eta_X,\Delta_X,\varepsilon_X)$ such that $\MM=(X,\mu_X,\eta_X)$ and $\OO=(X,\Delta_X,\varepsilon_X)$ are monoid and comonoid in $\bfC$ with compatible structures:
\begin{gather}
\begin{split} \label{BimodDiag}
\xymatrix{ X\otimes X\ar[r]^{\mu_X}\ar[d]_{\Delta_X\otimes\Delta_X}& X\ar[r]^{\Delta_X} & X\otimes X \\
X\otimes X\otimes X\otimes X\ar[rr]^{\id\otimes\sigma_{X,X}\otimes\id} && X\otimes X\otimes X\otimes X\ar[u]_{\mu_X\otimes\mu_X}
} \qquad
\xymatrix{ I_\bfC\ar@{=}[rd]\ar[r]^{\eta_X}& X\ar[d]^{\varepsilon_X} \\
 & I_\bfC
} \\
\xymatrix{ I_\bfC\ar[rr]^{\eta_X}\ar@{=}[d]&& X\ar[d]^{\Delta_X} \\
 I_\bfC\otimes I_\bfC\ar[rr]^{\eta_X\otimes\eta_X}&& X\otimes X
} \qquad\qquad
\xymatrix{ X\otimes X\ar[rr]^{\mu_X}\ar[d]_{\varepsilon_X\otimes\varepsilon_X}&& X\ar[d]^{\varepsilon_X} \\
 I_\bfC\otimes I_\bfC\ar@{=}[rr]&& I_\bfC
}
\end{split}
\end{gather}
These diagrams mean exactly that $\Delta_X\colon X\to X\otimes X$ and $\varepsilon_X\colon X\to I_\bfC$ are monoid morphisms $(X\otimes X,\mu_{X\otimes X},\eta_{X\otimes X})\to(X,\mu_X,\eta_X)$ and $(X,\mu_X,\eta_X)\to(I_\bfC,\id_\bfC,\id_\bfC)$ or, equivalently, that $\mu_X\colon X\otimes X\to X$ and $\varepsilon_X\colon I_\bfC\to X$ are comonoid morphisms $(X\otimes X,\Delta_{X\otimes X},\varepsilon_{X\otimes X})\to(X,\Delta_X,\varepsilon_X)$ and $(I_\bfC,\id_\bfC,\id_\bfC)\to(X,\Delta_X,\varepsilon_X)$. Bimonoids in $\bfC=(\bfC,\otimes)$ are objects of the symmetric monoidal category $\Bimon(\bfC,\otimes)=\Bimon(\bfC)=\Comon\big(\Mon(\bfC)\big)=\Mon\big(\Comon(\bfC)\big)$.

Any symmetric strong monoidal functor $F=(F,\varphi,\phi)\colon\bfC\to\bfD$ between symmetric monoidal categories $\bfC$ and $\bfD$ induces the symmetric strong monoidal functor
\begin{align} \label{BimonF}
 \Bimon(F)\colon\Bimon(\bfC)\to\Bimon(\bfD),
\end{align}
where $\Bimon(F)=\Mon\big(\Comon(F)\big)= \Comon\big(\Mon(F)\big)$.

\bb{Hopf monoids.} For arbitrary bimonoid $(X,\mu_X,\eta_X,\Delta_X,\varepsilon_X)$ define a convolution on the set $\End_\bfC(X)$. The convolution of two morphisms $\alpha\colon X\to X$ and $\beta\colon X\to X$ in $\bfC$ is the morphism $\alpha\star\beta:=\mu_X\cdot(\alpha\otimes\beta)\cdot\Delta_X\colon X\to X$. This is a multiplication on the set $\End_\bfC(X)$ with the neutral element $\eta_X\cdot\varepsilon_X$. A{\it~Hopf monoid} in $\bfC$ is a bimonoid $(X,\mu_X,\eta_X,\Delta_X,\varepsilon_X)$ in $\bfC$ such that the morphism $\id_X$ has an inverse $\zeta_X\colon X\to X$ with respect to the convolution, that is $\zeta_X\star\id_X=\id_X\star\zeta_X=\eta_X\cdot\varepsilon_X$. The morphism $\zeta_X\colon X\to X$ is called {\it antipode}. It is unique since the convolution is associative.
 
The functor~\eqref{BimonF} maps a Hopf monoid in $\bfC$ to a Hopf monoid in $\bfD$, the antipode $\zeta_X\colon X\to X$ is mapped to the antipode $F(\zeta_X)\colon FX\to FX$.

Bimonoids and Hopf monoids in $(\Vect,\otimes)$ are {\it bialgebras} and {\it Hopf algebras} respectively.

\bb{Opposite monoid and coopposite comonoid.} \label{bbOpMon}
Let $\MM=(X,\mu_X,\eta_X)$ be a monoid in a symmetric monoidal category $\bfC=(\bfC,\otimes)$. Let $\mu_X^\op:=\mu_X\cdot\sigma_{X,X}\colon X\otimes X\to X$ be {\it opposite multiplication}. Then $\MM^\op:=(X,\mu_X^\op,\eta_X)$ is also monoid in $\bfC$, it is called {\it monoid opposite to} $\MM$. Note that $\MM$ is commutative in the sense of p.~\ref{bbEHP} iff $\MM^\op=\MM$, but in general these are different monoids. If $f\colon X\to \wt X$ is a morphism of the monoids $\MM=(X,\mu_X,\eta_X)$ and $\wt\MM=(\wt X,\mu_{\wt X},\eta_{\wt X})$, then it is a morphism $\MM^\op\to\wt\MM^\op$ as well.

Dually, let $\OO=(X,\Delta_X,\varepsilon_X)$ be a comonoid in a monoidal category $\bfC=(\bfC,\otimes)$. The comonoid $\OO^\cop:=(X,\Delta_X^\cop,\varepsilon_X)\in\Comon(\bfC)$ with the so-called {\it coopposite comultiplication} $\Delta^\cop=\sigma_{X,X}\cdot\Delta$ and the same counit $\varepsilon_X$ is called {\it comonoid coopposite to} $\OO$.

For a bimonoid $\BB=(X,\mu_X,\eta_X,\Delta_X,\varepsilon_X)$ we can define the following tree bimonoids: {\it opposite bimonoid} $\BB^\op=(X,\mu_X^\op,\eta_X,\Delta_X,\varepsilon_X)$, {\it coopposite bimonoid} $\BB^\cop=(X,\mu_X,\eta_X,\Delta_X^\cop,\varepsilon_X)$ and {\it opposite coopposite bimonoid} $\BB^{\op,\cop}=(X,\mu_X^\op,\eta_X,\Delta_X^\cop,\varepsilon_X)$ (the commutativity of the corresponding diagrams~\eqref{BimodDiag} is check straightforwardly). Moreover, if the bimonoid $\BB$ is a Hopf monoid, then $\BB^{\op,\cop}$ is also a Hopf monoid with the same antipode.

\bb{Actions and coactions.} A {\it (left) action} of a monoid $(X,\mu_X,\eta_X)$ on an object $V\in\bfC$ is a morphism $a\colon X\otimes V\to V$ making the diagrams
\begin{align} \label{aDiag}
\xymatrix{ X\otimes X\otimes V\ar[d]_{\id_X\otimes a}\ar[rr]^{\mu_X\otimes\id_V} && X\otimes V\ar[d]^{a} \\
 X\otimes V\ar[rr]^{a} && V
} \qquad\qquad\qquad
\xymatrix{I_\bfC\otimes V\ar[rr]^{\eta_X\otimes\id_V}\ar@{=}[drr] && X\otimes V\ar[d]^{a} \\
 && V
}
\end{align}
commutative. A {\it (left) coaction} of a comonoid $(X,\Delta_X,\varepsilon_X)$ on an object $V\in\bfC$ is a morphism $\delta\colon V\to X\otimes V$ making the reversed diagrams commutative:
\begin{align} \label{deltaDiag}
\xymatrix{ V\ar[d]_{\delta}\ar[rr]^{\delta} && X\otimes V\ar[d]^{\Delta_X\otimes\id_V} \\
 X\otimes V\ar[rr]^{\id_X\otimes\delta} && X\otimes X\otimes V
} \qquad\qquad\qquad
\xymatrix{V\ar[r]^{\delta}\ar@{=}[dr] & X\otimes V\ar[d]^{\varepsilon_X\otimes\id_V} \\
 & I_\bfC\otimes V
}
\end{align}
The actions of a monoid $\MM=(X,\mu_X,\eta_X)\in\Mon(\bfC)$ or a comonoid $\OO=(X,\Delta_X,\varepsilon_X)\in\Comon(\bfC)$ form a category $\Lact(\MM)$ or $\Lcoact(\OO)$ respectively. Their objects are pairs $(V,a)$ and $(V,\delta)$. A morphism $(V,a)\to(V',a')$ in $\Lact(\MM)$ or $(V,\delta)\to(V',\delta')$ in $\Lcoact(\MM)$ is a morphism $f\colon V\to V'$ in $\bfC$ such that the diagram
\begin{align} \label{LactMor}
\xymatrix{X\otimes V\ar[d]_{\id_X\otimes f}\ar[r]^a & V\ar[d]^f \\
 X\otimes V'\ar[r]^{a'} & V'
} \qquad\text{or}\qquad
\xymatrix{V\ar[d]_{f}\ar[r]^\delta & X\otimes V\ar[d]^{\id_X\otimes f} \\
 V'\ar[r]^{\delta'} & X\otimes V'
}
\end{align}
commutes. For the case $\bfC=(\Vect,\otimes)$ the monoid $\MM$ is an algebra $\gR\in\Alg$ and the objects of $\Lact(\gR)$ are exactly $\gR$-modules, while $\Lcoact(\gC)$ is the category of comodules of a coalgebra $\gC\in\Comon(\Vect,\otimes)$.

\bb{Actions and coactions of a bimonoid.} \label{bbActBimon}
A bimonoid in a symmetric monoidal category $\bfC=(\bfC,\otimes)$ is $\BB=(X,\mu_X,\eta_X,\Delta_X,\varepsilon_X)\in\Bimon(\bfC)$. By considering it as the monoid $\MM=(X,\mu_X,\eta_X)$ we obtain the category $\Lact(\BB):=\Lact(\MM)$. The comonoid structure $(\Delta_X,\varepsilon_X)$ turns $\Lact(\BB)$ into a monoidal category in the following way~\cite{P}. Let $(V,a)$ and $(W,b)$ be two objects of $\Lact(\BB)$. Their monoidal product $(V,a)\otimes(W,b)$ is the object $V\otimes W$ with the action
\begin{align}
 X\otimes(V\otimes W)\xrightarrow{\Delta_X\otimes\id}X\otimes X\otimes V\otimes W \xrightarrow{\id\otimes\sigma_{X,V}\otimes\id}X\otimes V\otimes X\otimes W \xrightarrow{a\otimes b} V\otimes W.
\end{align}
If $f\colon V\to V'$ and $g\colon W\to W'$ are morphisms in $\Lact(\BB)$ of the forms $(V,a)\to (V',a')$ and $(W,b)\to(W',b')$, then their monoidal product $f\otimes g\colon V\otimes W\to V'\otimes W'$ is a morphism $(V,a)\otimes(W,b)\to (V',a')\otimes(W',b')$.
The unit object of $\Lact(\BB)$ is $I_\bfC$ with the action $X\otimes I_\bfC=X\xrightarrow{\varepsilon_X}I_\bfC$. If $\BB$ is cocommutative ($\BB^\cop=\BB$), then the symmetric structure on $\bfC$ gives a symmetric structure on the monoidal category $\Lact(\BB)$.

Dually, the structure $(\mu_X,\eta_X)$ turns the category $\Lcoact(\BB):=\Lcoact(X,\Delta_X,\varepsilon_X)$ into monoidal category. If $\delta$ and $\gamma$ are coactions of the comonoid $(X,\Delta_X,\varepsilon_X)$ on $V$ and $W$ respectively, then its coaction on $V\otimes W$ is the composition
\begin{align} \label{deltagamma}
 V\otimes W\xrightarrow{\delta\otimes\gamma}X\otimes V\otimes X\otimes W \xrightarrow{\id\otimes\sigma_{V,X}\otimes\id}X\otimes X\otimes V\otimes W\xrightarrow{\mu_X\otimes\id} X\otimes(V\otimes W).
\end{align}
The coaction on $I_\bfC$ is the morphism $I_\bfC\xrightarrow{\eta_X}X=X\otimes I_\bfC$. If $\BB$ is commutative ($\BB^\op=\BB$), then $\Lcoact(\BB)$ is a symmetric monoidal category.

\bb{Translation of actions and coactions.} \label{bbFact}
 Let $F\colon(\bfC,\otimes)\to(\bfD,\odot)$ be a lax monoidal functor with structure morphisms $\phi_{X,Y}\colon FX\odot FY\to F(X\otimes Y)$, $\varphi\colon I_\bfD\to FI_\bfC$. The functor $\Mon(F)\colon\Mon(\bfC)\to\Mon(\bfD)$ maps a monoid $\MM=(X,\mu_X,\eta_X)$ to the monoid $\wt\MM=\Mon(F)(\MM)=(FX,\mu_{FX},\eta_{FX})$, where $\mu_{FX}\colon FX\odot FX\to FX$, $\eta_{FX}\colon I_\bfD\to FX$ are~\eqref{muFX}. Let $a\colon X\otimes V\to V$ be an action of the monoid $\MM$ on $V$ and consider the morphism $\wt a\colon FX\odot FV\to FV$ defined as a composition
\begin{align}
 FX\odot FV\xrightarrow{\phi_{X,V}}F(X\otimes V)\xrightarrow{F(a)}FV.
\end{align}
It is straightforward to check that $\wt a$ is an action of $\wt\MM$ on $FV$. If $f\colon V\to V'$ is a morphism $(V,a)\to(V',a')$ in $\Lact(\MM)$ then $F(f)\colon FV\to FV'$ is a morphism between the corresponding actions in $\Lact(\wt\MM)$. Thus we obtain a functor
\begin{align} \label{FLact}
 \Lact(\MM)\to\Lact(\wt\MM) 
\end{align}
induced by the lax monoidal functor $F=(F,\varphi,\phi)$.

Dually, a colax monoidal functor $F=(F,\varphi,\phi)\colon\bfC\to\bfD$ induces the functor
\begin{align} \label{FLcoact}
 &\Lcoact(\OO)\to\Lcoact(\wt\OO)
\end{align}
for comonoids $\OO\in\Comon(\bfC)$ and $\wt\OO=\Comon(F)\OO\in\Comon(\bfD)$.

If the monoidal categories $\bfC=(\bfC,\otimes)$ and $\bfD=(\bfD,\odot)$ are symmetric, then any symmetric strong monoidal functor $F=(F,\varphi,\phi)\colon\bfC\to\bfD$ induces the strong monoidal functors
\begin{align} \label{FLactB}
 &\Lact(\BB)\to\Lact(\wt\BB), &
 &\Lcoact(\BB)\to\Lcoact(\wt\BB)
\end{align}
for bimonoids $\BB\in\Bimon(\bfC)$ and $\wt\BB=\Bimon(F)\BB\in\Bimon(\bfD)$.

\subsection{Monoids and groups in a category with finite products}

A category with finite products is an important case of monoidal category. Here we consider this case in details (see~\cite[\S~3.5, 3.6]{Mcl}).

\bb{Category with finite products.} \label{bbCatFinPr}
Let $\bfC$ be a category. Suppose that there exists a terminal object $E$ and a categorical product $X\times Y$ of any two objects $X,Y\in\bfC$. Then there exist all the {\it finite products} $\prod\limits_{i=1}^nX_i$, where $X_i\in\bfC$, $n\in\NN_0$ (for $n=0$ this product is equal to~$E$). In this case $(\bfC,\times)$ is a symmetric monoidal category~\cite[\S~7.7]{Mcl}. For morphisms $f\colon Y\to X_1$ and $g\colon Y\to X_2$ we denote by $(f,g)$ the unique morphism $Y\to X_1\times X_2$ making the diagram
\begin{align} \label{fgp1p2}
\xymatrix{
 & Y\ar[dl]_f\ar[d]^{(f,g)}\ar[dr]^g & \\
 X_1 & X_1\times X_2\ar[l]^{p_1}\ar[r]_{p_2} & X_2
}
\end{align}
commutative, where $p_1\colon X_1\times X_2\to X_1$ and $p_2\colon X_1\times X_2\to X_2$ are canonical projections of the product $X_1\times X_2$. More generally, for $f_i\colon Y\to X_i$ there is a unique morphism $f=(f_1,\ldots,f_n)\colon Y\to\prod\limits_{i=1}^nX_i$ such that $p_i\cdot f=f_i$, where $p_i\colon\prod\limits_{i=1}^nX_i\to X_i$ are canonical projections.

\bb{Diagonal morphism.} Any object $X$ of a category with finite products $\bfC$ is equipped with a structure of comonoid in a unique way. Indeed, since the object $E$ is terminal, there is a unique morphism $\varepsilon_X\colon X\to E$. The morphisms $X\times X\xrightarrow{\id_X\times\varepsilon_X}X\times E\cong X$ and $X\times X\xrightarrow{\varepsilon_X\times\id_X}E\times X\cong X$ coincide with the canonical projections $p_1$ and $p_2$ respectively, so the morphism $\Delta_X\colon X\to X\times X$ can be uniquely found from the commutative diagram
\begin{align}
\xymatrix{
 & X\ar@{=}[dl]\ar@{-->}[d]^{\Delta_X}\ar@{=}[dr] \\
X\times E & X\times X\ar[l]^{\id_X\times\varepsilon_X} \ar[r]_{\varepsilon_X\times\id_X} & E\times X
}
\end{align}
This is $\Delta_X=(\id_X,\id_X)\colon X\to X\times X$, it is called {\it diagonal morphism} for the object $X$. One can check that the composition of both morphisms $(\id_X\times\Delta_X)\cdot\Delta_X\colon X\to X\times X\times X$ and $(\Delta_X\times\id_X)\cdot\Delta_X\colon X\to X\times X\times X$ with all three canonical projections $X\times X\times X\to X$ are equal to $\id_X$. Thus they coincide with the morphism $(\id_X,\id_X,\id_X)$ and hence equal to each other. Note that the diagonal morphism is cocommutative: $\Delta_X^\op=\Delta_X$.

In this way we obtain a unique comonoid $(X,\Delta_X,\varepsilon_X)\in\Comon(\bfC,\times)$ for any object $X\in\bfC$. Moreover, any morphism $f\colon X\to Y$ in $\bfC$ is a morphism of comonoids $(X,\Delta_X,\varepsilon_X)\to(Y,\Delta_Y,\varepsilon_Y)$. This means that $\Comon(\bfC,\times)$ coincide with $\bfC$ as a category. Since the monoidal product in $\Comon(\bfC,\times)$ coincide with the direct product $\times$ in $\bfC$, we have the strict monoidal equivalence
\begin{align} \label{ComonEquiv}
 \Comon(\bfC,\times)=\cocComon(\bfC,\times)=(\bfC,\times).
\end{align}

\bbo{Monoids and groups in the category $\bfC$}. \label{bbMonCatFP}
A monoid $\MM=(X,\mu_X,\eta_X)$ in $(\bfC,\times)$ is called {\it monoid in the category} $\bfC$. Due to the equivalence~\eqref{ComonEquiv} any structure of monoid on $X\in\bfC$ is compatible with the structure $(\Delta_X,\varepsilon_X)$.
Hence the monoid $\MM$ has a unique structure of bimonoid $\BB=(X,\mu_X,\eta_X,\Delta_X,\varepsilon_X)\in\Bimon(\bfC,\times)$ given by the diagonal morphism $\Delta_X\colon X\to X\times X$ and the unique morphism $\varepsilon_X\colon X\to E$, so we obtain the strict monoidal equivalence
\begin{align}
 \Bimon(\bfC,\times)=\Mon(\bfC,\times).
\end{align}
Note also that any bimonoid $\BB\in\Bimon(\bfC,\times)$ is cocommutative.

The category $\Set$ is a category with finite products. The terminal object in $\Set$ is the one-point set $E=\{*\}$. For a set $X\in\Set$ the diagonal morphism $\Delta_X\colon X\to X\times X$ has the form $\Delta_X(x)=(x,x)$ and the morphism $\varepsilon_X\colon X\to E$ maps any element $x\in X$ to the unique element $*$ of the set $E$. A monoid $\MM=(X,\mu_X,\eta_X)$ in $(\Set,\times)$ is a monoid in the usual sense: $\mu_X(x,y)=xy$ is a multiplication map and $\eta_X(*)=e$ gives the neutral element $e\in X$. It becomes a bimonoid in a unique way: $\BB=(X,\mu_X,\eta_X,\Delta_X,\varepsilon_X)$. The convolution of two maps $\alpha,\beta\colon X\to X$ is the point-wise product: $(\alpha\star\beta)(x)=\alpha(x)\beta(x)$; the neutral element for this convolution is the map $\eta_X\cdot\varepsilon_X\colon x\mapsto e$. Hence the monoid $\MM=(X,\mu_X,\eta_X)\in\Mon(\Set,\times)$ is a group iff the bimonoid $\BB=(X,\mu_X,\eta_X,\Delta_X,\varepsilon_X)$ is a Hopf monoid in $(\Set,\times)$. The role of antipode is played by the map $\zeta_X(x)=x^{-1}$.

For general category with finite products $\bfC$ the monoid $\MM=(X,\mu_X,\eta_X)$ in $\bfC$ is called {\it group in the category} $\bfC$ if the corresponding bimonoid $\BB=(X,\mu_X,\eta_X,\Delta_X,\varepsilon_X)$ is a Hopf monoid in $(\bfC,\times)$. It is enough to require one of the conditions $\mu_X\cdot(\zeta_X\times\id_X)\cdot\Delta_X=\eta_X\cdot\varepsilon_X$ or $\mu_X\cdot(\id_X\times\zeta_X)\cdot\Delta_X=\eta_X\cdot\varepsilon_X$, since they are equivalent in this case.

\bb{Algebraic monoids and algebraic groups.} \label{bbAlgMon}
The category $\AlgSet$ also has all finite products. Indeed, the terminal object is $\AA^0$, the product of two algebraic sets $X\subset\AA^n$ and $Y\subset\AA^m$ is the set-theoretic product $X\times Y$ canonically embedded into $\AA^n\times\AA^m=\AA^{n+m}$. As we already mentioned in p.~\ref{bbExMonCat} the category $\CommAlg$ has all finite coproducts and hence $\AffSch=\CommAlg^\op$ has all finite products. Since the categorical product of the affine schemes $\Spec\gR$ and $\Spec\gS$ is $\Spec(\gR\otimes\gS)$, we have $A(X\times Y)=A(X)\otimes A(Y)$ for any $X,Y\in\AlgSet$.

An {\it (affine) algebraic monoid} is a monoid in $\AlgSet$. Explicitly, this is a usual monoid that has a structure of algebraic set $X$ such that the multiplication map $X\times X\to X$, $(x,y)\mapsto xy$, is a morphism in $\AlgSet$.

An {\it (affine) algebraic group} is a group in $\AlgSet$. This is a group with a structure of algebraic set $X$ such that the maps $X\times X\to X$, $(x,y)\mapsto xy$, and $X\to X$, $x\mapsto x^{-1}$, are morphisms in $\AlgSet$.

More generally, one can define {\it affine monoid/group scheme} as a monoid/group in $\AffSch$. Structure of monoid on an affine scheme $X=\Spec\gR$ is given by the homomorphisms $\Delta_\gR\colon\gR\to\gR\otimes\gR$ and $\varepsilon_\gR\colon\gR\to\KK$ making $\gR\in\CommAlg$ into a commutative bialgebra. This is a structure of affine group scheme iff this bialgebra $(\gR,\Delta_\gR,\varepsilon_\gR)$ is a Hopf algebra. Thus, a commutative bialgebra/Hopf algebra is the same as an affine monoid/group scheme (modulo direction of their morphisms).

\section{Internal hom and representations}
\label{sec3}

A linear representation of an algebra $\A\in\Mon(\Vect,\otimes)$ on a vector space $V\in\Vect$ is an algebra homomorphism $\rho\colon\A\to\End(V)$, but from the categorical point of view the object $\End(V)\in\Mon(\Set,\times)$ is a monoid without a structure of vector space. To require the linearity of $\rho$ we need to equip $\End(V)$ with a linear structure. This can be done by using the internal hom-functor. This functor and its generalised version (adjunction with a parameter) allow us to consider representations in more general monoidal categories.

\subsection{Internal (co)hom-functor and its generalisation}

\bbo{Closed categories.} A symmetric monoidal category $\bfC=(\bfC,\otimes)$ is called {\it closed} if for any $Y\in\bfC$ the functor $-\otimes Y\colon\bfC\to\bfC$ has a right adjoint $\bfhom(Y,-)\colon\bfC\to\bfC$. In this case there exists a unique bifunctor $\bfhom=\bfhom_\bfC\colon\bfC^\op\otimes\bfC\to\bfC$, which is called {\it (internal) hom-functor}, and an isomorphism
\begin{align} \label{homDef}
 \theta=\theta_{X,Y,Z}\colon\Hom\big(X,\bfhom(Y,Z)\big)\isoright\Hom\big(X\otimes Y,Z\big)
\end{align}
natural in $X,Y,Z\in\bfC$ (see details in~\cite[\S~6.1]{Bor2} or~\cite[\S~4.7]{Mcl}).
The natural isomorphism~\eqref{homDef} is unique up to an automorphism of the functor $\bfhom(-,Y)$ (or, equivalently, of the functor $-\otimes Y$).%
\footnote{If $G\colon\bfD\to\bfC$ and $G'\colon\bfD\to\bfC$ are two right adjoints of a functor $F\colon\bfC\to\bfD$ with adjunctions $\alpha_{X,Z}\colon\Hom(FX,Z)\isoright\Hom(X,GZ)$ and $\alpha'_{X,Z}\colon\Hom(FX,Z)\isoright\Hom(X,G'Z)$, then there exists a unique natural isomorphism $\gamma_Z\colon GZ\to G'Z$ such that $\alpha'_{X,Z}=(\gamma_Z)_*\cdot\alpha_{X,Z}$  (see the proof of~\cite[\S~4.1, Cor.~1]{Mcl}). In particular, by considering the case $G=G'$ we deduce that the adjunction $\alpha$ is unique up to an automorphism $\gamma$.
}
The object $\bfhom(X,Y)\in\bfC$ is called {\it (internal) hom-object}. Due to symmetricity of $\bfC$ the functors $-\otimes Y$ and $Y\otimes-$ are isomorphic, so their right adjoints are the same (up to an isomorphism of functors).

\bb{Adjunction with a parameter.} \label{bbParam}
The classical example of a generalisation of the notion of internal hom appears in the theory of bimodules. Let $\A,\B,\C\in\Alg$. Let $M$ and $N$ be $(\A,\B)$- and $(\B,\C)$-bimodules respectively. Their tensor product is the $(\A,\C)$-module $M\otimes_\B N$. We have a natural isomorphism
\begin{align}
 \Hom_{(\A,\B)}(M,\bfhom_\C(N,K))\cong\Hom_{(\A,\C)}(M\otimes_\B N,K),
\end{align}
where $K$ is a $(\A,\C)$-bimodule, $\Hom_{(\A,\B)}(-,-)$ is the usual `external' Hom-functor in the category of the $(\A,\B)$-bimodules and $\bfhom_\C(N,K)$ is the set of right $\C$-module homomorphisms $N\to K$ equipped with the structure of $(\A,\B)$-module in a natural way.

In general, let $F\colon\bfC\times\bfP\to\bfC'$ be a bifunctor such that $F(-,Y)\colon\bfC\to\bfC'$ has a right adjoint $G_Y\colon\bfC'\to\bfC$ for any $Y\in\bfP$. Then there exists a unique bifunctor $G\colon\bfP^\op\times\bfC'\to\bfC$ such that $G_Y=G(Y,-)$ and the adjunction isomorphism
\begin{align} \label{FXYZG}
 \Hom_{\bfC'}\big(F(X,Y),Z\big)\cong\Hom_\bfC\big(X,G(Y,Z)\big)
\end{align}
is natural in $X\in\bfC$, $Y\in\bfP$, $Z\in\bfC'$ (see~\cite[\S~4.7, Th.~3]{Mcl}). This is so-called {\it adjunction with a parameter} with the parameter category $\bfP$.

\bb{Generalised internal hom.} \label{bbGenerIntHom}
If one uses the notation $F(X,Y)=X\otimes Y$ for the bifunctor $F$ in~\eqref{FXYZG}, then we can denote $G(Y,Z)$ by $\bfhom(Y,Z)$ in a generalised sense. The adjunction with a parameter has the form~\eqref{homDef} where $X\in\bfC$, $Y\in\bfP$, $Z\in\bfC'$.

In some important cases the categories $\bfC$, $\bfP$ and $\bfC'$ are full subcategories of a symmetric monoidal category such that $F(X,Y)=X\otimes Y\in\bfC'$ for any $X\in\bfC$ and $Y\in\bfP$. In particular, such situation takes place if $\bfP\subset\bfC=\bfC'$ for a symmetric monoidal category $\bfC$.

\bb{Evaluation.} \label{bbEval}
Let $Y\in\bfP$ and $Z\in\bfC'$. Substitute $X=\bfhom(Y,Z)\in\bfC$ to~\eqref{homDef} and take $\id_X$ in the left hand side, then we obtain so-called {\it evaluation morphism}
\begin{align} \label{evalYZ}
 \eval_{Y,Z}\colon\bfhom(Y,Z)\otimes Y\to Z
\end{align}
in the right hand side, that is $\eval_{Y,Z}=\theta(\id_{\bfhom(Y,Z)})$. This morphism is the counit of the adjunction $\big(-\otimes Y, \bfhom(Y,-), \theta^{-1}\big)$ and, in particular, it is natural in $Z\in\bfC'$; the isomorphism~\eqref{homDef} can be expressed through the evaluation as
\begin{align} \label{thetaEval}
 &\theta_{X,Y,Z}(f)=\eval_{Y,Z}\cdot(f\otimes\id_Y), &&f\in\Hom\big(X,\bfhom(Y,Z)\big),
\end{align}
(see~\cite[\S~4.1]{Mcl}).

\begin{Prop} \label{Propthetafg}
 Let $X,\wt X\in\bfC$, $Y\in\bfP$, $Z\in\bfC'$ and $f\colon X\to\bfhom(Y,Z)$, $g\colon\wt X\to X$ be morphisms in $\bfC$. Then
\begin{align} \label{thetafg}
 \theta_{\wt X,Y,Z}(f\cdot g)=
 \theta_{X,Y,Z}(f)\cdot(g\otimes\id_Y).
\end{align}
\end{Prop}

\noindent{\bf Proof.} By virtue of~\eqref{thetaEval} we derive
\begin{align*}
 \theta_{\wt X,Y,Z}(f\cdot g)=\eval_{Y,Z}\cdot\big((f\cdot g)\otimes\id_Y\big)= \eval_{Y,Z}\cdot(f\otimes\id_Y)\cdot(g\otimes\id_Y)=
 \theta_{X,Y,Z}(f)\cdot(g\otimes\id_Y),
\end{align*}
where we used the functoriality of $-\otimes-\colon\bfC\times\bfP\to\bfC'$ with respect to the first argument. \qed

\bb{Internal composition.} \label{bbIntComp}
Remind that the usual `external' composition is the map $\Hom(Y,Z)\times\Hom(X,Y)\to\Hom(X,Z)$ (morphism in $\Set$). One can define the analogous morphism for the internal hom objects as a morphism in $\bfC$, where the role of $\times$ is played by $\otimes$ (see~\cite[sec.~6.1]{Bor2}).

Let us consider the case of full subcategory $\bfP\subset\bfC=\bfC'$. For arbitrary objects $X,Y\in\bfP$ and $Z\in\bfC$ consider the morphism
\begin{align} \label{evaleval}
 \bfhom(Y,Z)\otimes\bfhom(X,Y)\otimes X\xrightarrow{\id\otimes\eval_{X,Y}}
 \bfhom(Y,Z)\otimes Y\xrightarrow{\eval_{Y,Z}}Z.
\end{align}
By applying the map
\begin{align*}
 \Hom\big(\bfhom(Y,Z)\otimes\bfhom(X,Y)\otimes X,Z\big)\xrightarrow{\theta^{-1}}\Hom\big(\bfhom(Y,Z)\otimes\bfhom(X,Y),\bfhom(X,Z)\big)
\end{align*}
to the morphism~\eqref{evaleval} we obtain a morphism
\begin{align*}
 &c_{X,Y,Z}\colon\bfhom(Y,Z)\otimes\bfhom(X,Y)\to\bfhom(X,Z), &&X,Y\in\bfP, Z\in\bfC.
\end{align*}
It is called {\it (internal) composition morphism}. This composition is associative in the sense that $c_{A,C,D}\cdot(\id\otimes c_{A,B,C})=c_{A,B,D}\cdot(c_{B,C,D}\otimes\id)$ for any $A,B,C\in\bfP$, $D\in\bfC$.

\bb{Internal end.} \label{bbIntEnd}
 Suppose again that $\bfP$ is a full subcategory of a symmetric monoidal category $\bfC=(\bfC,\otimes,I)=\bfC'$. By substituting $X=I$ and $Z=Y\in\bfP$ to \eqref{homDef} and taking the identification isomorphism $I\otimes Y=Y$ in the right hand side we obtain a morphism $u_Y:=\theta^{-1}(\id_Y)\colon I\to\bfhom(Y,Y)$. The diagram
\begin{align}
\xymatrix{
 I\otimes\bfhom(X,Y)\ar[d]_{u_Y\otimes\id}\ar@{=}[r] & \bfhom(X,Y)\ar@{=}[r] & \bfhom(X,Y)\otimes I\ar[d]^{\id\otimes u_X} \\
 \bfhom(Y,Y)\otimes\bfhom(X,Y)\ar[ru]_{c_{X,Y,Y}} & & \bfhom(X,Y)\otimes\bfhom(X,X)\ar[lu]^{c_{X,X,Y}}
}
\end{align}
commutes for any $X,Y\in\bfP$. In particular, for each $X\in\bfP$ the morphisms $c_X:=c_{X,X,X}$ and $u_X$ give a structure of monoid on the object $\bfhom(X,X)$. We denote this monoid by $\bfend(X)$.

\bb{Examples of closed categories.} \label{bbExIntHom}
\begin{itemize}
 \item In $(\Set,\times)$ the internal hom-object $\bfhom(X,Y)$ coincides with the `external' Hom-object $\Hom(X,Y)$. The evaluation is the map $\eval_{X,Y}\colon\Hom(X,Y)\times X\to Y$ that evaluates $f\in\Hom(X,Y)$ on an element of $X$. Internal composition coincides with the usual composition and $\bfend(X)=\End(X)$ is the usual monoid of the maps $X\to X$.
 \item For the objects $V$ and $W$ of the monoidal category $(\Vect,\otimes)$ the object $\bfhom(V,W)$ is the vector space of linear maps $V\to W$. In this case it coincides with $\Hom(V,W)$ as a set. The evaluation morphism $\eval_{V,W}\colon\bfhom(V,W)\otimes V\to W$ acts as $\eval(f\otimes v)=f(v)$, $f\colon V\to W$, $v\in V$. The composition $\bfhom(W,Z)\otimes\bfhom(V,W)\to\bfhom(V,Z)$ is given by the usual composition in the sense that it maps $f\otimes g$ to $f\cdot g$. The monoid $\bfend(V)\in\Mon(\Vect,\otimes)$ is the algebra of linear operators on $V$.
 \item The monoidal category $(\GrVect,\otimes)$ is also closed. As a vector space the internal hom in this category coincides with the internal hom in $(\Vect,\otimes)$, but it has additionally the structure of grading: the $k$-th component of the object $\bfhom(V,W)$ consists of linear maps $f\colon V\to W$ such that $f(V_l)\subset W_{k+l}$. By considering the zero component of the graded space $\bfhom(V,W)$ as a set we can identify this component with $\Hom(V,W)$ with the zero component of $\bfhom(V,W)$ (for this case the internal hom does not coincide with the external one as a set). The evaluation and composition look the same as for $(\Vect,\otimes)$. The monoid $\bfend(V)$ is the graded algebra of linear operators on $V$. This example is directly generalized to the case of the category of $\AA$-graded vector spaces for any abelian monoid $\AA\in\cMon(\Set,\times)$.
\end{itemize}

\bbo{Coclosed categories.} We also need dual notions to the notions of internal hom and end. Let $\bfC=(\bfC,\otimes)$ be a symmetric monoidal category such that for any object $Y\in\bfC$ the functor $-\otimes Y\colon\bfC\to\bfC$ has a left adjoint $\bfcohom(Y,-) \colon\bfC\to\bfC$, then there is a unique bifunctor  $\bfcohom=\bfcohom_\bfC\colon\bfC^\op\otimes\bfC\to\bfC$ with a unique (up to an automorphism) natural transformation
\begin{align}
 \vartheta=\vartheta_{X,Y,Z}\colon\Hom\big(\bfcohom(Y,X),Z\big)\isoright\Hom\big(X,Z\otimes Y)\big).
\end{align}
A symmetric monoidal category $\bfC=(\bfC,\otimes)$ satisfying this condition is called {\it coclosed}.
Let us call the object $\bfcohom(Y,X)$ {\it (internal) cohom-object}. It coincides with the internal hom-object $\bfhom(Y,X)$ in the opposite category $\bfC^\op$. The bifunctor $\bfcohom_\bfC\colon\bfC^\op\otimes\bfC\to\bfC$ coincides with the bifunctor $(\bfhom_{\bfC^\op})^\op$, it is called {\it (internal) cohom-functor}.

\bb{Generalised cohom.} \label{bbGenCohom}
 Let us consider the dual version of the adjunction with a parameter from p.~\ref{bbParam}. Denote the bifunctor $F\colon\bfC\times\bfP\to\bfC'$ as $F(X,Y)=X\otimes Y$. If $-\otimes Y\colon\bfC\to\bfC'$ has a left adjoint for any fixed $Y\in\bfP$, then we obtain a bifunctor $\bfcohom\colon\bfP^\op\times\bfC'\to\bfC$ with an isomorphism
\begin{align} \label{FXYZGcohom}
 \vartheta=\vartheta_{X,Y,Z}\colon\Hom_{\bfC'}\big(\bfcohom(Y,X),Z\big)\isoright\Hom_\bfC\big(X,Z\otimes Y\big)
\end{align}
natural in $X\in\bfC$, $Y\in\bfP$, $Z\in\bfC'$. The {\it coevaluation} is defined as the morphism
\begin{align} \label{coevXY}
 \coev_{Y,X}=\vartheta(\id_{\bfhom(Y,X)})\colon X\to\bfcohom(Y,X)\otimes Y,
\end{align}
which is natural in $X\in\bfC$.
Then the isomorphism~\eqref{FXYZGcohom} can be expressed as
\begin{align} \label{varethetacoev}
 &\vartheta(f)=(f\otimes\id_Z)\cdot\coev_{Y,X}, &&f\colon\bfcohom(Y,X)\to Z.
\end{align}

\bb{Internal coend.} \label{bbGenCoend}
 Suppose that $\bfP$ is a full subcategory of a symmetric monoidal category $\bfC=\bfC'=(\bfC,\otimes)$. By using categorical duality to the morphisms defined in p.~\ref{bbIntComp} and \ref{bbIntEnd} we obtain the {\it cocomposition}
\begin{align}
 d_{X,Y,Z}=\vartheta^{-1}\big((\id\otimes\coev_{Z,Y})\cdot\coev_{Y,X}\big)\colon\bfcohom(Z,X)\to\bfcohom(Y,X)\otimes\bfcohom(Z,Y) \label{dXYZ}
\end{align}
and morphisms $v_Y\colon\bfcohom(Y,Y)\to I$, where $X\in\bfC$, $Y,Z\in\bfP$. The pair $(d_Y:=d_{Y,Y,Y},v_Y)$ turns $\bfcohom(Y,Y)$ into a comonoid in $\bfC$. Denote it by $\bfcoend(Y)$ (not to confuse with `Coend of a functor').

\bb{Internal cohom for vector spaces.} Since the strong monoidal functor $(-)^*$ is an equivalence between $\FVect$ and $\FVect^\op$, the monoidal category $(\FVect,\otimes)$ is coclosed. The internal cohom-functor can be extended for the case $\bfP=\FVect$, $\bfC=\bfC'=\Vect$.

\begin{Prop} \label{PropCohomVect}
 For any $W\in\FVect$ the functor $-\otimes W\colon\Vect\to\Vect$ has a left adjoint $\bfcohom(W,-)\colon\Vect\to\Vect$, which coincides with the internal hom-functor:
\begin{align}
 \bfcohom(W,V)=\bfhom(W,V), &&V\in\Vect.
\end{align}
\end{Prop}

\noindent{\bf Proof.}  Let $(w_i)_{i=1}^m$, $(w^i)_{i=1}^m$ be dual bases of $W$ and $W^*$, that is $w^i(w_j)=\delta^i_j$. Note that for any $v\in V$ and $\xi\in W^*$ the element $\xi\otimes v\in W^*\otimes V$ can be considered as the linear operator $W\to V$ that maps $w\in W$ to $\xi(w)v\in V$ (in particular, $(w^i\otimes v)(w_j)=\delta^i_jv$). This gives a natural isomorphism $\bfhom(W,V)\cong W^*\otimes V$. Let us define the linear operator $\eta_V\colon V\to\bfhom(W,V)\otimes W=(W^*\otimes V)\otimes W$ by the formula $\eta_V(v)=\sum_{i=1}^m(w^i\otimes v)\otimes w_i$. It can be also defined in the form $\eta_V(v)=(\id_W\otimes\sigma_{W,V})(1_{\bfend(W)}\otimes v)$, where $1_{\bfend(W)}=u_W(1)\in\bfend(W,W)=W^*\otimes W$ is the unity of the algebra $\bfend(W)$, so it does not depend on the choice of the basis. Let us show that for any $V,Z\in\Vect$ and a linear operator $f\colon V\to Z\otimes W$ there is a unique $h\colon\bfhom(W,V)\to Z$ such that the diagram
\begin{align}
\xymatrix{
V\ar[rr]^{\eta_V\qquad}\ar[drr]_{f}&&\bfhom(W,V)\otimes W\ar[d]^{h\otimes\id_W} \\
 && Z\otimes W
}
\end{align}
commute. By decomposing the image $f$ on $v\in V$ we obtain $f(v)=\sum_{i=1}^m f^i(v)\otimes w_i$, where $f^i(v)$ is the value of $f(v)\in Z\otimes W=\bfhom(W^*,Z)$ on $w^i\in W^*$. The commutativity of the diagram implies that $h(w^i\otimes v)=f^i(v)=f(v)(w^i)$, so $h$ is unique. It exists due to the linearity of $f\colon V\to Z\otimes W$ and $f(v)\colon W^*\to Z$. We obtain the universal arrow $\eta_V$ from $V$ to the functor $-\otimes W$ and it is straightforward to check that $\eta_V$ is natural in $V$. Hence the functor $\bfhom(W,-)$ is a left adjoint of $-\otimes W$ (see~\cite[\S~4.1, Th.~2~(i)]{Mcl}). \qed

\subsection{(Co)representations of (co)monoids}
\label{secRepMon}

\bbo{Representation of a usual monoid or a group.} Remind that for any object $V$ of a category $\bfC$ the set $\End(V)=\Hom(V,V)$ is equipped with a structure of usual monoid (monoid in $\Set$). In a wide generality we can say that a representation is a monoid homomorphism $\rho\colon\MM\to\End(V)$. In particular, if $\MM$ is a group, then $\rho(m)$ is an automorphism of $V$ for any $m\in\MM$, so we obtain a representation of a group $\MM$ by automorphisms of the object $V\in\bfC$.

If the monoids $\MM$ and $\End(V)$ are equipped with some additional structure of the same type, one usually requires $\rho$ to preserve this structure.

A morphism from a representation $\rho\colon\MM\to\End(V)$ to a representation $\rho'\colon\MM\to\End(V')$ is a morphism $f\colon V\to V'$ in $\bfC$ such that $f\cdot\rho(m)=\rho'(m)\cdot f$ for any $m\in\MM$. In the pointless form this condition reads $f_*\cdot\rho=f^*\cdot\rho'$.

\bb{Representation of an algebra.} \label{bbRepAlg}
For example, let $V\in\Vect$; if the monoid $\MM$ has additionally a structure of vector space compatible with the monoid structure, then both $\MM$ and $\End(V)$ are algebras and we require $\rho\colon\MM\to\End(V)$ to be an algebra homomorphism, so we obtain the notion of `representation of an algebra'. More precisely, one should write there $\bfend(V)$ instead of $\End(V)$. As we will see below this case can be generalised to an arbitrary closed symmetric monoidal category (see p.~\ref{bbClCat}).

 Remind that the representations of a group $\GG$ (a group in $\Set$) or of a Lie algebra $\mathfrak{g}$ can be considered as particular cases of representations of algebras, these are the representations of the corresponding group algebra $\KK[\GG]$ or of the universal enveloping algebra $U(\mathfrak{g})$ respectively. In the same way representations of a monoid $\MM\in\Mon(\Set,\times)$ can be identified with representations of the algebra $\KK[\MM]$ consisting of the formal sums $\sum\limits_{m\in\MM}\alpha_m m$, $\alpha_m\in\KK$.

\bb{Representation of an algebraic monoid/group.} Another classical example is a structure of algebraic set on a monoid/group $\MM$. Suppose that this structure is compatible with the structure of monoid/group in the sense that $\MM$ is an algebraic monoid/group (see~p.~\ref{bbAlgMon}). For $V\in\FVect$ the monoid $\End(V)$ has a structure of algebraic monoid. In this way we obtain the representations of an algebraic monoid/group $\MM$ on a vector space $V$. In particular, the representation $\rho\colon\MM\to\End(\KK^n)$ gives an action $a\colon\MM\times\AA^n\to\AA^n$ of a monoid/group $\MM$ on the object $\AA^n\in\AlgSet$. In other words, `linear' representations of a monoid/group $\MM\in\Mon(\AlgSet,\times)$ is a particular case of an action of $\MM$.

\bb{The case of a closed category.} \label{bbClCat}
Let $\bfC=(\bfC,\otimes)$ be a closed symmetric monoidal category and $\MM=(X,\mu_X,\eta_X)\in\Mon(\bfC)$ be a monoid in $\bfC$. Define {\it representation} of $\MM$ on an object $V\in\bfC$ as a morphism $\rho\colon\MM\to\bfend(V)$ in the category $\Mon(\bfC)$. This is a morphism $\rho\colon X\to\bfhom(V,V)$ in $\bfC$ such that the diagrams
\begin{align} \label{rhoDiag}
\xymatrix{
 X\otimes X\ar[rr]^{\mu_X}\ar[d]_{\rho\otimes\rho}& & X\ar[d]^\rho \\
 \bfhom(V,V)\otimes\bfhom(V,V)\ar[rr]^{\qquad c_{V}}& &\bfhom(V,V)
}\qquad\qquad
\xymatrix{
I\ar[r]^{\eta_X}\ar[dr]_{u_V} & X\ar[d]^\rho \\
 & \bfhom(V,V)
}
\end{align}
are commutative.

\bb{The case of generalised internal hom.}
Now consider more general case. Let $\bfC=(\bfC,\otimes)$ be a symmetric monoidal category and let $\bfP$ be its full subcategory such that the functor $-\otimes Y\colon\bfC\to\bfC$ has a right adjoint for any $Y\in\bfP$, so we have the generalised $\bfhom\colon\bfP^\op\times\bfC\to\bfC$ in the sense of p.~\ref{bbGenerIntHom} (the case $\bfP\subset\bfC'=\bfC$). We get a monoid $\bfend(V)\in\Mon(\bfC)$ for any $V\in\bfP$. In this case we define {\it representation} of a monoid $\MM=(X,\mu_X,\eta_X)\in\Mon(\bfC)$ on an object $V\in\bfP$ as a morphism $\rho\colon\MM\to\bfend(V)$ in $\Mon(\bfC)$. Again, this is a morphism $\rho\colon X\to\bfhom(V,V)$ making the diagrams~\eqref{rhoDiag} commutative.

\bb{Representation as an action.} By substituting $Y=Z=V\in\bfP$ to~\eqref{homDef} we obtain the bijection
\begin{align} \label{thetaXVV}
 \theta\colon\Hom\big(X,\bfhom(V,V)\big)\isoright\Hom\big(X\otimes V,V\big).
\end{align}
Let us prove that the representations and actions on $V$ are in one-to-one correspondence via this bijection.

\begin{Lem} \label{LemRepAct}
 A morphism $\rho\colon X\to\bfhom(V,V)$ is a representation of the monoid $\MM$ iff $a=\theta(\rho)\colon X\otimes V\to V$ is an action of $\MM$.
\end{Lem}

\noindent{\bf Proof.} We check that the commutativity of~\eqref{rhoDiag} is equivalent to the commutativity of~\eqref{aDiag}. Prop.~\ref{Propthetafg} implies $\theta(\rho\cdot\mu_X)=\theta(\rho)\cdot(\mu_X\otimes\id_V)=a\cdot(\mu_X\otimes\id_V)\colon X\otimes X\otimes V\to V$. By virtue of the same proposition we have the equalities $\theta\big(c_{V}\cdot(\rho\otimes\rho)\big)= \eval_{V,V}\cdot(\id_{\bfhom(V,V)} \otimes\eval_{V,V})\cdot(\rho\otimes\rho\otimes\id_V)= \eval_{V,V}\cdot (\rho\otimes\id_V)\cdot(\id_X\otimes\eval_{V,V})\cdot (\id_X\otimes\rho\otimes\id_V)$ of morphisms $X\otimes X\otimes V\to V$. Due to~\eqref{thetaEval} one can rewrite the latter one in the form $\theta\big(c_{V}\cdot(\rho\otimes\rho)\big)=a\cdot(\id_X\otimes a)$.
 Thus $\rho\cdot\mu_X=c_{V}\cdot(\rho\otimes\rho)$ iff $a\cdot(\mu_X\otimes\id_V)=a\cdot(\id_X\otimes a)$. Analogously, we derive $\theta(\rho\cdot\eta_X)=a\cdot(\eta_X\otimes\id_V)$. By taking into account $\theta(u_V)=\id_V$ we deduce the equivalence of the condition $\rho\cdot\eta_X=u_V$ to the condition $a\cdot(\eta_X\otimes\id_V)=\id_V$. \qed

\bb{Morphisms of representations.} \label{bbMorRep}
Denote by $\Rep_\bfP(\MM)$ the category of the pairs $(V,\rho)$, where $V\in\bfP$ and $\rho$ is a representation of $\MM$ on $V$; a morphism $(V,\rho)\to(V',\rho')$ in this category is defined as a morphism $f\colon V\to V'$ in $\bfP$ satisfying $\bfhom(\id_V,f)\cdot\rho=\bfhom(f,\id_{V'})\cdot\rho'$. The last equation can be written in the diagram form as
\begin{align} \label{rhof}
\xymatrix{
X\ar[rr]^{\rho}\ar[d]_{\rho'} && \bfhom(V,V)\ar[d]^{\bfhom(\id_V,f)} \\
\bfhom(V',V')\ar[rr]^{\bfhom(f,\id_{V'})} && \bfhom(V,V')
}
\end{align}
To calculate $\bfhom(\id_V,f)$ and $\bfhom(f,\id_{V'})$ we will use the following general formula.

\begin{Prop} \label{Propthetahom}
 Let $f\colon V'\to V$ and $g\colon Z\to Z'$ be morphisms in $\bfP$ and $\bfC$ respectively. Then
\begin{align} \label{thetahom}
 \theta\big(\bfhom(f,g)\big)=g\cdot\eval_{V,Z}\cdot(\id_{\bfhom(V,Z)}\otimes f).
\end{align}
\end{Prop}

\noindent{\bf Proof.} From the naturality of $\theta_{X,V,Z}$ in $V$ and $Z$ we obtain
\begin{align} \label{diagthetahom}
\xymatrix{
 \Hom\big(X,\bfhom(V,Z)\big)\ar[r]^{\quad\theta}\ar[d]^{\bfhom(f,g)_*} &
\Hom\big(X\otimes V,Z\big)\ar[d]^{(\id_X\otimes f)^*\cdot g_*} \\
 \Hom\big(X,\bfhom(V',Z')\big)\ar[r]^{\quad\theta} &
\Hom\big(X\otimes V',Z'\big)
}
\end{align}
Let $X=\bfhom(V,Z)$. By taking $\id_X$ in the left upper corner of the diagram~\eqref{diagthetahom} we obtain the formula~\eqref{thetahom}. \qed

By virtue of Lemma~\ref{LemRepAct} we can regard the objects of the category $\Rep_\bfP(\MM)$ as objects of $\Lact(\MM)$. The definitions of morphisms in these categories coincide due to the following statement.

\begin{Lem} \label{LemRepMor}
 Let $a=\theta(\rho)\colon X\otimes V\to V$ and $a'=\theta(\rho)\colon X\otimes V'\to V'$ be actions of $\MM$ corresponding to the representations $\rho\colon \MM\to\bfend(V)$ and $\rho'\colon\MM\to\bfend(V')$, where $V,V'\in\bfP$. Let $f\colon V\to V'$ be a morphism in $\bfP$. Then the commutativity of~\eqref{rhof} is equivalent to the commutativity of the left diagram~\eqref{LactMor}.
\end{Lem}

\noindent{\bf Proof.} Application of the formulae~\eqref{thetafg} and~\eqref{thetahom} gives
\begin{align*}
 \theta\big(\bfhom(f,\id_{V'})\cdot\rho'\big)=
\eval_{V',V'}\cdot(\id_{\bfhom(V',V')}\otimes f)\cdot(\rho'\otimes&\id_V)= \notag \\
\eval_{V',V'}\cdot(\rho'\otimes\id_{V'})\cdot(\id_{X}\otimes f)=
\theta(\rho')\cdot(\id_{X}\otimes f)=
 &a'\cdot(\id_X\otimes f), \\
 \theta\big(\bfhom(\id_V,f)\cdot\rho\big)=f\cdot\eval_{V,V}\cdot (\rho\otimes\id_V)=f\cdot\theta(\rho)=&f\cdot a,
\end{align*}
where we also used~\eqref{thetaEval}. \qed

The Lemmas~\ref{LemRepAct} and \ref{LemRepMor} have the following corollary.

\begin{Th} \label{ThRepLact}
Let $\bfP$ be a full subcategory of a symmetric monoidal category $\bfC=(\bfC,\otimes)$. Suppose that the functor $-\otimes V\colon\bfC\to\bfC$ has a right adjoint for any $V\in\bfP$. Let $\MM=(X,\mu_X,\eta_X)\in\Mon(\bfC)$. Then the morphisms~\eqref{thetaXVV} gives a fully faithful functor embedding the category $\Rep_\bfP(\MM)$ into $\Lact(\MM)$. In particular, if the symmetric monoidal category $\bfC=(\bfC,\otimes)$ is closed then we obtain an equivalence between the categories $\Rep_\bfC(\MM)$ and $\Lact(\MM)$.
\end{Th}

\bbo{Corepresentations of comonoids.} \label{bbCorep}
Let us dualise the notions and results of p.~\ref{bbClCat}--\ref{bbMorRep}. Let $\bfP$ be a full subcategory of a symmetric monoidal category $\bfC=(\bfC,\otimes)$ such that the functor $-\otimes V$ has left adjoint for any $V\in\bfP$, so we have the generalised cohom-functor $\bfcohom\colon\bfP^\op\times\bfC\to\bfC$ with the adjunction~\eqref{FXYZGcohom}. We call by {\it corepresentation} of a comonoid $\OO=(X,\Delta_X,\varepsilon_X)\in\Comon(\bfC)$ on $V\in\bfP$ a morphism $\omega\colon\bfcoend(V)\to\OO$ in $\Comon(\bfC)$. This is a morphism $\omega\colon\bfcohom(V,V)\to X$ in $\bfC$ such that the diagrams
\begin{align*}
\xymatrix{
 \bfcohom(V,V)\ar[d]_{d_{V}}\ar[rr]^{\qquad\omega}& & X\ar[d]^{\Delta_X} \\
 \bfcohom(V,V)\otimes\bfcohom(V,V)\ar[rr]^{\qquad\qquad\omega\otimes\omega}& &X\otimes X
}\qquad\qquad
\xymatrix{
\bfcohom(V,V)\ar[r]^{\qquad\omega}\ar[dr]_{v_V} & X\ar[d]^{\varepsilon_X} \\
 &I 
}
\end{align*}
are commutative. Denote by $\Corep_\bfP(\OO)$ the category whose objects are the pairs $(V,\omega)$, where $V\in\bfP$ and $\omega$ is a corepresentation of $\OO$ on $V$, and morphisms $(V,\omega)\to(V',\omega')$ are morphisms $f\colon V\to V'$ in $\bfP$ such that
\begin{align} \label{omegaf}
\xymatrix{
 \bfcohom(V',V)\ar[rrr]^{\bfcohom(f,\id_V)}\ar[d]^{\bfcohom(\id_{V'},f)} &&& \bfcohom(V,V)\ar[d]^{\omega} \\
\bfcohom(V',V')\ar[rrr]^{\qquad\omega'} &&& X
}
\end{align}

By dualising Prop.~\ref{Propthetahom} and Theorem~\ref{ThRepLact} we obtain the following statements.

\begin{Prop} \label{Propthetacohom}
 Let $f\colon V'\to V$ and $g\colon Z\to Z'$ be morphisms in $\bfP$ and $\bfC$ respectively. Then
\begin{align} \label{thetacohom}
 \vartheta\big(\bfcohom(f,g)\big)=(\id_{\bfcohom(V',Z')}\otimes f)\cdot\coev_{V',Z'}\cdot g.
\end{align}
\end{Prop}

\begin{Th} \label{ThCorepLcoact}
 Let $\OO=(X,\Delta_X,\varepsilon_X)\in\Comon(\bfC)$. The isomorphisms
\begin{align} \label{thetaP}
 &\vartheta_{V,V,X}\colon\Hom\big(\bfcohom(V,V),X\big)\cong\Hom\big(V,X\otimes V\big), &&V\in\bfP,
\end{align}
give a fully faithful functor $\Corep_\bfP(\OO)\hookrightarrow\Lcoact(\OO)$. In other words, the morphism $\omega\colon\bfcohom(V,V)\to X$ is a corepresentation of $\OO$ iff $\delta=\vartheta(\omega)$ is a coaction of $\OO$ on $V\in\bfP$, and a morphism $f\colon V\to V'$ in $\bfP$ makes the diagram~\eqref{omegaf} to be commutative for some corepresentations $\omega,\omega'$ iff it makes the right diagram~\eqref{LactMor} commutative with the coactions $\delta=\vartheta(\omega)$, $\delta'=\vartheta(\omega')$.

In particular, if the symmetric monoidal category $\bfC=(\bfC,\otimes)$ is coclosed then we obtain an equivalence between the categories $\Corep_\bfC(\OO)$ and $\Lcoact(\OO)$.
\end{Th}

\subsection{Translation of (co)representations under monoidal functors}
\label{secTransRep}

As we saw in p.~\ref{bbFact} a (co)lax functor translates (co)actions to (co)actions. Theorems~\ref{ThRepLact} and \ref{ThCorepLcoact} imply that such functor translates the corresponding (co)representations to each other. Here we describe the translation of (co)representations explicitly.

\bb{Translation of internal hom.} Let $\bfC=(\bfC,\otimes)$ and $\bfD=(\bfD,\odot)$ be symmetric monoidal categories. Let $\bfP\subset\bfC$ and $\bfQ\subset\bfD$ be their full subcategories such that the functors $-\otimes V\colon\bfC\to\bfC$ and $-\odot W\colon\bfD\to\bfD$ have right adjoints for each $V\in\bfP$ and $W\in\bfQ$, so we have generalised internal hom-functors $\bfhom\colon\bfP^\op\times\bfC\to\bfC$ and $\bfhom\colon\bfQ^\op\times\bfD\to\bfD$.

Let $F\colon(\bfC,\otimes)\to(\bfD,\otimes)$ be a lax monoidal functor with the monoidal structure morphisms $\phi_{X,Y}\colon FX\odot FY\to F(X\otimes Y)$, $\varphi\colon I_\bfD\to FI_\bfC$. Suppose $F(\bfP)\subset\bfQ$. For arbitrary $V\in\bfP$ and $Z\in\bfC$ consider the composition
\begin{align} \label{phiFeval}
 F\big(\bfhom(V,Z)\big)\odot FV\xrightarrow{\phi}
 F\big(\bfhom(V,Z)\otimes V\big)\xrightarrow{F(\eval_{V,Z})}FZ.
\end{align}
By applying the adjunction
\begin{align}
 \theta^{-1}\colon\Hom\big(F\big(\bfhom(V,Z)\big)\odot FV,FZ\big)\isoright
\Hom\big(F\big(\bfhom(V,Z)\big),\bfhom(FV,FZ)\big)
\end{align}
to~\eqref{phiFeval} we obtain the following morphism in $\bfD$:
\begin{align} \label{PhiDef}
 \Phi_{V,Z}=\theta^{-1}\big(F(\eval_{V,Z})\cdot\phi_{\bfhom(V,Z),V}\big)\colon
F\big(\bfhom(V,Z)\big)\to\bfhom(FV,FZ).
\end{align}
Thus we get a collection of morphisms $\Phi_{V,Z}$, $V\in\bfP$, $Z\in\bfC$, making the diagram
\begin{align} \label{PhiDiag}
\xymatrix{
 F\big(\bfhom(V,Z)\big)\odot FV\ar[rr]^{\Phi_{V,Z}\odot\id_{FV}}\ar[d]_\phi &&
 \bfhom(FV,FZ)\odot FV\ar[d]^{\eval_{FV,FZ}} \\
 F\big(\bfhom(V,Z)\otimes FV\big)\ar[rr]^{F(\eval_{VZ})} && FZ
}
\end{align}
commute (the equivalence of the commutativity of~\eqref{PhiDiag} and the definition~\eqref{PhiDef} follows from the formula~\eqref{thetaEval}). The naturality of $\phi_{X,Y}$ in $X\in\bfC$ allows to prove the following properties of $\Phi_{V,Z}$.

\begin{Prop} \label{PropPhiphi}
 For any objects $X,Z\in\bfC$, $V\in\bfP$ and a morphism $f\colon X\to\bfhom(V,Z)$ in $\bfC$ we have the formula
\begin{align} \label{Phif}
 \theta\big(\Phi_{V,Z}\cdot F(f)\big)=F\big(\theta(f)\big)\cdot\phi_{X,V}.
\end{align}
\end{Prop}

\noindent{\bf Proof.} Note that~\eqref{Phif} is equivalent to the commutativity of the diagram
\begin{align} \label{PhifDiag}
\xymatrix{
FX\odot FV\ar[d]_\phi\ar[rr]^{F(f)\odot\id_{FV}\qquad} && F\big(\bfhom(V,Z)\big)\odot FV\ar[rr]^{\Phi_{V,Z}\odot\id_{FV}} &&
 \bfhom(FV,FZ)\odot FV\ar[d]^{\eval_{FV,FZ}} \\
F(X\otimes V)\ar[rr]^{F(f\otimes\id_{V})\qquad} && F\big(\bfhom(V,Z)\otimes FV\big)\ar[rr]^{F(\eval_{VZ})} && FZ
}
\end{align}
By adding the vertical arrow $\phi$ to the centre we obtain two diagrams. Commutativity of the left one follows from the naturality of $\phi$, while the right one is exactly the diagram~\eqref{PhiDiag}. \qed

\noindent{\bf Warning.} Even when $F$ is strong monoidal, i.e. all $\phi_{X,Y}$ and $\phi$ are isomorphisms, we can not guarantee that $\Phi_{V,Z}$ are also isomorphisms.

\bb{Translation of the internal end.} Remind that the lax monoidal functor $F\colon\bfC\to\bfD$ induces the functor $\Mon(F)\colon\Mon(\bfC)\to\Mon(\bfD)$ which translates each monoid $\MM=(X,\mu_X,\eta_X)\in\Mon(\bfC)$ to the monoid  $\Mon(F)\MM=(FX,\mu_{FX},\eta_{FX})\in\Mon(\bfD)$, where  $\mu_{FX}\colon FX\odot FX\to FX$ and $\eta_{FX}\colon I_\bfD\to FX$ are the compositions~\eqref{muFX}. In particular, if $V\in\bfP$ the monoid $\bfend(V)$ is mapped to the monoid $\Mon(F)\big(\bfend(V)\big)$. The latter is the object $F\big(\bfhom(V,V)\big)$ with the structure morphisms $F(c_{V})\cdot\phi_{\bfhom(V,V),\bfhom(V,V)}$ and $F(u_V)\cdot\varphi$.

On the other hand, since $FV\in\bfQ$ we have the monoid
\begin{align}
 \bfend(FV)=\big(\bfhom(FV,FV),c_{FV},u_{FV}\big)\in\Mon(\bfD).
\end{align}
We get a relationship between these two monoids by means of the morphisms~\eqref{PhiDef}.

\begin{Prop} \label{PropFend}
 For any $V\in\bfP$ the morphism $\Phi_{V,V}\colon F\big(\bfhom(V,V)\big)\to\bfhom(FV,FV)$ in $\bfD$ gives a morphism $\Mon(F)\big(\bfend(V)\big)\to\bfend(FV)$ in $\Mon(\bfD)$. In particular, if $\Phi_{V,V}$ is an isomorphism in $\bfD$ then the monoids $\Mon(F)\big(\bfend(V)\big)$ and $\bfend(FV)$ are isomorphic to each other (as objects of $\Mon(\bfD)$).
\end{Prop}

\noindent{\bf Proof.} The first sentence of the proposition will be proved in p.~\ref{bbFrep} in more general settings (see~Remark~\ref{RemFend}). The second sentence follows from the first one and Prop.~\ref{PropMonIso}. \qed

\bb{Translation of representations.} \label{bbFrep}
Since $F(\bfP)\subset\bfQ$, the restriction of the functor~\eqref{FLact} gives the functor
\begin{align} \label{FRep}
 \Rep_\bfP(\MM)\to\Rep_\bfQ(\wt\MM),
\end{align}
where $\wt\MM=\Mon(F)\MM$. Let us describe it explicitly.

\begin{Prop} \label{PropFRep}
 The functor~\eqref{FRep} maps an object $(V,\rho)\in\Rep_\bfP(\MM)$ to $(FV,\wt\rho)$, where $\wt\rho\colon\wt\MM\to\bfend(FV)$ is the composition
\begin{align} \label{FrhoPhi}
 FX\xrightarrow{F\rho}F\big(\bfhom(V,V)\big)\xrightarrow{\Phi_{V,V}}\bfhom(FV,FV).
\end{align}
\end{Prop}

\noindent{\bf Proof.} Let $a=\theta(\rho)$. The functor~\eqref{FLact} maps $(V,a)$ to $(FV,\wt a)$, where $\wt a=F(a)\cdot\phi_{X,V}$. By applying Prop.~\ref{PropPhiphi} to $f=\rho$ we derive $\theta(\Phi_{V,V}\cdot F\rho)=F\big(\theta(\rho)\big)\cdot\phi=F(a)\cdot\phi=\wt a$, so the action $\wt a$ corresponds to the representation~\eqref{FrhoPhi}. \qed

\begin{Rem} \normalfont \label{RemFend}
 The fact that~\eqref{FrhoPhi} is a morphism of monoids $\wt\MM\to\bfend(FV)$ follows immediately from Prop.~\ref{PropFend} and functoriality of $\Mon(F)$. Conversely, by applying Prop.~\ref{PropFRep} to the case $X=\bfhom(V,V)$, $\rho=\id_X$ and by taking into account the statements written in p.~\ref{bbFact} we see that $\Phi_{V,V}$ is a representation. This implies the first sentence of Prop.~\ref{PropFend}.
\end{Rem}

\bbo{Translation of internal cohom, coend and of corepresentations.} In the same setting suppose that the functors $-\otimes V\colon\bfC\to\bfC$ and $-\odot W\colon\bfD\to\bfD$ have left adjoints for each $V\in\bfP$ and $W\in\bfQ$ instead of right adjoints, so we have generalised internal cohom-functors $\bfcohom\colon\bfP^\op\times\bfC\to\bfC$ and $\bfcohom\colon\bfQ^\op\times\bfD\to\bfD$. 

Let $F=(F,\varphi,\phi)\colon(\bfC,\otimes)\to(\bfD,\odot)$ be a colax monoidal functor such that $F(\bfP)\subset\bfQ$. For arbitrary $V\in\bfP$ and $Z\in\bfC$ consider the composition
\begin{align} \label{Fcoevphi}
 FZ\xrightarrow{F(\coev_{V,Z})} F\big(\bfcohom(V,Z)\otimes V\big)\xrightarrow{\phi}F\big(\bfcohom(V,Z)\big)\odot FV.
\end{align}
By applying the adjunction we obtain
\begin{align}
 \Phi_{V,Z}=\vartheta^{-1}\big(\phi_{\bfcohom(V,Z),V}\cdot F(\coev_{V,Z})\big)\colon
\bfcohom(FV,FZ)\to F\big(\bfcohom(V,Z)\big),
\end{align}
where $V\in\bfP$, $Z\in\bfC$. These morphisms satisfies
\begin{align}
 \vartheta\big(F(f)\cdot\Phi_{V,Z}\big)=\phi_{X,V}\cdot F\big(\vartheta(f)\big)
\end{align}
for any $X,Z\in\bfC$, $V\in\bfP$ and $f\in\Hom_\bfC\big(\bfcohom(V,Z),X\big)$.

\begin{Prop}
For any $V\in\bfP$ the morphism $\Phi_{V,V}$ gives a morphism of comonoids
\begin{align}
 \Phi_{V,V}\colon\bfcoend(FV)\to\Comon(F)\big(\bfcoend(V)\big).
\end{align}
In particular, if $\Phi_{V,V}$ is an isomorphism in $\bfD$, then these monoids are isomorphic.

Let $\OO=(X,\Delta_X,\varepsilon_X)\in\Mon(\bfC)$ and $\wt\OO=\Comon(F)\OO$. Then the restriction of the functor~\eqref{FLcoact} is the functor
\begin{align} \label{FCorep}
 \Corep_\bfP(\OO)\to\Corep_\bfQ(\wt\OO).
\end{align}
It maps the object $(V,\omega)$ to $(FV,\wt\omega)$, where $\wt\omega\colon\bfcoend(FV)\to\wt\OO$ is a corepresentation of $\wt\OO$ on $FV$ given by the composition
\begin{align}
 \bfcoend(FV)\xrightarrow{\Phi_{V,V}}F\big(\bfcoend(V)\big)\xrightarrow{F\omega} FX.
\end{align}
\end{Prop}

\bbo{Translation of monoidal product of (co)representations.} \label{bbTrTP}
Let $F\colon\bfC\to\bfD$ be a symmetric strong monoidal functor such that $F(\bfP)\subset\bfQ$. Let $\BB=(X,\mu_X,\eta_X,\Delta_X,\varepsilon_X)$ be a bimonoid in $\bfC$ and $\wt\BB=\Bimon(F)\BB$ be the corresponding bimonoid in $\bfD$. Then the restriction of the functors~\eqref{FLactB} are the functors~\eqref{FRep} and~\eqref{FCorep}, where we have $\MM=(X,\mu_X,\eta_X)\in\Mon(\bfC)$, $\wt\MM=\Mon(F)\MM\in\Mon(\bfD)$ and $\OO=(X,\Delta_X,\varepsilon_X)\in\Comon(\bfC)$, $\wt\OO=\Comon(F)\OO\in\Comon(\bfD)$. Since the functors~\eqref{FLactB} are symmetric strong monoidal, the restriction gives the symmetric strong monoidal functors
\begin{align} \label{FRepB}
 &\Rep_\bfP(\BB)\to\Rep_\bfQ(\wt\BB), &
 &\Corep_\bfP(\BB)\to\Corep_\bfQ(\wt\BB).
\end{align}

For example, let $\rho$ and $\pi$ be representations of the monoid $\MM$ on $W\in\bfP$ and $Z\in\bfP$. Their monoidal product in $\Rep_\bfP(\BB)$ is a representation $\tau$ on $W\otimes Z$. The functor~\eqref{FRep} maps them to representations $\wt\rho$, $\wt\pi$ and $\wt\tau$ on $FW$, $FZ$ and $F(W\otimes Z)$ respectively. The monoidal product of $\wt\rho$ and $\wt\pi$ in $\Rep_\bfQ(\wt\BB)$ is a representation $\lambda$ on $FW\odot FZ$. The isomorphism $\phi_{W,Z}\colon FW\odot FZ\isoright F(W\otimes Z)$ gives the isomorphism $\big(FW\odot FZ,\lambda\big)\isoright\big(F(W\otimes Z),\wt\tau\big)$ in $\Rep_\bfQ(\wt\BB)$. In the same way one can describe the isomorphism of the corresponding corepresentations.

\begin{Rem} \normalfont
 The strong monoidality of the functors~\eqref{FRepB} follows from the strong monoidality of the functors~\eqref{FLactB}. Alternatively, one can prove this in a direct way by using the formula~\eqref{thetahom} and Prop.~\ref{PropMonIso}.
\end{Rem}

\section{Quantum linear spaces}
\label{sec4}

{\bf Idea.} Consider the sets $\KK^n$ and $\KK^m$ as vector (linear) spaces, that is as objects of $\Vect$ or $\FVect$. The morphisms between these vector spaces are linear maps $\KK^n\to\KK^m$ (matrices $m\times n$ over $\KK$). If we (partially) forget linear structure of these spaces we obtain the affine spaces $\AA^n$ and $\AA^m$, which are the objects of $\AlgSet\subset\AffSch$. Since $\AffSch=\CommAlg^\op$ and $A(\AA^n)=\KK[x^1,\ldots,x^n]$, the morphisms $\AA^n\to\AA^m$ is in one-to-one correspondence with the algebra homomorphisms $f\colon\KK[y^1,\ldots,y^m]\to\KK[x^1,\ldots,x^n]$. Namely, if a homomorphism $f$ is given by the images $f(y^i)=P_i(x^1,\ldots,x^n)\in\KK[x^1,\ldots,x^n]$ for some polynomials $P_i$, then the corresponding morphism $\Phi\colon\AA^n\to\AA^m$ is the map $\Phi(x^1,\ldots,x^n)=\big(P_1(x^1,\ldots,x^n),\ldots,P_m(x^1,\ldots,x^n)\big)$. Now let us recall back the linear structure of $\KK^n=\AA^n$ and $\KK^m=\AA^m$. The map $\Phi$ preserves the linear structure iff all $P_i$ are homogeneous polynomials of order $1$, that is $P_i(x^1,\ldots,x^n)=\sum_{j=1}^n a_{ij}x^j$ for some $a_{ij}\in\KK$. This means exactly that the homomorphism $f$ preserves the grading of the algebras $\KK[y^1,\ldots,y^m]$ and $\KK[x^1,\ldots,x^n]$ (the $k$-th component of $\KK[x^1,\ldots,x^n]$ is the space of all the homogeneous polynomials of order $k$). Any linear map $\KK^n\to\KK^m$ is arisen in this way: it is given by the $m\times n$ matrix $(a_{ij})$.

Hence a natural candidates to a quantum analogue of a vector spaces are graded algebras or at least some of them, considered as objects of (a subcategory of) the category $\GrAlg^\op$. In~\cite{Manin87,Manin88} Yuri Manin proposed a subcategory $\FQA^\op$ to generalise the finite-dimensional vector spaces. He introduced the notion `{\it quantum linear space}' for this case.

\subsection{Operations with quadratic algebras}

\bbo{Manin's binary operations.} \label{bbManinOpns}
 Let $\A,\B\in\QA$ be arbitrary quadratic algebras (here we do not suppose that the quadratic algebras are finitely generated). They can be presented in the form
\begin{align}
 &\A=TV/(R), & \B=TW/(S),
\end{align}
where $V,W\in\Vect$ and $(R)\subset TV$, $(S)\subset TW$ are ideals generated by subspaces $R\subset V^{\otimes2}$, $S\subset W^{\otimes2}$ respectively. Define the following operations~\cite{Manin87,Manin88}.
\begin{itemize}
 \item {\it Manin white product:} $\A\circ\B=T(V\otimes W)/(R_{\mathrm w})$, where
\begin{align*}
 &R_{\mathrm w}=(\id_V\otimes\sigma_{V,W}\otimes\id_W)(R\otimes W\otimes W+V\otimes V\otimes S)\subset V\otimes W\otimes V\otimes W.
\end{align*}
 \item {\it Manin black product:} $\A\bullet\B=T(V\otimes W)/(R_{\mathrm b})$, where
\begin{align*}
 &R_{\mathrm b}=(\id_V\otimes\sigma_{V,W}\otimes\id_W)(R\otimes S)\subset V\otimes W\otimes V\otimes W.
\end{align*}
 \item {\it (Even) tensor product:} $\A\otimes\B=T(V\oplus W)/\big(R\oplus[V,W]\oplus S\big)$, where
\begin{multline*}
 [V,W]=\{v\otimes w-w\otimes v\mid v\in V,w\in W\}\subset(V\otimes W)\oplus(W\otimes V)\subset(V\oplus W)^{\otimes2}.
\end{multline*}
 \item {\it Odd tensor product:} $\A\ootimes\B=T(V\oplus W)/\big(R\oplus[V,W]_+\oplus S\big)$, where
\begin{multline*}
 [V,W]_+=\{v\otimes w+w\otimes v\mid v\in V,w\in W\}\subset(V\otimes W)\oplus(W\otimes V)\subset(V\oplus W)^{\otimes2}.
\end{multline*}
\end{itemize}

\bbo{Properties of Manin's binary operations}~\cite{Manin88}. The Manin white product and the even tensor product coincide with the Manin product $\A\circ\B$ and the usual tensor product $\A\otimes\B$ defined in p.~\ref{bbExMonCat} for general graded algebras. All four operations are bifunctors $\QA\times\QA\to\QA$ equipping $\QA$ with monoidal structures. The unit objects of $(\QA,\circ)$ and $(\QA,\bullet)$ are $\KK[u]$ and $\KK[\epsilon]/(\epsilon^2)$ respectively. The unit object of $(\QA,\otimes)$ and $(\QA,\ootimes)$ is $\KK$. The subcategory $\FQA\subset\QA$ inherits all four monoidal structures. The inclusion $R_{\mathrm b}\subset R_{\mathrm w}$ defines an epimorphism $\A\bullet\B\to\A\circ\B$.

\bb{The functor $(-)_1$.}
If we have a graded algebra $\A=\bigoplus\limits_{k\in\NN_0}\A_1$ we can take its first graded component $\A_1$. In this way we obtain the functor $(-)_1\colon\GrAlg\to\Vect$, which maps a graded homomorphism $f\in\Hom_{\GrAlg}(\A,\B)$ to its first order component $f_1\colon\A_1\to\B_1$. Restriction gives functors $(-)_1\colon\QA\to\Vect$ and $(-)_1\colon\FQA\to\FVect$, which translate a quadratic algebra $\A=TV/(R)$ to the vector space $V$. If $\A\in\QA$, then a graded homomorphism $f\colon\A\to\B$ can be uniquely restored by the component $f_1$. Hence the functors $(-)_1\colon\QA\to\Vect$ and $(-)_1\colon\FQA\to\FVect$ are faithful. By applying these functors to the result of the Manin's binary operations we obtain the following natural isomorphisms:
\begin{align} \label{AB1}
 &(\A\circ\B)_1=(\A\bullet\B)_1=\A_1\otimes\B_1, &
 &(\A\otimes\B)_1=(\A\ootimes\B)_1=\A_1\oplus\B_1.
\end{align}
This means that the functors $(-)_1\colon\QA\to\Vect$ and $(-)_1\colon\FQA\to\FVect$ are equipped with the following strong monoidal structures
\begin{align}
 &(\QA,\circ)\to(\Vect,\otimes), &
 &(\FQA,\circ)\to(\FVect,\otimes), \\
 &(\QA,\bullet)\to(\Vect,\otimes), &
 &(\FQA,\bullet)\to(\FVect,\otimes), \\
 &(\QA,\otimes)\to(\Vect,\oplus), &
 &(\FQA,\otimes)\to(\FVect,\oplus), \\
 &(\QA,\ootimes)\to(\Vect,\oplus), &
 &(\FQA,\ootimes)\to(\FVect,\oplus).
\end{align}
The fullness of the functor $(-)_1\colon\QA\to\Vect$ and the left isomorphism~\eqref{AB1} imply that the epimorphism $\A\bullet\B\to\A\circ\B$ is natural in $\A,\B\in\QA$.

Note also that the functors $T\colon\Vect\to\QA$ and $T\colon\FVect\to\FQA$ considered in p.~\ref{bbExMonFunc} are left adjoints of $(-)_1\colon\QA\to\Vect$ and $(-)_1\colon\FQA\to\FVect$ respectively.

We see that the binary operations $\circ$ and $\bullet$ are related with the tensor product of vector spaces, while $\otimes$ and $\ootimes$ are related with the direct sum $\oplus$. Below we introduce two embedding of $\FVect$ into $\FQA^{\op}$: `even' and `odd'. The former one will relate the monoidal products $\otimes$ and $\oplus$ on $\FVect$ to $\circ$ and $\otimes$ on $\FQA$, while the odd embedding will relate $\otimes$ and $\oplus$ to $\bullet$ and $\ootimes$ respectively.

\bb{Purely even/odd quadratic super-algebras.}
Let $\ZZ_2=\ZZ/2\ZZ=\{\bar0,\bar1\}$. A {\it quadratic super-algebra} over $\KK$ is a quadratic algebra $\A=TV/(R)$ with the structure of $\ZZ_2$-grading compatible with multiplication and with $\NN_0$-grading in the sense that each component $\A_k$ is a direct sum of its subspaces of even and odd elements: $\A=\bigoplus\limits_{k\in\NN_0}\A_k=\bigoplus\limits_{k\in\NN_0}(\A_k)_{\bar0}\oplus(\A_k)_{\bar1}$. Such algebra is generated by (basis) elements of the subspaces $V_{\bar0}=(A_1)_{\bar0}$ and $V_{\bar1}=(A_1)_{\bar1}$.
We call it {\it commutative} if $ab=-ba$ for $a,b\in V_{\bar1}$ and $ab=ba$ for $a\in V_{\bar0}$, $b\in V$ (this is exactly the commutativity defined in p.~\ref{bbEHP} for the monoidal category of super-vector spaces~\cite[\S~1.1, 1.2]{Del}). We call a quadratic super-algebra $\A=TV/(R)$ {\it purely even} or {\it purely odd} if $V=V_{\bar0}$ or $V=V_{\bar1}$ respectively. The commutative purely even and odd quadratic super-algebras over $\KK$ form categories, which we denote by $\CommQSAe$ and $\CommQSAo$.  They can be identified with subcategories of $\QA$ consisting of commutative $\A\in\QA$ (i.e. $ab=ba$ $\;\forall\,a,b\in\A$) and of $\A\in\QA$ satisfying $a_1b_1=-b_1a_1$ $\;\forall\,a_1,b_1\in\A_1$ respectively. The additional condition that they are finitely generated gives the subcategories
\begin{align*}
 &\FCommQSAe=\CommQSAe\cap\FQA, \\ &\FCommQSAo=\CommQSAo\cap\FQA.
\end{align*}

\bb{The functors $S$ and $\Lambda$.} \label{bbSL}
For a vector space $V\in\Vect$ {\it symmetric algebra} and {\it exterior algebra} are the quadratic algebras
\begin{align}
 &SV=TV/\big(\{v_1\otimes v_2-v_2\otimes v_1\mid v_1,v_2\in V\}\big), \\
 &\Lambda V=TV/\big(\{v_1\otimes v_2+v_2\otimes v_1\mid v_1,v_2\in V\}\big).
\end{align}
We obtain the functors $S\colon\Vect\to\QA$ and $\Lambda\colon\Vect\to\QA$ and their `finite-dimensional' versions $S\colon\FVect\to\FQA$ and $\Lambda\colon\FVect\to\FQA$. The functors $S$ and $\Lambda$ are left adjoint functors to the restrictions of the functor $(-)_1\colon\QA\to\Vect$ to the subcategories $\CommQSAe$ and $\CommQSAo$ respectively. Analogously, the finite-dimensional versions are left adjoints to the restrictions of $(-)_1\colon\FQA\to\FVect$ to the subcategories $\FCommQSAe$ and $\FCommQSAo$.

Denote by $S^kV$ and $\Lambda^kV$ the $k$-th graded components of $S V$ and $\Lambda V$. In these notations we can write the following natural isomorphisms
\begin{align}
 &S^k(V\oplus W)\cong\bigoplus_{l=0}^k(S^lV)\otimes(S^{k-l}W), &
 &\Lambda^k(V\oplus W)\cong\bigoplus_{l=0}^k(\Lambda^lV)\otimes(\Lambda^{k-l}W)
\end{align}
in $\Vect$. Summation over $k\in\NN_0$ gives the natural isomorphisms
\begin{align} \label{Soplus}
 &S(V\oplus W)\cong(SV)\otimes(SW), & &\Lambda(V\oplus W)\cong(\Lambda V)\otimes(\Lambda W)
\end{align}
in $\QA$.
They are the structure morphisms of the {\bf strong} monoidal functors
\begin{align}
 &S\colon(\Vect,\oplus)\to(\QA,\otimes), &
 &S\colon(\FVect,\oplus)\to(\FQA,\otimes), \label{SFunctor} \\
 &\Lambda\colon(\Vect,\oplus)\to(\QA,\ootimes), &
 &\Lambda\colon(\FVect,\oplus)\to(\FQA,\ootimes). \label{LambdaFunctor}
\end{align}
On the other hand the natural transformations
\begin{align}
 &S(V\otimes W)\to(SV)\circ(SW), &
 &(\Lambda V)\bullet(\Lambda W)\to\Lambda(V\otimes W) 
\end{align}
are not isomorphisms even for finite-dimensional $V$ and $W$. This does not give strong monoidal functors. Nevertheless, we obtain {\bf colax} monoidal functors
\begin{align}
 &S\colon(\Vect,\otimes)\to(\QA,\circ), &
 &S\colon(\FVect,\otimes)\to(\FQA,\circ) \label{Sotimes}
\end{align}
and {\bf lax} monoidal functors
\begin{align}
 &\Lambda\colon(\Vect,\otimes)\to(\QA,\bullet), &
 &\Lambda\colon(\FVect,\otimes)\to(\FQA,\bullet).
\end{align}
The functors~$S$ and $\Lambda$ are faithful since their compositions with the functor $(-)_1$ coincide with the identical functor: $(SV)_1=V$, $(\Lambda V)_1=V$. Moreover, any morphism $SV\to SW$ or $\Lambda V\to\Lambda W$ is uniquely given by its first component and hence have the form $Sf$ or, respectively, $\Lambda f$ for some linear map $f\colon V\to W$, so that the functors~$S$ and $\Lambda$ are fully faithful. 

\bb{Koszul duality for quadratic algebras}~\cite{Manin87,Manin88}. Consider a finitely generated quadratic algebra $\A=TV/(R)$, where $V\in\FVect$ and $R\subset V^{\otimes2}$. Define its {\it Koszul dual quadratic algebra} as
\begin{align} \label{AKoszul}
 &\A^!=T(V^*)/(R^\bot), &&R^\bot=\{\xi\in V^*\otimes V^*\mid\xi(r)=0\;\forall\,r\in R\}.
\end{align}
The operation~\eqref{AKoszul} is a contravariant functor $(-)^!\colon\FQA\to\FQA$. We have the following natural isomorphisms
\begin{align}
 &(\A^!)^!\cong\A,
 &&(\A\circ\B)^!\cong\A^!\bullet\B^!,
 &&(\A\bullet\B)^!\cong\A^!\circ\B^!, \label{KoszulDualExch} \\
 &(\A^!)_1=(\A_1)^*,
 &&(\A\otimes\B)^!\cong\A^!\ootimes\B^!,
 &&(\A\ootimes\B)^!\cong\A^!\otimes\B^!.
\end{align}
We see that this functor is an involutive anti-autoequivalence, which interchanges the Manin's operations, and that the diagrams
\begin{align}
\xymatrix{
\FQA\ar[r]^{(-)^!}\ar[d]_{(-)_1} & \FQA^\op\ar[d]^{(-)_1} \\
 \FVect\ar[r]^{(-)^*} & \FVect^\op
}\qquad\qquad
\xymatrix{
\FQA^\op\ar[r]^{(-)^!}\ar[d]_{(-)_1} & \FQA\ar[d]^{(-)_1} \\
 \FVect^\op\ar[r]^{(-)^*} & \FVect
}
\end{align} 
are commutative.

\bb{The functors $S^*$ and $\Lambda^*$}. \label{bbSQVS}
 Define contravariant functors $S^*\colon\FVect\to\FQA$ and $\Lambda^*\colon\FVect\to\FQA$ by the formulae $S^*(V)=S(V^*)$, $\Lambda^*(V)=\Lambda(V^*)$, $V\in\FVect$. As functors $\FVect\to\FQA^\op$ they are the compositions of $(-)^*\colon\FVect\to\FVect^\op$ with the functors $S^\op\colon\FVect^\op\to\FQA^\op$ and $\Lambda^\op\colon\FVect^\op\to\FQA^\op$ respectively, which are the opposite functors to ones defined in~p.~\ref{bbSL}.

The algebra $S^*(V)$ is the algebra of polynomial functions on $V\in\FVect$. By choosing a basis $(x_i)_{i=1}^n$ in $V$ we obtain an isomorphism $V\cong\KK^n$ and hence $S^*(V)\cong S^*(\KK^n)=\KK[x^1,\ldots,x^n]$, where $(x^i)_{i=1}^n$ is the dual basis of $V^*\cong(\KK^n)^*$.

As a fully faithful functor the functor $S^*\colon\FVect\to\FQA^\op$ embeds the finite-dimensional vector spaces into the category $\FQA^\op$. This justifies the term `quantum linear spaces' for the objects of $\FQA^\op$, but this is only finite-dimensional quantum analogue.

The functor $\Lambda^*$ also embeds $\FVect$ into $\FQA^\op$, but in a different way. For $V=\KK^n$ the algebra $\Lambda^*(\KK^n)$ is isomorphic to the Grassmann algebra with generators $\psi^1,\ldots,\psi^n$ and relations $\psi^i\psi^j+\psi^j\psi^i=0$.

Thus, the functors $S^*$ and $\Lambda^*$ are two different embeddings of the category $\FVect$ into $\FQA^\op$. They identify the category $\FVect$ with some full subcategories of the categories $\FCommQSAe^\op\subset\FQA^\op$ and $\FCommQSAo^\op\subset\FQA^\op$ respectively. Hence the embeddings $S^*$ and $\Lambda^*$ can be interpreted as follows. Any vector space $V$ can be considered as a purely even or as a purely odd super-vector space. In the former case we apply $S^*$ to consider the vector space as quantum linear vector space $S^*V=S(V^*)$, in the latter case we apply $\Lambda^*$.

Note that $S^*(V)=S(V^*)=(\Lambda V)^!$ and $\Lambda^*(V)=\Lambda(V^*)=(SV)^!$. This means that the functor $(-)^!$ maps a vector space to its dual with change of its parity. The operations $\A\circ\B$ and $\A\otimes\B$ are analogues of tensor product and direct sum of vector spaces for the case of the quantum linear spaces $\A,\B\in\FQA^\op$ considered as purely even. The operations $\A\bullet\B$ and $\A\ootimes\B$ are also analogues of the tensor product and direct sum, but in this case we need to consider $\A$ and $\B$ as purely odd.

\bb{Monoidal structures of $S^*$.} \label{bbLaxMonS}
 Let us chose the `even' interpretation: consider $\FVect$ as the subcategory of $\CommQSAe^\op\subset\FQA^\op$ by means of $S^*$. Since the monoidal functor~\eqref{SFunctor} can be regarded as strict monoidal, the operation $\oplus$ in $\FVect$ coincides with the operation $\otimes$ in $\FQA$ (see the formula~\eqref{Soplus}). On the other hand, the functor~\eqref{Sotimes} is colax monoidal, so its opposite $S^\op\colon(\FVect^\op,\otimes)\to(\FQA^\op,\circ)$ is lax monoidal. By composing the latter one with the strong monoidal functor $(-)^*\colon(\FVect,\otimes)\to(\FVect^\op,\otimes)$ in a proper order we obtain a lax monoidal functor. Thus we have
\begin{align}
&\text{{\bf strong} monoidal functor} &&S^*\colon(\FVect,\oplus)\to(\FQA^\op,\otimes), \label{strongSs} \\
&\text{{\bf lax} monoidal functor} &&S^*\colon(\FVect,\otimes)\to(\FQA^\op,\circ). \label{laxSs}
\end{align}
In particular, the tensor product of vector spaces does not coincide with the result of the Manin white product after the embedding $S^*$, because the monoidal functor~\eqref{laxSs} is not strong monoidal.

\begin{Rem} \normalfont
 One can define a functor $S^*\colon\Vect\to\QA^\op$ as the composition of the functors $(-)^*\colon\Vect\to\Vect^\op$ and $S^\op\colon\Vect^\op\to\QA^\op$. It has a structure of strong monoidal functor $(\Vect,\oplus)\to(\QA^\op,\otimes)$, but it does not get any monoidal structure $(\Vect,\otimes)\to(\QA^\op,\circ)$, since $(-)^*\colon(\Vect,\otimes)\to(\Vect^\op,\otimes)$ is colax monoidal, while $S^\op\colon(\Vect^\op,\otimes)\to(\QA^\op,\circ)$ is lax monoidal.
\end{Rem}

\bb{Coproduct of quadratic algebras.} \label{bbCop}
 The product of quantum linear spaces is the coproduct in the category $\FQA$ or, more generally, in $\QA$. For any two quadratic algebras $\A=TV/(R)\in\QA$ and $\A'=TV'/(R')\in\QA$ their coproduct exists and has the form
\begin{align}
\A\amalg\A':=T(V\oplus V')/(R\oplus R').
\end{align}
Indeed, let $\B=TW/(S)\in\QA$ and $f\colon\A\to\B$, $f'\colon\A'\to\B$ be two morphisms in $\QA$. Their first order components satisfy $(f_1\otimes f_1)R\subset S$, $(f'_1\otimes f'_1)R'\subset S$. There exists a unique linear map $h_1\colon V\oplus V'\to W$ such that the diagram
\begin{align}
\xymatrix{
 V\ar[r]\ar[rd]_{f_1} & V\oplus V'\ar[d]_{h_1} & V'\ar[l]\ar[ld]^{f'_1} \\
 & W
}
\end{align}
commutes. Since $(h_1\otimes h_1)(R\oplus R')=(f_1\otimes f_1)(R)+(f'_1\otimes f'_1)(R')\subset S$, there exists a morphism $h\colon\A\amalg\A'\to\B$ with the first order component $h_1$. Hence there is a unique morphism $h$ making the diagram
\begin{align}
\xymatrix{
 \A\ar[r]\ar[rd]_{f} & \A\amalg\A'\ar[d]_{h} & \A'\ar[l]\ar[ld]^{f'} \\
 & \B
}
\end{align}
commutative.

Note that the inclusion $R\oplus R'\subset R\oplus[V,V']\oplus R'$ implies the natural epimorphism
\begin{align}
 \A\amalg\A'\twoheadrightarrow\A\otimes\A'.
\end{align}

\begin{Rem} \normalfont
 The coproduct in the category the connected affinely generated graded algebras is different. Namely, the coproduct of $\A=\bigoplus\limits_{k\in\NN_0}\A_k$ and $\A'=\bigoplus\limits_{k\in\NN_0}\A'_k$ in this category is the graded algebra $\A\coprod\limits_{c.a.g.}\A'$ with the components
\begin{align}
 \Big(\A\coprod\limits_{c.a.g.}\A'\Big)_k=\bigoplus_{p=1}^{k+1}\bigoplus_{l_1=0}^k\bigoplus_{l_2,\ldots,l_p=1\atop l_1+l_2+\ldots+l_p=k}^k(\A_{l_1}\otimes\A'_{l_2}\otimes\A_{l_3}\otimes\cdots).
\end{align}
\end{Rem}

\bbo{Functor $T^*$.} \label{bbT}
To extend the tensor product of vector spaces in an exact way one can consider another embedding such as the fully faithful functor $T^*\colon\FVect\to\FQA^\op$ defined as the composition of $(-)^*\colon\FVect\to\FVect^\op$ and $T^\op\colon\FVect^\op\to\FQA^\op$, that is $T^*(V)=T(V^*)$. It is a strong monoidal functor $T^*\colon(\FVect,\otimes)\to(\FQA^\op,\circ)$, but in this case the direct sum $\oplus$ is not strong-monoidally related with the Manin's operation $\otimes$ with quadratic algebras. Indeed, $T(V)\otimes T(W)$ is a quotient algebra of $T(V\oplus W)$ over the ideal generated by the subset $\{v\otimes w-w\otimes v\mid v\in V,\, w\in W\}$, so $T\colon(\FVect,\oplus)\to(\FQA,\otimes)$ is a colax (not strong) monoidal functor.

The functor $T$ preserves the finite coproducts, hence it has a structure of strong monoidal functor $(\Vect,\oplus)\to(\QA,\amalg)$. If $\A,\A'\in\FQA$ then $\A\amalg\A'\in\FQA$, so we also obtain the strong monoidal functor $T\colon(\FVect,\oplus)\to(\FQA,\amalg)$.
Thus we have:
\begin{align}
&\text{{\bf lax} monoidal functor} &&T^*\colon(\FVect,\oplus)\to(\FQA^\op,\otimes), \label{laxTstFVectplus} \\
&\text{{\bf strong} monoidal functor} &&T^*\colon(\FVect,\oplus)\to(\FQA^\op,\amalg), \label{TstFVectplus} \\
&\text{{\bf strong} monoidal functor} &&T^*\colon(\FVect,\otimes)\to(\FQA^\op,\circ), \label{TstFVectotimes} \\
&\text{{\bf strong} monoidal functor} &&T^*\colon(\Vect,\oplus)\to(\QA^\op,\amalg), \\
&\text{{\bf colax} monoidal functor} &&T^*\colon(\Vect,\otimes)\to(\QA^\op,\circ),
\end{align}
where $T^*\colon\Vect\to\QA^\op$ is the composition of the functors $(-)^*\colon\Vect\to\Vect^\op$ and $T^\op\colon\Vect^\op\to\QA^\op$.

\bb{Dequantisation functor.} Composition of the functors $(-)^*$ and $(-)_1$ gives the functor $(-)^*_1\colon\FQA^\op\to\FVect$, $\A\mapsto(\A)_1^*=(\A_1)^*$. By composing it with any of the functors $S^*$, $\Lambda^*$ and $T^*$, which embed $\FVect$ into $\FQA^\op$, we obtain the identical functor:
\begin{align} \label{DeqSLT}
 &\big(S^*(V)\big)_1^*=V, &
 &\big(\Lambda^*(V)\big)_1^*=V, &
 &\big(T^*(V)\big)_1^*=V.
\end{align}
The functor $(-)_1^*\colon\FQA^\op\to\FVect$ is faithful, but not full. It has the following structures of {\bf strong} monoidal functor:
\begin{align}
 &(\FQA^\op,\otimes)\to(\FVect,\oplus), &
 &(\FQA^\op,\circ)\to(\FVect,\otimes), \\
 &(\FQA^\op,\ootimes)\to(\FVect,\oplus), &
 &(\FQA^\op,\bullet)\to(\FVect,\otimes), \\
 &(\FQA^\op,\amalg)\to(\FVect,\oplus).
\end{align}
This functor extracts the classical part from a quadratic algebra. In other words, this is a classical limit for quantum linear spaces.

It can be extended to the functor $(-)^*_1\colon\QA^\op\to\Vect$, which has a lax monoidal structure $(-)_1^*\colon(\QA^\op,\circ)\to(\Vect,\otimes)$, but the identities~\eqref{DeqSLT} do not hold for infinite-dimensional vector spaces $V$.

\subsection{Manin matrices}
\label{secMM}

To define Manin matrices in terms of matrix idempotents we first describe the relationship between quadratic algebras and general idempotent operators.

\bb{Quadratic algebras for idempotents.} \label{bbIdem}
 Remind that an {\it idempotent} is an element $A$ of a ring $\gR$ satisfying $A^2=A$. Let $W$ be a vector space and $R,R'\subset W$ be its subspaces such that $W=R\oplus R'$. The composition $W\twoheadrightarrow W/R'\cong R\hookrightarrow W$ is an idempotent in the algebra $\bfend(W)$. Moreover, the idempotents $A\in\bfend(W)$ are in one-to-one correspondence with the decompositions $W=R\oplus R'$ into subspaces $R,R'\subset W$ such that $R=\Im A$ and $R'=\Im(1-A)$, where $1=\id_{V\otimes V}$. The operator $1-A\in\bfend(W)$ is the {\it dual} idempotent, it corresponds to the decomposition $W=R'\oplus R$.

Let $\A\in\QA$ have the form $\A=TV/(R)$ for a vector space $V\in\Vect$ and a subspace $R\subset V^{\otimes2}$. As any subspace the subspace $R$ can be given as an image of an idempotent $A\in\bfend(V^{\otimes2})$, since $R$ always has a complement $R'\subset V^{\otimes2}$. Thus any connected quadratic algebra is isomorphic to $TV/(\Im A)$ for some idempotent $A\in\bfend(V^{\otimes2})$.

Denote $\Xi_A(\KK)=TV/(R')$ where $R'=\Im(1-A)$. Any $\A\in\QA$ is isomorphic to $\Xi_A(\KK)$ for some idempotent operator $A$. More generally, any semi-connected quadratic algebra is isomorphic to $\Xi_A(\gR):=\gR\otimes\Xi_A(\KK)$ for some $A$ and $\gR\in\Alg$.

If $V\in\FVect$ then $(V^{\otimes2})^*=(V^*)^{\otimes2}$, so any (idempotent) operator $A\colon V^{\otimes2}\to V^{\otimes2}$ gives the (idempotent) transpose operator $A^*\colon(V^*)^{\otimes2}\to(V^*)^{\otimes2}$. Any quadratic algebra $\A\in\FQA$ or $\A\in\FQAsc$ isomorphic to the algebra $\gX_A(\KK):=TV^*/(\Im A^*)=\Xi_{1-A^*}(\KK)$ or $\gX_A(\gR):=\gR\otimes\gX_A(\KK)=\Xi_{1-A^*}(\gR)$ respectively for some idempotent $A\in\bfend(V^{\otimes2})$. The Koszul duality in these notations has the form:
\begin{align} \label{gXKoszul}
 &\gX_A(\KK)^!=\Xi_A(\KK), &&\Xi_A(\KK)^!=\gX_A(\KK).
\end{align}
The semi-connected quadratic algebras $\gX_A(\gR)$ and $\Xi_A(\gR)$ can be considered as values of functors $\gX_A,\Xi_A\colon\Alg\to\FQAsc$ coinciding with $-\otimes\A$ and $-\otimes\A^!$, where $\A=\gX_A(\KK)$.

For the idempotent $A_V=\dfrac{1-\sigma_{V,V}}2$ (anti-symmetrizer of $V\otimes V$) we obtain
\begin{align}
 &\gX_{A_V}(\KK)=S^*(V), &
 &\Xi_{A_V}(\KK)=\Lambda V, &
\end{align}
The functors $T^*$ and $T$ give the quadratic algebras for the trivial idempotents $0,1\in\bfend(V^{\otimes2})$:
\begin{align}
 &\gX_{0}(\KK)=T^*(V), &
 &\Xi_{1}(\KK)=TV.
\end{align}

\bbo{The case of matrix idempotents.}
A basis $(v_i)_{i=1}^n$ of a vector space $V\in\FVect$ gives the isomorphism $V\cong\KK^n$, $v\mapsto(x^1,\ldots,x^n)$, where $x^i$ are components of $v\in V$ defined by the formula $v=\sum_{i=1}^nx^iv_i$. It induces the algebra isomorphism $\bfend(V)\cong\bfend(\KK^n)$, which maps operator on $V$ to its $n\times n$ matrix in the basis $(v_i)$. The space $V\otimes V$ has the basis $(v_i\otimes v_j)_{i,j=1}^n$, so the operators on $V\otimes V$ correspond to the $n^2\times n^2$ matrices, whose entries have two pairs of indices. Namely, the operator $A\in\bfend(V^{\otimes2})$ corresponds to the matrix $A=(A^{ij}_{kl})$ defined by the formula $A(v_k\otimes v_l)=\sum_{i,j=1}^nA^{ij}_{kl}v_i\otimes v_j$. The transpose operator $A^*\in\bfend(V^*\otimes V^*)$ has the form $A^*(v^i\otimes v^j)=\sum_{k,l=1}^nA^{ij}_{kl}v^k\otimes v^l$, where $(v^i)_{i=1}^n$ is the dual basis.

Since any $V\in\FVect$ has a finite basis $(v_i)_{i=1}^n$, any algebra $\A=TV/(R)\in\FQA$ is isomorphic to $\Xi_A(\KK)$ for some matrix idempotent $A\in\bfend(\KK^n\otimes\KK^n)$. The latter is a matrix $A=(A^{ij}_{kl})$ such that $\sum_{k,l=1}^nA^{ij}_{kl}A^{kl}_{ab}=A^{ij}_{ab}$. The algebra $\A=\Xi_A(\KK)$ is generated by $v_1,\ldots,v_n$ with the quadratic commutation relations $\sum_{i,j=1}^nA^{ij}_{kl}v_jv_j=v_kv_l$, where we identified the basis $(v_i)$ with the standard basis of $\KK^n$.

The same algebra is isomorphic to $\gX_A(\KK)=TW/(R)$ for another matrix idempotent $A$ such that $W=(\KK^n)^*$ and $R=\Im A^*$. Let $(v^i)$ be the basis of $(\KK^n)^*$ dual to the standard basis of $\KK^n$, then $\gX_A(\KK)$ is generated by $v^1,\ldots,v^n$ with the relations $\sum_{k,l=1}^n A^{ij}_{kl}v^kv^l=0$.

Consider, for instance, the anti-symmetrizer $A_n=A_{\KK^n}\colon\KK^n\otimes\KK^n\to\KK^n\otimes\KK^n$. It has the entries $(A_n)^{ij}_{kl}=\frac12(\delta^i_k\delta^j_l-\delta^i_l\delta^j_k)$. The corresponding quadratic algebras $\gX_A(\KK)$ and $\Xi_A(\KK)$ are generated by $v^i$ and $v_i$ with the relations $v^iv^j=v^jv^i$ and $v_iv_j=-v_jv_i$. The former one is the polynomial algebra $\KK[x^1,\ldots,x^n]$ with commutative variables $x^i=v^i$, while the latter one is the Grassmann algebra with the anti-commutative variables $v_i$.

\bb{Manin matrix for a pair of idempotents.} \label{bbManinMatr}
Any $n\times m$ matrix $M$ with entries $M_{ij}=M^i_j\in\gR$ can be considered as an operator $\KK^m\to\gR\otimes\KK^n$. In terms of standard bases $(v_i)_{i=1}^n$ and $(w_j)_{j=1}^m$ of the spaces $\KK^n$ and $\KK^m$ it acts as $Mw_j=\sum_{i=1}^nM^i_j\otimes v_i$. Denote by $M^{(a)}$ the operator $M$ acting on the $a$-th matrix tensor factor $\KK^m$ of $\gR\otimes\KK^m\otimes\cdots\otimes\KK^m$. In the case of two such factors we have
\begin{align}
 &M^{(1)}(r\otimes w_j\otimes w_l)=\sum_{i=1}^n M^i_j r\otimes v_i\otimes w_l, &&r\in\gR, \\
 &M^{(2)}(r\otimes w_j\otimes w_l)=\sum_{k=1}^n M^k_l r\otimes w_j\otimes v_k,&&r\in\gR.
\end{align}

\begin{Def} \normalfont \cite{S}
 Let $A\in\bfend(\KK^n\otimes\KK^n)$ and $B\in\bfend(\KK^m\otimes\KK^m)$ be matrix idempotents. The matrix $M$ satisfying 
\begin{align}
 AM^{(1)}M^{(2)}(1-B)=0
\end{align}
 is called {\it Manin matrix for the pair of idempotents $(A,B)$} or simply {\it $(A,B)$-Manin matrix}.
In the case $A=B$ we will say {\it $B$-Manin matrix} instead of $(B,B)$-Manin matrix.
\end{Def}

The matrix $M$ gives the operator $(\KK^n)^*\to\gR\otimes(\KK^m)^*$ acting as $v^i\mapsto\sum_{j=1}^m M^i_j\otimes w^j$, where $(v^i)$ and $(w^j)$ are bases dual to $(v_i)$ and $(w_j)$. Any linear map $\gX_A(\KK)_1\to\gX_B(\gR)_1$ has such form for some matrix $M$. In fact, this map gives a homomorphism of graded algebras $\gX_A(\KK)\to\gX_B(\gR)$ iff $M$ is an $(A,B)$-Manin matrix (see~\cite{S}), we denote this homomorphism of graded algebras by $f_M\colon\gX_A(\KK)\to\gX_B(\gR)$.

In the same way any linear map $\Xi_B(\KK)_1\to\Xi_A(\gR)_1$ has the form $\KK^m\to\gR\otimes\KK^n$, $w_j\mapsto\sum_{i=1}^nM^i_jv_i$. It gives a graded homomorphism $\Xi_B(\KK)\to\Xi_A(\gR)$ iff $M=(M^i_j)$ is an $(A,B)$-Manin matrix, we denote this homomorphism by $f^M\colon\Xi_B(\KK)\to\Xi_A(\gR)$.

Thus we have bijections between the sets $\Hom\big(\gX_A(\KK),\gX_B(\gR)\big)$, $\Hom\big(\Xi_B(\KK),\Xi_A(\gR)\big)$ and the set of $(A,B)$-Manin matrices.

\bb{Usual Manin matrices.} The Manin matrices introduced in the papers~\cite{CF,CFR} are $n\times m$ Manin matrices for the pair of idempotents $(A_n,A_m)$. They correspond to the graded homomorphisms of polynomial algebras $\KK[x^1,\ldots,x^n]\to\gR[y^1,\ldots,y^m]$ or, equivalently, to the graded homomorphisms $S(V^*)\to\gR\otimes S(W^*)$, where $V$ and $W$ are $n$- and $m$-dimensional vector spaces.

\subsection{Semi-linear spaces and their quantum analogue}
\label{secQSLS}

To interpret representations of algebraic monoids and groups on (finite-dimensional) linear spaces in terms of Subsection~\ref{secRepMon} one needs to include both algebraic sets and linear spaces to a bigger monoidal category. More generally, one needs to do the same for the affine schemes. In this way we will get categories of spaces which are partially linear in some sense. We call them semi-linear spaces.

\bb{Semi-linear algebraic sets.}
 Describe the first case of semi-linear spaces. For an algebraic set $X\in\AlgSet$ and a linear space $V\in\FVect$ consider their set-theoretic product $X\times V$ as an algebraic set. Define morphisms between two such products $X\times V$ and $Y\times W$ as morphisms $F\colon X\times V\to Y\times W$ in $\AlgSet$ that have the form
\begin{align} \label{Fxv}
 F(x,v)&=(\varphi(x),f(x,v)),\qquad\quad x\in X,\quad v\in V,\qquad\qquad\varphi\colon X\to Y,\qquad f\colon X\times V\to W, \notag \\
 &f(x,\alpha v+v')=\alpha f(x,v)+f(x,v'),\quad\quad x\in X,\quad v,v'\in V.
\end{align}
In other words, the morphism $F$ is given by arbitrary morphisms $\varphi\in\Hom_{\AlgSet}(X,Y)$ and $f\in\Hom_{\AlgSet}(X\times V,W)$ such that $f(x,-)\colon V\to W$ is linear for any $x\in X$. In this way we obtain a category $\SLAlgSet$ with morphisms that partially satisfy the linearity condition, we call them {\it semi-linear maps}. We call the objects $X\times V\in\SLAlgSet$ {\it semi-linear algebraic sets}.

\bb{Monoidal product of semi-linear algebraic sets.} Let $X,X'\in\AlgSet$ and $V,V'\in\FVect$. Define monoidal product of $X\times V,X'\times V'\in\SLAlgSet$ as the object $(X\times X')\times(V\otimes V')$. For morphisms $F\colon X\times V\to Y\times W$ and $F'\colon X'\times V'\to Y'\times W'$ in $\SLAlgSet$, given by maps $\varphi\colon X\to Y$, $f\colon X\times V\to W$ and $\varphi'\colon X'\to Y'$, $f'\colon X'\times V'\to W'$, their monoidal product is $F''=F\otimes F'\colon(X\times X')\times(V\otimes V')\to (Y\times Y')\times(W\otimes W')$ given by $\varphi''=\varphi\times\varphi'\colon X\times X\to Y\times Y$ and $f''\colon(X\times X')\times(V\otimes V')\to W\otimes W'$, $f''(x,x',v\otimes v')=f(x,v)\otimes f'(x',v')$. In this way we obtain a symmetric monoidal category $(\SLAlgSet,\otimes)$ with the unit object $0\times\KK$.

The categories $(\AlgSet,\times)$ and $(\FVect,\otimes)$ can be regarded as monoidal subcategories of $(\SLAlgSet,\otimes)$ embedded via the faithful symmetric strict monoidal functors $X\mapsto X\times\KK$ and $V\mapsto\{0\}\times V$. The latter functor is full, the former is not, i.e. $\FVect$ is a full subcategory of $\SLAlgSet$, while the subcategory $\AlgSet\subset\SLAlgSet$ is not full.

\bb{Trivial vector bundles.} \label{bbTVB}
The category $\SLAlgSet$ is equivalent to the category of trivial vector bundles over algebraic sets. The semi-linear algebraic set $X\times V$ corresponds to the trivial bundle over $X$ with the fibre $V$. A semi-linear map $F\colon X\times V\to Y\times W$ given by the formula~\eqref{Fxv} is the vector bundle morphism $(\varphi,F)$ in the sense that the diagram
\begin{align}
 \xymatrix{
X\times V\ar[d]\ar[r]^F & Y\times W\ar[d] \\
X\ar[r]^\varphi & Y
}
\end{align}
is commutative (the vertical arrows mean the projections) and the morphism $F$ induces a linear map on each fibre (see~\cite[III,~\S~1]{Lang}).
Category embeddings $\AlgSet\hookrightarrow\SLAlgSet$ and $\FVect\hookrightarrow\SLAlgSet$ identify an algebraic set $X$ with the $1$-dimensional trivial bundle over it, and a vector space $V$ -- with the trivial bundle over $\{0\}=\Spec\KK$ with the fibre $V$.

The main difference between the notions of semi-linear algebraic set and trivial vector bundle is that we consider the former ones not over the structure of algebraic set, but together with this structure. In particular, the monoidal product of two semi-linear algebraic sets $X\times V$ and $X\times W$ is not the same as the tensor product of the corresponding bundles over $X$. One can consider more general vector bundles in the same manner, but the category $\SLAffSch$ is enough for out purposes, this is a minimal monoidal category that includes $(\AlgSet,\times)$ and $(\Vect,\otimes)$.

\bb{Algebras of functions on semi-linear algebraic sets.} Let $X\in\AlgSet$ and $V\in\FVect$. The algebra of regular functions on $X\times V$ as on algebraic set is $A(X)\otimes SV^*$. It inherits the grading from $SV^*$ as a tensor product of an algebra with a graded algebra (see p.~\ref{bbGr}).
 Let us denote $S^*(X\times V)=A(X)\otimes SV^*\in\GrAlg$. For a morphism of the form~\eqref{Fxv} we obtain the morphism of graded algebras with the reverse direction: $S^*(F)\colon A(Y)\otimes SW^*\to A(X)\otimes SV^*$, it is defined as
\begin{align}
 &S^*(F)\colon s\mapsto\varphi^*(s)=s\cdot\varphi\in A(X), &
 &\varphi^*(s)(x)=s\big(\varphi(x)\big), \label{SFvarphi} \\
 &S^*(F)\colon\xi\mapsto f^*(\xi)=\xi\cdot f\in A(X)\otimes V^*, &
  &f^*(\xi)(x,v)=\xi\big(f(x,v)\big), \label{SFf}
\end{align}
where $s\in A(Y)$, $\xi\in W^*$, $x\in X$, $v\in V$. The formula~\eqref{SFvarphi} gives the usual embedding $\AlgSet\hookrightarrow\CommAlg^\op$. Let $(v_j)$, $(v^j)$ and $(w_i)$, $(w^i)$ be dual bases of $V$, $V^*$ and $W$, $W^*$ respectively. Then $f(x,v_j)=\sum_i f^i_j(x)w_i$ for some $f^i_j\in A(X)$, so the map~\eqref{SFf} has the form $f^*(w^i)=\sum_j f^i_j\otimes v^j\in A(X)\otimes V^*$.

We obtain a lax monoidal functor $S^*\colon(\SLAlgSet,\otimes)\to(\FQAsc^\op,\circ)$. In particular, any algebraic set $X$ as the object $X\times\KK$ corresponds to the graded algebra $S^*(X)=A(X)\otimes\KK[u]$ and we have the natural isomorphism $S^*(X\times V)\cong S^*(X)\circ S^*(V)$.

\bb{Semi-linear affine schemes.} \label{bbSLAS}
Let us generalise the semi-linear spaces to the case of affine schemes.
Among the commutative semi-connected quadratic algebras consider the algebras of the form $\gR\otimes SV^*$ where $\gR\in\CommAlg$, $V\in\FVect$. Such algebra can be interpreted as the `function algebra' on the affine scheme $X\times V$, where $X=\Spec\gR$. This is an affine scheme which corresponds to the algebra $\gR\otimes SV^*$ and has an additional structure given by the grading of this algebra. We call such schemes {\it semi-linear affine schemes}. By definition they are objects of the category $\FQAsc^\op$ corresponding to the graded algebras $\gR\otimes SV^*$, where $\gR\in\CommAlg$, $V\in\FVect$. The embedding $\SLAffSch\hookrightarrow\FQAsc^\op$ extends the functor $S^*\colon\SLAlgSet\to\FQAsc^\op$, $S^*(X\times V)=\gR\otimes SV^*$, since $\SLAlgSet$ is a subcategory of $\SLAffSch$, so we also denote it by $S^*$.

Let $\gS\in\CommAlg$, $Y=\Spec\gS$ and $W\in\FVect$.
A morphism $X\times V\to Y\times W$ in $\SLAffSch$ is given by a graded homomorphism $\gS\otimes SW^*\to\gR\otimes SV^*$ generated by an algebra homomorphism $\alpha\colon\gS\to\gR$ and a linear map $t\colon W^*\to\gR\otimes V^*$. In the order $1$ it has the form
\begin{align}
 &(\alpha,t)\colon\gS\otimes W^*\to\gR\otimes V^*, &&(\alpha,t)(s\otimes\xi)=\alpha(s)t(\xi). 
\end{align}
 For a morphism in $\SLAlgSet$ of the form~\eqref{Fxv} we have $\alpha=\varphi^*$ and $t=f^*$, $(\alpha,t)=F^*$.

\begin{Rem} \normalfont
The semi-linear affine scheme $X\times V$ can be identified with a quasi-coherent sheaf corresponding to the free $\gR$-module $\gR\otimes V$ (or rather to its dual $\gR\otimes V^*$). The category $\SLAffSch$ is equivalent to the category of such free sheaves ($\gR$-modules) with appropriate morphisms. This generalises the equivalence described in p.~\ref{bbTVB}.
\end{Rem}

Let $\gR,\gR'\in\CommAlg$, $X=\Spec\gR$, $X'=\Spec\gR'$, $V,V'\in\FVect$. Define the monoidal product of the semi-linear affine schemes $X\times V$ and $X'\times V'$ as the semi-linear affine scheme $(X\times X')\times(V\otimes V')$ with the function algebra $\gR\otimes\gR'\otimes S^*(V\otimes V')$. Note that the embedding $\SLAffSch\hookrightarrow\FQAsc^\op$ is a fully faithful symmetric lax monoidal functor $S^*\colon(\SLAffSch,\otimes)\to(\FQAsc^\op,\circ)$.

\bb{Semi-linear algebraic monoids and semi-linear affine monoid schemes.} \label{bbSLAM}
A monoid in $(\SLAlgSet,\otimes)$ is a semi-linear algebraic set $X\times V\in\SLAlgSet$ with the semi-linear maps $\mu_{X\times V}\colon(X\times X)\times(V\otimes V)\to X\times V$, $\eta_{X\times V}\colon\KK\to X\times V$ which have the form
\begin{align}
 &\mu_{X\times V}(x,y,u\otimes v)=\big(xy,f(x,y,u\otimes v)\big), &\eta_{X\times V}(1)=\big(e,1_V\big),
\end{align}
where $xy=\mu_X(x,y)$, $e=\eta_X(0)$ are defined by a monoid $(X,\mu_X,\eta_X)\in\Mon(\AlgSet,\times)$ and $f\colon X\times X\times(V\otimes V)\to V$ is a morphism in $\AlgSet$ linear in the last argument such that $f(e,x,1_V\otimes v)=v=f(x,e,v\otimes1_V)$, $f\big(x,yz,u\otimes f(y,z,v\otimes w)\big)=f\big(xy,z,f(x,y,u\otimes v)\otimes w\big)$, $x,y,z\in X$, $u,v,w\in V$. Denote by $V_x$ the copy of $V$ considered as a fibre over $x\in X$, then for any $x,y\in X$ we obtain the linear map $f_{x,y}\colon V_x\otimes V_y\to V_{xy}$, $f_{x,y}(u\otimes v)=f(x,y,u\otimes v)$, and the last condition is equivalent to the commutativity of the diagram
\begin{align}
\xymatrix{
 V_x\otimes V_y\otimes V_z\ar[rr]^{\quad\id\otimes f_{y,z}}\ar[d]_{f_{x,y}\otimes\id} &&
 V_x\otimes V_{yz}\ar[d]^{f_{x,yz}} \\
 V_{xy}\otimes V_z\ar[rr]^{\quad f_{xy,z}} && V_{xyz}
}
\end{align}
while $1_V$ belongs to $V_e$. If $X=\{e\}$ is a one-point set, then to equip $V=X\times V$ with a monoid structure is the same as to choose the algebra $(V,f_{e,e},1)\in\Mon(\Vect,\otimes)$. If $V=\KK$, then a monoid structure on $X\times V=X\times\KK$ is given by an algebraic monoid $(X,\mu_X,e)\in\Mon(\AlgSet,\times)$ with collection of coefficients $p_{x,y}\in\KK$, while $q\in\KK\backslash\{0\}$ satisfying $p_{x,yz}p_{y,z}=p_{xy,z}p_{x,y}$ and $p_{e,x}=p_{x,e}=q^{-1}$ (the linear map $f_{x,y}\colon\KK\otimes\KK\to\KK$ is a multiplication by the coefficient $p_{x,y}$, while $1_V=q$). The case of a usual algebraic monoid corresponds to the values $p_{x,y}=q=1$. In this way we see that the categorical embeddings $\FVect\hookrightarrow\SLAlgSet$ and $\AlgSet\hookrightarrow\SLAlgSet$ induce the fully faithful functor $\Mon(\FVect,\otimes)\to\Mon(\SLAlgSet,\otimes)$ and the faithful functor $\Mon(\AlgSet,\times)\to\Mon(\SLAlgSet,\otimes)$, which are symmetric strong monoidal due to p.~\ref{bbMonF}.

Let $\gR=A(X)$. The homomorphisms $\alpha=\mu_X^*\colon\gR\to\gR\otimes\gR$ and $\beta=\eta_X^*\colon\gR\to\KK$ gives the structure of coalgebra on $\gR$. The linear maps $t=f^*\colon V^*\to\gR\otimes\gR\otimes V^*\otimes V^*$ and $\varepsilon_{V^*}=1_V\colon V^*\to\KK$ make the diagrams
\begin{gather} \label{DiagAlphat}
\xymatrix{
 \gR\otimes V^*\ar[rr]^{(\alpha,t)}\ar[d]^{(\alpha,t)} &&
 \gR\otimes\gR\otimes V^*\otimes V^*\ar[d]^{\sigma^{(23)}} \\
\gR\otimes\gR\otimes V^*\otimes V^*\ar[d]^{\sigma^{(23)}} && \gR\otimes V^*\otimes\gR\otimes V^*\ar[d]^{\id\otimes\id\otimes(\alpha,t)} \\
 \gR\otimes V^*\otimes\gR\otimes V^*\ar[d]^{(\alpha,t)\otimes\id\otimes\id} &&
\gR\otimes V^*\otimes\gR\otimes\gR\otimes V^*\otimes V^*\ar[d]^{\sigma^{(34)}\sigma^{(23)}} \\
 \gR\otimes\gR\otimes V^*\otimes V^*\otimes\gR\otimes V^*\ar[rr]^{\sigma^{(34)}\sigma^{(45)}} &&
\gR\otimes\gR\otimes\gR\otimes V^*\otimes V^*\otimes V^*
} \\
\xymatrix{
 &&& \gR\otimes V^*\ar[d]^{(\alpha,t)}\ar@{=}[dlll]\ar@{=}[drrr] &&& \\
 \KK\otimes\gR\otimes\KK\otimes V^* &&&
 \gR\otimes\gR\otimes V^*\otimes V^*
\ar[lll]_{\beta\otimes\id\otimes\varepsilon_{V^*}\otimes\id}
\ar[rrr]^{\id\otimes\beta\otimes\id\otimes\varepsilon_{V^*}} &&&
 \gR\otimes\KK\otimes V^*\otimes\KK
} \notag
\end{gather}
to be commutative, where $\sigma^{(i,i+1)}=\id^{\otimes(i-1)}\otimes \sigma\otimes\id\otimes\cdots\otimes\id$. Since $(\gR,\alpha,\beta)$ is a comonoid in $(\CommAlg,\otimes)$, it is enough to check the commutativity of the diagrams~\eqref{DiagAlphat} by substituting the elements of the form $1\otimes\xi\in\gR\otimes V^*$.

Now consider monoids in $(\SLAffSch,\otimes)$, which we call {\it semi-linear affine monoid schemes}. Let $\gR\in\CommAlg$, $V\in\FVect$ and $X=\Spec\gR$. A structure of monoid on $X\times V$ is defined by the morphisms $\mu_{X\times V}\colon(X\times V)\otimes(X\times V)\to X\times V $ and $\eta_{X\times V}\colon\KK\to X\times V$, which are given by graded homomorphisms $\gR\otimes S(V^*)\to\gR\otimes\gR\otimes S(V^*)\otimes S(V^*)$ and $\gR\otimes S(V^*)\to\KK[u]$. The zero and first order components of these homomorphisms are $\alpha\colon\gR\to\gR\otimes\gR$, $(\alpha,t)\colon\gR\otimes V^*\to\gR\otimes\gR\otimes V^*\otimes V^*$ and $\beta\colon\gR\to\KK$, $(\beta,\varepsilon_{V^*})\colon\gR\otimes V^*\to\KK$, where $t\colon V^*\to\gR\otimes\gR\otimes V^*\otimes V^*$ and $\varepsilon_{V^*}\colon V^*\to\KK$ are linear maps. They give a structure of a monoid on $X\times V$ iff $(\gR,\alpha,\beta)\in\Comon(\CommAlg,\otimes)$ and the diagrams~\eqref{DiagAlphat} commute.
The embedding $\AffSch\hookrightarrow\SLAffSch$ induces the faithful symmetric strong monoidal functor $\Mon(\AffSch,\times)\to\Mon(\SLAffSch,\otimes)$.

\bb{Quantum semi-linear spaces.} \label{bbQSLS}
The `algebra of functions' on a semi-linear space $X\times V$ is the tensor product of a quadratic algebra $S^*V$ with a commutative algebra that is the `algebra of functions' on $X$.  By following the Manin's idea described in the beginning of Section~\ref{sec3} we call {\it quantum semi-linear space} any semi-connected quadratic algebra $\gR\otimes\A$ (where $\gR\in\Alg$, $\A\in\QA$) considered as an object of the dual category $\QAsc^\op$. The finite-dimensional case corresponds to the objects of $\FQAsc^\op$. The quantum linear spaces correspond to the connected case $\gR=\KK$. In Section~\ref{sec5} we will construct Quantum Representation Theory as the theory of representations of monoids in $\bfC=(\QAsc^\op,\circ)$ on the objects of a monoidal subcategory $\bfP=\FQA^\op$.

The composition of the embeddings $\AffSch\hookrightarrow\SLAffSch\hookrightarrow\FQAsc^\op$ is a faithful symmetric strong monoidal functor $(\AffSch,\times)\to(\FQAsc^\op,\circ)$, its dual is the functor $(\CommAlg,\otimes)\to(\FQAsc,\circ)$, which maps $\gR$ to $\gR\otimes\KK[u]$. Thus, we consider an affine non-commutative space $\gR\in\Alg^\op$ as the quantum semi-linear space $\gR\otimes\KK[u]\in\FQA^\op$.

Any morphism $f\colon\B\to\A$ in $\QAsc$ between arbitrary semi-connected quadratic algebras $\B=\gS\otimes\big(TW/(S)\big)$ and $\A=\gR\otimes\big(TV/(R)\big)$ is uniquely defined by its zero and first order components $\alpha=f_0\colon\gS\to\gR$ and $f_1\colon\gS\otimes W\to\gR\otimes V$.
We have $f_1=(\alpha,t)$ in the sense that $f_1(s\otimes w)=\alpha(s)t(w)=t(w)\alpha(s)$, where $t\colon W\to\gR\otimes V$ is defined by the formula $t(w)=f_1(1\otimes w)$. Since $\A_2$ is the quotient of $(\gR\otimes V)\otimes_\gR(\gR\otimes V)=\gR\otimes V\otimes V$ over $\gR\otimes R$, the linear map $f_2\colon\B_2\to\A_2$ is induced by
\begin{align*}
 &\gS\otimes W\otimes W=(\gS\otimes W)\otimes_\gS(\gS\otimes W) \xrightarrow{f_1\otimes_\gS f_1}
(\gR\otimes V)\otimes_\gR(\gR\otimes V)=\gR\otimes V\otimes V, \\
 &s\otimes w\otimes w'=(s\otimes w)\otimes_{\gS}(1\otimes w')\mapsto 
\big(\alpha(s)t(w)\big)\otimes_{\gR}t(w')=\alpha(s)\big(t(w)\otimes_\gR t(w')\big), \\
&1\otimes w\otimes w'\mapsto\big(t(w)\otimes_\gR t(w')\big)=(t\dototimes t)(w\otimes w'),
\end{align*}
where $t\dototimes t\colon W\otimes W\to\gR\otimes V\otimes V$ is the composition
\begin{align*}
W\otimes W\xrightarrow{t\otimes t}
\gR\otimes V\otimes\gR\otimes V\xrightarrow{\id_\gR\otimes\sigma_{V,\gR}\otimes\id_V}
\gR\otimes\gR\otimes V\otimes V\xrightarrow{\mu_\gR\otimes\id_V\otimes\id_V}
 \gR\otimes V\otimes V.
\end{align*}
An algebra homomorphism $\alpha\colon\gS\to\gR$ and a linear map $t\colon W\to\gR\otimes V$ define a homomorphism $f\colon\B\to\A$ iff
\begin{align} \label{alphatcomm}
 &\alpha(s)t(w)=t(w)\alpha(s)\qquad\forall\,s\in\gS,\,w\in W
\end{align}
and $(t\dototimes t)S\subset\gR\otimes R$. In particular, $(\B,\Delta,\varepsilon)$ is a comonoid in $(\QAsc,\circ)$ iff the zero and first order components of the graded homomorphisms $\Delta\colon\B\to\B\circ\B$ and $\varepsilon\colon\B\to\KK[u]$ satisfy the following conditions: $\alpha=\Delta_0\colon\gS\to\gS\otimes\gS$ and $\beta=\varepsilon_0\colon\gS\to\KK$ are algebra homomorphisms, $\Delta_1(s\otimes w)=\alpha(s)t(w)$ and $\varepsilon_1(s\otimes w)=\varepsilon_0(s)\varepsilon_W(w)$ for some linear maps $t\colon W\to\gS\otimes\gS\otimes W\otimes W$, $\varepsilon_W\colon W\to\KK$ subjected to the commutativity condition~\eqref{alphatcomm}, the inclusion conditions $(t\dototimes t)S\subset\gS\otimes\gS\otimes\sigma^{(23)}(S\otimes W^{\otimes 2}+W^{\otimes 2}\otimes S)$, $(\varepsilon_W\otimes\varepsilon_W)S=0$ and making the diagrams~\eqref{DiagAlphat} to be commutative (where $V^*=W$ and $\gR=\gS$).

A coaction of the comonoid $\OO=(\A,\Delta,\varepsilon)\in\Comon(\QAsc,\circ)$ on an object $\B\in\QAsc$ is a graded homomorphism $\delta\colon\B\to\A\circ\B$. Its zero component $\delta_0\colon\B_0\to\A_0\otimes\B_0$ should be coaction of $(\A_0,\Delta_0,\varepsilon_0)\in\Comon(\Alg,\otimes)$ on $\B_0\in\Alg$ and the first order component is $\delta_1=(\delta_0,p)$ for some linear map $p\colon W\to\A_0\otimes\B_0\otimes V\otimes W$. To obtain the condition for $p$ one needs to take first order components of the morphisms in the commutative diagrams~\eqref{deltaDiag}, this leads to diagrams similar to~\eqref{DiagAlphat}.

\subsection{Internal hom for quantum linear spaces}

The next ingredient, which we need to construct quantum representations, is a quantum analogue of the internal hom of vector spaces. This is an internal hom in $(\FQA^\op,\circ)$, i.e. an internal cohom in $(\FQA,\circ)$. We need also its generalised version in the sense of p.~\ref{bbGenerIntHom}.

\bb{Internal cohom and universal Manin matrix.} \label{bbIntHomManin}
It was shown in~\cite{Manin87,Manin88} that the category $\FQA$ is coclosed with respect to the Manin product, namely, proved the existence of a natural isomorphism $\Hom(\A\bullet\B^!,\C)\cong\Hom(\A,\B\circ\C)$. Hence the internal cohom for $\A,\B\in\FQA$ can be defined as%
\footnote{Manin denoted the cohom-functor $\bfcohom\colon\FQA^\op\times\FQA\to\FQA$ as $\underline{\hom}$.
}
\begin{align}
 \bfcohom(\B,\A)=\B^!\bullet\A. \label{cohomBA}
\end{align}
Let us calculate it for $\A=\gX_A(\KK)$ and $\B=\gX_B(\KK)$, where $A\in\bfend(V^{\otimes2})$ and $B\in\bfend(W^{\otimes2})$ are idempotents and $V=\KK^n$, $W=\KK^m$. Since $\B^!=\Xi_B(\KK)=TW/\big(\Im(1-B)\big)$ and $\A=TV^*/(\Im A^*)$ we have
\begin{align}
 \B^!\bullet\A=T(W\otimes V^*)/\big(\sigma^{(23)}(\Im(1-B)\otimes\Im A^*)\big),
\end{align}
where $\sigma^{(23)}=\id\otimes\sigma\otimes\id$ (the indices $2,3$ mean that $\sigma$ acts on the second and third tensor factors). Let $(v_i)$ and $(w_j)$ be bases of $V$ and $W$ respectively. Let $(v^i)$ be basis of $V^*$ dual to $(v_i)$. Then the basis of $W\otimes V^*=\Hom(V,W)$ consists of $\M^i_j=w_j\otimes v^i$. The subspace $\sigma^{(23)}(\Im(1-B)\otimes\Im A^*)\subset W\otimes V^*\otimes W\otimes V^*$ is spanned by the vectors
 $$\sigma^{(23)}\sum_{k,l,s,t}(1-B)^{kl}_{ij}A^{pq}_{st}w_k\otimes w_l\otimes v^s\otimes v^t=\sum_{k,l,s,t}A^{pq}_{st}(\M_k^s\otimes\M_l^t)(1-B)^{kl}_{ij},$$
where $(1-B)^{kl}_{ij}=\delta^k_i\delta^l_j-B^{kl}_{ij}$. Therefore, the algebra $\bfcohom(\B,\A)$ is generated by $\M^i_j$ with the quadratic relation
\begin{align}
 A\M^{(1)}\M^{(2)}(1-B)=0, \label{Muniv}
\end{align}
where $\M$ is the matrix with entries $\M^i_j$. The relation~\eqref{Muniv} means exactly that $\M$ is an $(A,B)$-Manin matrix over $\gR=\bfcohom\big(\gX_B(\KK),\gX_A(\KK)\big)$. We call it {\it universal $(A,B)$-Manin matrix}. Note also that
\begin{multline}
 \bfcohom\big(\gX_B(\KK),\gX_A(\KK)\big)=\gX_B(\KK)^!\bullet\gX_A(\KK)=\Xi_A(\KK)^!\bullet\Xi_B(\KK)=\\
 =\bfcohom\big(\Xi_A(\KK),\Xi_B(\KK)\big). \label{cohomgXAB}
\end{multline}

\bbo{Connection with the Manin matrices.} In~\cite{S} we denote by $\gU_{A,B}$ the algebra~\eqref{cohomgXAB} considered as an object of $\Alg$. This algebra was interpreted there via left adjoints to the functors $\gX_B\colon\Alg\to\FQAsc$ and $\Xi_A\colon\Alg\to\FQAsc$. In particular, we derived the natural bijections
\begin{align}
 \Hom_\GrAlg\big(\gX_A(\KK),\gX_B(\gR)\big)\cong\Hom_\Alg(\gU_{A,B},\gR)\cong\Hom_\GrAlg\big(\Xi_B(\KK),\Xi_A(\gR)\big), \label{SetABManin}
\end{align}
which can be described as follows.
The left and right sets consists of the graded homomorphisms $f_M$ and $f^M$ for $(A,B)$-Manin matrices $M$ over $\gR$ (see p.~\ref{bbManinMatr}). A formula $\M^i_j\mapsto M^i_j$ gives a homomorphism $\gU_{A,B}\to\gR$ (in $\Alg$) iff $M=(M^i_j)$ is an $(A,B)$-Manin matrix.

Consider the case when $\gR$ is also a quadratic algebra: $\gR=\gX_C(\KK)$ for some idempotent $C$. Note that $\gX_B\big(\gX_C(\KK)\big)_k=\gX_C(\KK)\otimes\gX_B(\KK)_k$, so the algebra $\gX_B\big(\gX_C(\KK)\big)$ does not coincide with $\gX_C(\KK)\otimes\gX_B(\KK)$ as graded algebra in the sense of the Manin's operation $\otimes$ (see p.~\ref{bbManinOpns}). Consider the subset of the set~\eqref{SetABManin} corresponding to the $(A,B)$-Manin matrices $M$ with entries of the first order: $M^i_j\in\gX_C(\KK)_1$. As a subset of $\Hom_\GrAlg\Big(\gX_A(\KK),\gX_B\big(\gX_C(\KK)\big)\Big)$ it consists of the graded homomorphisms $\gX_A(\KK)\to\gX_B\big(\gX_C(\KK)\big)$ which factors through the morphism $\gX_C(\KK)\circ\gX_B(\KK)\to\gX_B\big(\gX_C(\KK)\big)$. This subset can be identified with the set $\Hom_\GrAlg\big(\gX_A(\KK),\gX_C(\KK)\circ\gX_B(\KK)\big)$. As a subset of $\Hom_\Alg\Big(\gU_{A,B},\gX_B\big(\gX_C(\KK)\big)\Big)$ it consists of homomorphisms $\gU_{A,B}\to\gX_C(\KK)$ which preserve the grading, where the grading of the algebra $\gU_{A,B}$ arises, if one considers it as the cohom-object $\bfcohom\big(\gX_B(\KK),\gX_A(\KK)\big)$. This leads to the adjunction from the beginning of p.~\ref{bbIntHomManin}, in the notations of this paragraph it takes the form
\begin{align}
 \Hom\big(\gX_A(\KK),\gX_C(\KK)\circ\gX_B(\KK)\big)\cong\Hom\Big(\bfcohom\big(\gX_B(\KK),\gX_A(\KK)\big),\gX_C(\KK)\Big),
\end{align}
this isomorphism is induced by the isomorphism~\eqref{SetABManin}.

\bb{Internal hom for finite-dimensional vector spaces.} Let us regard how the internal hom of quantum linear spaces agrees with the internal hom of vector spaces. Remind that the internal hom-object $\bfhom(W,V)$ in $\Vect$ is the set $\Hom(W,V)$ equipped with the usual structure of vector space (see p.~\ref{bbExIntHom}). In the case $W\in\FVect$ we have natural linear isomorphism $\bfhom(W,V)=W^*\otimes V$. By applying the strong monoidal functor $(-)^*_1\colon(\FQA^\op,\circ)\to(\FVect,\otimes)$ to the internal hom~\eqref{cohomBA} we obtain exactly $\big(\bfcohom(\B,\A)\big)^*_1=W^*\otimes V=\bfhom(W,V)$, where $W=\B_1^*$, $V=\A_1^*$. In the notations of Subsection~\ref{secTransRep} it means that the morphisms $\Phi_{\B,\A}$ for the strong monoidal functor $(-)^*_1\colon(\FQA^\op,\circ)\to(\FVect,\otimes)$ are isomorphisms.

On the other hand the quadratic algebra $S^*\bfhom(W,V)$ does not coincide with the algebra $\bfcohom(S^*W,S^*V)=(S^*W)^!\bullet S^* V=\Lambda(W)\bullet S^*(V)$, which means that $\Phi_{W,V}$ for the lax monoidal functor $S^*\colon(\FVect,\otimes)\to(\FQA^\op,\circ)$ are not isomorphisms.

\bb{Coevaluation and cocomposition morphisms.}
 Let $(v_i)$, $(w_j)$ be the bases of $V,W\in\FVect$ and let $(v^i)$, $(w^j)$ be the dual bases of $V^*$, $W^*$. The elements $v^i$, $w^j$ and $\M^i_j=w_j\otimes v^i$ are generators of the algebras $\A=TV^*/(R)$, $\B=TW^*/(S)$ and $\bfcohom(\B,\A)$. In these terms the coevaluation~\eqref{coevXY} (for $X=\A$, $Y=\B$) has the form (see~\cite{Manin88})
\begin{align} \label{coevAB}
 &\coev\colon \A\to\bfcohom\big(\B,\A\big)\circ\B, & &v^i\mapsto\sum_j\M^i_j\otimes w^j.
\end{align}
 Let $\C\in\FQA$ and $(z_l)$ be a basis of $Z=\C_1^*$. Then $\N^j_l=z_l\otimes w^j$ and $\K^i_l=z_l\otimes v^i$ are the generators of the algebras $\bfcohom(\C,\B)$ and $\bfcohom(\C,\A)$. The cocomposition morphism~\eqref{dXYZ}  (for $X=\A$, $Y=\B$, $Z=\C$) reads (see~\cite{S})
\begin{align} \label{dKMN}
 &d=d_{\A,\B,\C}\colon\bfcohom(\C,\A)\to\bfcohom(\B,\A)\circ\bfcohom(\C,\B), &&\K^i_l\mapsto\sum_{j=1}^m \M^i_j\otimes\N^j_l.
\end{align}

Let $B=(B^{ij}_{kl})_{i,j,k,l=1}^m$ be a matrix idempotent such that the elements $\sum\limits_{i,j=1}^m B^{ij}_{kl}w_i\otimes w_j$ span the subspace $S\subset W\otimes W$, then $\B=\gX_B(\KK)$ and $\bfcohom(\B,\B)$ is the quadratic algebra generated by $\M^i_j=w_j\otimes w^j$ with the relations $B\M^{(1)}\M^{(2)}(1-B)=0$, where $\M=(\M^i_j)_{i,j=1}^m$ is the universal $B$-Manin matrix. The comonoid $\bfcoend(\B)\in\Comon(\FQA,\circ)$ is the quadratic algebra $\bfcohom(\B,\B)$ with the structure morphisms $d_\B\colon\bfcoend(\B)\to\bfcoend(\B)\circ\bfcoend(\B)$ and $v_\B\colon\bfcoend(\B)\to\KK[u]$ defined on generators as
\begin{align} \label{dvMMM}
 &d_\B(\M^i_j)=\sum_{k=1}^m\M^i_k\otimes\M^k_j, &&v_\B(\M^i_j)=\delta^i_j u.
\end{align}

Since $\B^!=\Xi_B(\KK)=TW/\big(\Im(1-B)\big)=\gX_{1-B^\top}(\KK)$ is generated by $w_1,\ldots,w_m\in W$, the algebra $\bfcoend(\B^!)$ is generated by the elements $\wt\M_{ij}=w^j\otimes w_i\in W^*\otimes W$. The latter ones are entries of the universal $(1-B^\top)$-Manin matrix $\wt\M=(\wt\M_{ij})$ transposed to the universal $B$-Manin matrix $\M$. Due to the formulae~\eqref{dvMMM} this implies
\begin{align} \label{dBcop}
 &d_{\B^!}=d_\B^\cop, &&v_{\B^!}=v_\B.
\end{align}

\bbo{Infinite-dimensional case.} \label{bbInfcohom}
Let $\A,\wt\A\in\QA$ and $\B\in\FQA$. Then $\A=TV/(R)$, $\wt\A=T\wt V/(\wt R)$ and $\B=TW/(S)$ for some $V,\wt V\in\Vect$, $W\in\FVect$, $R\subset V^{\otimes2}$, $\wt R\subset\wt V^{\otimes2}$, $S\subset W^{\otimes2}$. Let us prove that there exists $\bfcohom(\B,\A)$ in the sense of p.~\ref{bbGenerIntHom}, where the role of parameter is played by $\B\in\FQA$.

\begin{Prop} \label{PropQA}
 For any $\B\in\FQA$ the functor $-\circ\B\colon\QA\to\QA$ has a left adjoint $\bfcohom(\B,-)\colon\QA\to\QA$. We have $\bfcohom(\B,\A)=TV'/(R')$ where $V'=\bfhom(W,V)$ and
\begin{align}
 R'=\{U\in\bfhom(W,V)^{\otimes2}\mid U(W^{\otimes2})\subset R,U(S)=0\}.
\end{align}
\end{Prop}

\noindent{\bf Proof.} Let $(w_i)_{i=1}^m$ and $(w^i)_{i=1}^m$ be dual bases of $W^*$ and $(W^*)^*=W$ respectively. These bases give identifications $W=(\KK^m)^*$ and $W^*=\KK^m$. Let $B\in\bfend(\KK^m\otimes\KK^m)$ be a matrix idempotent such that $S=\Im B^*$, then $\B=\gX_B(\KK)=TW/(\Im B^*)$. By identifying $W^*\otimes W$ with $\bfend(W)$ we obtain
\begin{align}
 \sum_{i,j}w_i\otimes w^i\otimes w_j\otimes w^j =\sum_{i,j,k,l}\big(B_{ij}^{kl}+(1-B)_{ij}^{kl}\big)w_k\otimes w^i\otimes w_l\otimes w^j,
\end{align}
where the coefficients $B^{ij}_{kl}\in\KK$ and $(1-B)^{ij}_{kl}=\delta^i_k\delta^j_l-B^{ij}_{kl}$ are entries of the operators $B$  and $1-B$ respectively.
Consider the linear map $\eta_V\colon V\to\bfhom(W,V)\otimes W$ defined in the proof of Prop.~\ref{PropCohomVect}, that is $\eta_V(v)=\sum_i(w_i\otimes v)\otimes w^i$. It generates a morphism $\eta\colon\A\to\bfcohom(\B,\A)\circ\B$ with the first order component $\eta_1=\eta_V$. Indeed, for any element $r=\sum_s x_s\otimes y_s\in R$ we have
\begin{multline} \label{sigmaetaetar}
 \sigma^{(23)}(\eta_V\otimes\eta_V)r=\sum_{s,i,j}(w_i\otimes x_s)\otimes (w_j\otimes y_s)\otimes w^i\otimes w^j=\\
 =\sum_{s,i,j,k,l}\big(B_{ij}^{kl}+(1-B)_{ij}^{kl}\big)w_k\otimes x_s\otimes w_l\otimes y_s\otimes w^i\otimes w^j.
\end{multline}
Since $B_{ij}^{kl}w^i\otimes w^j\in S$ and $(1-B)_{ij}^{kl}w_k\otimes x_s\otimes w_l\otimes y_s\in R'$, the element~\eqref{sigmaetaetar} belongs to $\bfhom(W,V)^{\otimes2}\otimes S+R'\otimes W^{\otimes2}$.

Now we need to prove that for any $\wt\A\in\QA$ and $f\in\Hom(\A,\wt\A\circ\B)$ there exists a unique $h\in\Hom\big(\bfcohom(\B,\A),\wt\A\big)$ making the diagram
\begin{align} \label{fhetaDiag}
\xymatrix{
\A\ar[rr]^{\eta\qquad}\ar[drr]_f&&\bfcohom(\B,\A)\circ\B\ar[d]^{h\circ\id_\B} \\
 && \wt\A\circ\B
}
\end{align}
commutative. For purposes of the next paragraph we do it for a slightly more general case $\wt\A\in\QAsc$. Let $\wt\gR=\wt\A_0$, then $\wt\A=\wt\gR\otimes T\wt V/(\wt R)$ for some vector space $\wt V\in\Vect$ and subspace $\wt R\subset\wt V\otimes\wt V$. We have $\wt\A_1=\wt\gR\otimes\wt V$ and $\wt\A=T_{\wt\gR}\A_1/(\wt\gR\otimes\wt R)$.
Note that any graded homomorphism $g\colon\A\to\wt\A$ is uniquely determined by its first order component $g_1\colon V\to\wt\gR\otimes\wt V$ and, conversely, this component gives the whole graded homomorphism $g\colon\A\to\wt\A$ if it is linear and satisfies $(g_1\dototimes g_1)R\subset\wt\gR\otimes\wt R$, where $g_1\dototimes g_1\colon V\otimes V\to\gR\otimes\wt V\otimes\wt V$ is defined as in p.~\ref{bbQSLS}, that is $g_1\dototimes g_1=(\mu_{\wt\gR}\otimes\id_{\wt V}\otimes\id_{\wt V})\cdot\sigma^{(23)}\cdot(g_1\otimes g_1)$. By using the natural isomorphism $(\gR\otimes\A)\circ\B=\gR\otimes(\A\circ\B)$ we obtain $\wt\A\circ\B=\wt\gR\otimes T(\wt V\otimes W)/(R_{\mathrm w})$, where $R_{\mathrm w}=\sigma^{(23)}(\wt V^{\otimes2}\otimes S+\wt R\otimes W^{\otimes2})$. In particular, $(\wt\A\circ\B)_1=\wt\gR\otimes\wt V\otimes W=\bfhom(W^*,\wt\gR\otimes\wt V)$.

Due to Prop.~\ref{PropCohomVect} there exists a unique linear map $h_1\colon\bfhom(W,V)\to\wt\gR\otimes\wt V$ such that the diagram
\begin{align}
\xymatrix{
V\ar[rr]^{\eta_1\qquad}\ar[drr]_{f_1}&&\bfhom(W,V)\otimes W\ar[d]^{h_1\otimes\id_W} \\
 && \wt\gR\otimes\wt V\otimes W
}
\end{align}
commute, it has the form $h_1(w_i\otimes v)=f_1(v)(w_i)$. This implies the uniqueness of $h$. Let us prove its existence. Let $U\in R'$, it has the form $U=\sum\limits_{s,i,j}U_s^{ij}w_i\otimes x_s\otimes w_j\otimes y_s$ for some $U_s^{ij}\in\KK$ and $x_s,y_s\in V$ such that $x_s\otimes y_s\in R$. Since $U(S)=0$ we have $\sum\limits_{s,i,j}U_s^{ij}B_{ij}^{kl}x_s\otimes y_s=0$, which implies $U=\sum\limits_{s,i,j,k,l}U_s^{ij}(1-B)_{ij}^{kl}w_k\otimes x_s\otimes w_l\otimes y_s$. Hence
\begin{align}
 (h_1\dototimes h_1)U=\sum_{s,i,j,k,l}U_s^{ij}(1-B)_{ij}^{kl}(f_1\dototimes f_1)(x_s\otimes y_s)(w_k\otimes w_l),
\end{align}
where $(f_1\dototimes f_1)(x_s\otimes y_s)$ belongs to $\wt\gR\otimes\wt V\otimes W\otimes\wt V\otimes W=\bfhom(W^*\otimes W^*,\wt\gR\otimes\wt V\otimes\wt V)$.
Since $(f_1\dototimes f_1)R\subset R_{\mathrm w}$, we obtain
\begin{align*}
(f_1\otimes f_1)(x_s\otimes y_s)\in(f_1\otimes f_1)R\subset\wt\gR\otimes(\id_V\otimes\sigma_{V,W}\otimes\id_W) (V\otimes V\otimes S+\wt R\otimes W\otimes W).
\end{align*}
The formula $(1-B^*)(\Im B^*)=0$ implies $\sum_{k,l}(1-B)_{ij}^{kl}(w_k\otimes w_l)(S)=0$, so  we derive $(h_1\otimes h_1)U\in\wt\gR\otimes\wt R$.

We proved that $\eta$ is a universal arrow from $\B$ to the functor $-\circ\B\colon\QA\to\QA$, hence this functor has a left adjoint $\bfcohom(\B,-)$, which acts on an object $\A$ as defined in the proposition (see~\cite[\S~4.1, Th.~2~(ii)]{Mcl}). \qed

\begin{Rem} \normalfont
 The proof of Prop.~\ref{PropQA} can be directly generalised to the case of arbitrary connected affinely generated graded algebras $\A,\wt\A,\B$ such that $\B$ is finely generated: $\dim\B_1<\infty$.
\end{Rem}

\bbo{Semi-connected case.} The generalised internal cohom found in the previous paragraph can be extended to the case $\bfC=(\QAsc,\circ)$, $\bfP=\FQA$. Remind that any $\D\in\QAsc$ has the form $\D=\gR\otimes\A$ for $\gR=\D_0\in\Alg$ and $\A\in\QA$.

\begin{Prop} \label{PropQAsc}
 For any $\B\in\FQA$ the functor $-\circ\B\colon\QAsc\to\QAsc$ has a left adjoint $\bfcohom(\B,-)\colon\QAsc\to\QAsc$. Let $\D=\gR\otimes\A$ for some $\gR\in\Alg$ and $\A\in\QA$. Then $\bfcohom(\B,\D)=\gR\otimes\bfcohom(\B,\A)$, where $\bfcohom(\B,\A)\in\QA$ is defined as in Prop.~\ref{PropQA}.
\end{Prop}

\noindent{\bf Proof.} Let $\wt\A\in\QAsc$. For an arbitrary graded homomorphism $\wh f\colon\gR\otimes\A\to\wt\A$ consider its zero component $\alpha=\wh f_0\in\Hom_{\Alg}(\gR,\wt A_0)$ and its restriction $f\in\Hom_{\GrAlg}(\A,\wt\A)$ to the subalgebra $\A\subset\gR\otimes\A$. They satisfy the commutativity condition $\alpha(r)f(a_1)=f(a_1)\alpha(r)$ $\;\forall\,r\in\gR,\,a_1\in\A_1$. Any pair of such morphisms subjected to the commutativity condition uniquely defines a morphism $\wh f=(\alpha,f)\in\Hom_{\GrAlg}(\gR\otimes\A,\wt\A)$ (see~\cite[\S~2.5]{S}, cf.~p.~\ref{bbQSLS}).

Let $\wh f=(\alpha,f)\colon\D\to\wt\A$ be a graded homomorphism to an algebra $\wt\A\in\QAsc$, where $\alpha\colon\gR\to\wt\A_0$ and $f\colon\A\to\wt\A$ are morphisms subjected to the commutativity condition.
As it was proved in p.~\ref{bbInfcohom} there exists a unique graded homomorphism $h\colon\bfcohom(\B,\A)\to\wt\A$ making the digram~\eqref{fhetaDiag} commutative. By substituting $f(v)=\sum_ih_1(w_i\otimes v)\otimes w^i$, $v\in V$, to the commutativity condition $\alpha(r)f(v)=f(v)\alpha(r)$ we obtain the commutativity condition for the pair $(\alpha,h)$, so it gives a graded homomorphism $\wh h=(\alpha,h)\colon\gR\otimes\bfcohom(\B,\A)\to\wt\A$ and we obtain the commutative diagram
\begin{align} \label{alphafhetaDiag}
\xymatrix{
 \gR\otimes\A\ar[rr]^{\id_\gR\otimes\eta\quad\qquad}\ar[drr]_{\wh f=(\alpha,f)}
 &&\gR\otimes\bfcohom(\B,\A)\circ\B\ar[d]^{\wh h\circ\id_B=(\alpha,h\circ\id_\B)} \\
 && \wt\A\circ\B
}
\end{align}
The uniqueness of $h$ implies the uniqueness of $\wh h$. \qed

\bb{Internal cohom on the morphisms.}
Prop.~\ref{PropQAsc} implies the existence of the generalised internal cohom-functor $\bfcohom\colon\FQA^\op\times\QAsc\to\QAsc$ and gives its values on objects. Let us calculate this functor on morphisms. Let $\B,\B'\in\FQA$, $\gR,\gR'\in\Alg$ and $\A,\A'\in\QA$. Set $W=\B_1$, $W'=\B'_1$, $V=\A_1$, $V'=\A'_1$. Consider arbitrary morphisms $f\colon\B'\to\B$ and $g\colon\gR\otimes\A\to\gR'\otimes\A'$. They are uniquely determined by the components $f_1\colon W'\to W$, $\alpha=g_0\colon\gR\to\gR'$ and $g_1=(\alpha,t)\colon\gR\otimes V\to\gR'\otimes V'$, where $t\colon V\to\gR'\otimes V'$. Then the morphism $\bfcohom(f,g)\colon\gR\otimes\bfcohom(\B,\A)\to\gR'\otimes\bfcohom(\B',\A')$ is uniquely determined by the components
\begin{align}
 &\bfcohom(f,g)_0\colon\gR\to\gR',  \label{cohomfg0} \\
 &\bfcohom(f,g)_1\colon\gR\otimes\bfhom(W,V)\to\gR'\otimes\bfhom(W',V').  \label{cohomfg1}
\end{align}

\begin{Prop}
 The components~\eqref{cohomfg0},~\eqref{cohomfg1} equal
\begin{align}
 &\bfcohom(f,g)_0=\alpha, && \bfcohom(f,g)_1=(\alpha,t_*f_1^*), \label{cohomfg01}
\end{align}
where $g_0=\alpha$, $g_1=(\alpha,t)$ and the operator $t_*f_1^*\colon\bfhom(W,V)\to\gR'\otimes\bfhom(W',V')$ is the composition
\begin{align}
 \bfhom(W,V)\xrightarrow{f_1^*}\bfhom(W',V)\xrightarrow{t_*}
\bfhom(W',\gR'\otimes V')=\gR'\otimes\bfhom(W',V').
\end{align}
\end{Prop}

\noindent{\bf Proof.} Due to the formulae~\eqref{varethetacoev} and~\eqref{thetacohom} we obtain
\begin{align}
 \big(\bfcohom(f,g)\circ\id\big)\cdot\coev=\vartheta\big(\bfcohom(f,g)\big)= 
\big(\id\circ f\big)\cdot\coev\cdot g,
\end{align}
where $\coev$ coincides with the horizontal arrow in~\eqref{alphafhetaDiag}. This gives us the commuting diagram
\begin{align}
\xymatrix{
\gR\otimes\A\ar[r]^{g}\ar[d]^{\id\otimes\eta} & \gR'\otimes\A'\ar[rr]^{\id\otimes\eta\qquad\qquad} &&
\gR'\otimes\bfcohom(\B',\A')\circ\B'\ar[d]^{\id\otimes\id\circ f} \\
\gR\otimes\bfcohom(\B,\A)\circ\B\ar[rrr]^{\bfcohom(f,g)\circ\id} &&&
\gR'\otimes\bfcohom(\B',\A')\circ\B
}
\end{align}
By taking zero and first order components in this diagram we obtain~\eqref{cohomfg01}. \qed

Consider the case $\A,\A'\in\FQA$. Up to isomorphisms we have $\A=\gX_A(\KK)$, $\A'=\gX_{A'}(\KK)$, $\B=\gX_B(\KK)$, $\B=\gX_{B'}(\KK)$ for some matrix idempotents $A,A',B,B'$. Then the morphisms $f\colon\B'\to\B$ and $g\colon\gR\otimes\A\to\gR'\otimes\A'$ have the form
\begin{align}
 &f=f_K\colon\gX_{B'}(\KK)\to\gX_{B}(\KK) &&\text{and} &
 &g=(\alpha,f_M)\colon\gX_A(\gR)\to\gX_{A'}(\gR')
\end{align}
respectively for a $(B',B)$-Manin matrix $K$ over $\KK$, an $(A,A')$-Manin matrix $M$ over $\gR'$ and a homomorphism $\alpha\colon\gR\to\gR'$ such that $\alpha(r)M^i_j=M^i_j\alpha(r)$ $\;\forall\,r\in\gR$. Let $\M$ and $\N$ be the universal $(A,B)$- and $(A',B')$-Manin matrices respectively. Their entries are generators of $\bfcohom\big(\gX_B(\KK),\gX_A(\KK)\big)$ and $\bfcohom\big(\gX_{B'}(\KK),\gX_{A'}(\KK)\big)$. The graded homomorphism
 $$\bfcohom(f,g)\colon\gR\otimes\bfcohom\big(\gX_B(\KK),\gX_A(\KK)\big)\to \gR'\otimes\bfcohom\big(\gX_{B'}(\KK),\gX_{A'}(\KK)\big)$$
 acts on the generators as
\begin{align}
 &\bfcohom(f,g)\colon r\otimes1\mapsto\alpha(r)\otimes1, \qquad r\in\gR, \\
 &\bfcohom(f,g)\colon 1\otimes\M^i_j\mapsto\sum_{a,b} M^i_a\otimes\N^a_b K^b_j. \label{cohomdfM}
\end{align}
Note that the matrix $M\N K$ appeared in~\eqref{cohomdfM} is a $(A,B)$-Manin matrix, which follows directly from~\cite[Prop.~2.26]{S}.

\section{Quantum Representation Theory}
\label{sec5}

By taking the category $\bfC=\QAsc^\op$ with the Manin product `$\circ$' and parameter category $\bfP=\FQA^\op$ we obtain a class of representations on quantum linear spaces. {\it Quantum Representation Theory} presented here investigates representations for this case and, more generally, for $\bfC=(\GrAlg^\op,\circ)$. It can be regarded as a generalisation of Representation Theory on the usual vector spaces (classical representation theory). The latter one can be embedded into Quantum Representation Theory in two ways: by the functor $S^*$ or $T^*$. The binary and duality operations with representations can be generalised to the quantum case.

\subsection{Quantum representations}

First we describe objects which we are going to `represent'. In the classical representation theory these were algebraic monoids and algebras. In the quantum case we need to consider monoids in $(\QAsc^\op,\circ)$ or, more generally, in $(\GrAlg^\op,\circ)$. We call these monoids {\it quantum algebras}.

\bb{Quantum monoids.} \label{bbQM}
A quantum analogue of algebraic monoids is a comonoid in the monoidal category $(\Alg,\otimes)$. This is an algebra $\gR\in\Alg$ with comultiplication $\Delta_\gR\colon\gR\to\gR\otimes\gR$ and counit $\varepsilon_\gR\colon\gR\to\KK$.
These are algebra homomorphisms satisfying the equations $(\varepsilon_\gR\otimes\id_\gR)\cdot\Delta_\gR=\id_\gR=(\id_\gR\otimes\varepsilon_\gR)\cdot\Delta_\gR$ and the coassociativity condition $(\Delta_\gR\otimes\id_\gR)\cdot\Delta_\gR=(\id_\gR\otimes\Delta_\gR)\cdot\Delta_\gR$. Such comonoids $(\gR,\Delta_\gR,\varepsilon_\gR)$ are exactly the bialgebras: $\Comon(\Alg,\otimes)=\Bimon(\Vect,\otimes)$.

To define representations of a quantum monoid we need to consider it as a quantum algebra, so we should translate it by means of the (not full) embedding $\Alg\hookrightarrow\QAsc$, $\gR\mapsto\gR\otimes\KK[u]$.
This is a strong monoidal functor $(\Alg,\otimes)\hookrightarrow(\QAsc,\circ)$, it translates the bialgebra $(\gR,\Delta_\gR,\varepsilon_\gR)$ to the comonoid
\begin{align} \label{CCgR}
 &\OO_\gR=(\gR\otimes\KK[u],\Delta,\varepsilon)\in\Comon(\QAsc,\circ), &&
 &\Delta=\Delta_\gR\otimes\id_{\KK[u]}, &
 &\varepsilon=\varepsilon_\gR\otimes\id_{\KK[u]}. 
\end{align}
Component-wise we have $(\gR\otimes\KK[u])_k=\gR$, $\Delta_k=\Delta_\gR$ and $\varepsilon_k=\varepsilon_\gR$. The comonoid $\OO_\gR$ is a quantum monoid considered as a comonoid in $(\QAsc,\circ)$. Note, however, that not all the comonoid structure on the semi-connected graded algebra $\gR\otimes\KK[u]$ has this form (because the embedding functor is not full).

\bb{Comonoids in $(\GrAlg,\circ)$.} \label{bbComonGrAlg}
A linear version of an algebraic monoid is an algebra. Its quantum analogue is a comonoid in $(\QA,\circ)$. A quantum semi-linear monoid is a comonoid in $(\QAsc,\circ)$. For wider generality we describe comonoids in $(\GrAlg,\circ)$, which we call {\it quantum algebras}. 
Such comonoid is a graded algebra $\A=\bigoplus\limits_{k\in\NN_0}\A_k$ with graded homomorphisms $\Delta\colon\A\to\A\circ\A$ and $\varepsilon\colon\A\to\KK[u]$, their components are linear maps $\Delta_k\colon\A_k\to\A_k\otimes\A_k$ and $\varepsilon_k\colon\A_k\to\KK$ satisfying the conditions $(\A_k,\Delta_k,\varepsilon_k)\in\Comon(\Vect,\otimes)$ $\;\forall\,k\in\NN_0$ and $\Delta_{k+l}(a_lb_l)=\Delta_k(a_k)\Delta_l(b_l)$ $\;\forall\,k,l\in\NN_0$, $a_k\in\A_k$, $b_l\in\A_l$.

By composing $\Delta$ and $\varepsilon$ with the embedding $\A\circ\A\hookrightarrow\A\otimes\A$ and evaluation $\KK[u]\to\KK$ at $u=1$ we get a comultiplication $\Delta_\A\colon\A\to\A\otimes\A$ and a counit $\varepsilon_\A\colon\A\to\KK$ respectively (they are not graded). The conditions on the graded homomorphisms $\Delta$ and $\varepsilon$ imply that $(\A,\Delta_\A,\varepsilon_\A)$ is a bialgebra. Conversely, let a graded algebra $\A\in\GrAlg$ has a bialgebra structure: $(\A,\Delta_\A,\varepsilon_\A)\in\Comon(\Alg,\otimes)$, then it gives a comonoid $(\A,\Delta,\varepsilon)$ in $(\GrAlg,\circ)$ iff
\begin{align} \label{DeltaAk}
 &\Delta_\A(\A_k)\subset\A_k\otimes\A_k &&k\in\NN_0.
\end{align}

\bbo{Comonoids in $(\QAsc,\circ)$ and $(\QA,\circ)$.}
If $\A\in\QAsc$, then it is enough to require~\eqref{DeltaAk} for $k=0$ and $k=1$:
\begin{align} \label{Delta01}
 &\Delta_\A(\A_0)\subset\A_0\otimes\A_0, &
 &\Delta_\A(\A_1)\subset\A_1\otimes\A_1.
\end{align}
In this case $\Delta$ and $\varepsilon$ are completely determined by their zero and first order components $\Delta_0=\alpha$, $\Delta_1=(\alpha,t)$, $\varepsilon_0=\beta$, $\varepsilon_1=\beta\otimes\varepsilon_V$ (see~p.~\ref{bbQSLS} for details).

In the connected case $\A\in\QA$ we have $\A_0=\KK$, so it is enough to impose the condition~\eqref{DeltaAk} for $k=1$, the graded homomorphisms $\Delta$ and $\varepsilon$ are uniquely determined by their first order components $\Delta_1=\Delta_V\colon V\to V\otimes V$ and $\varepsilon_1=\varepsilon_V\colon V\to\KK$, where $V=\A_1\in\Vect$. Then the commutativity of the diagrams~\eqref{DiagAlphat} reduces to the requirement that these components gives a structure of coalgebra $(V,\Delta_V,\varepsilon_V)\in\Comon(\Vect,\otimes)$. The comonoid structure on a connected quadratic algebra $\A=TV/(R)$ is given by a coalgebra $(V,\Delta_V,\varepsilon_V)\in\Comon(\Vect,\otimes)$ such that $(\Delta_V\otimes\Delta_V)R\subset\sigma^{(23)}(R\otimes V^{\otimes2}+V^{\otimes2}\otimes R)$ and $(\varepsilon_V\otimes\varepsilon_V)R=0$.

\bb{Quantum representations.}
 Let $\OO=(\A,\Delta,\varepsilon)$ be a comonoid in $(\QAsc,\circ)$, that is a semi-connected quadratic algebra $\A$ with a bialgebra structure $(\Delta_\A,\varepsilon_\A)$ satisfying~\eqref{Delta01}. Due to the propositions~\ref{PropQA}, \ref{PropQAsc} we have the category $\Corep_{\FQA}(\OO)$, which consists of corepresentations of $\OO$ on finitely generated quadratic algebras $\B\in\FQA$. These corepresentations can be regarded as representations of the corresponding monoid in $(\QAsc^\op,\circ)$ on finite-dimensional quantum linear spaces. Let us define a slightly more general notion.

\begin{Def} \normalfont
 A {\it quantum representation} is a corepresentation $\omega\colon\bfcoend(\B)\to \OO$ of a comonoid $\OO=(\A,\Delta,\varepsilon)\in\Comon(\GrAlg,\circ)$ on a quadratic algebra $\B\in\FQA$. A morphism between corepresentations $\omega\colon\bfcoend(\B)\to\OO$ and $\omega'\colon\bfcoend(\B')\to\OO$ is a morphism $f\colon\B\to\B'$ in $\FQA$ such that the diagram
\begin{align} \label{omegafC}
\xymatrix{
 \bfcohom(\B',\B)\ar[rrr]^{\bfcohom(f,\id_\B)}\ar[d]^{\bfcohom(\id_{\B'},f)} &&& \bfcohom(\B,\B)\ar[d]^{\omega} \\
\bfcohom(\B',\B')\ar[rrr]^{\qquad\omega'} &&& \A
}
\end{align}
commute. The objects $(\B,\omega)$, where $\B\in\FQA$ and $\omega\colon\bfcoend(\B)\to\OO$ is a corepresentation, form a category with such defined morphisms. We denote it by $\Corep_{\FQA}(\OO)$.
\end{Def}

Note that in the case $\OO\in\Comon(\QAsc,\circ)$ the category $\Corep_\FQA(\OO)$ coincides with one defined in p.~\ref{bbCorep}, so in this case Theorem~\ref{ThCorepLcoact} implies that it is a full subcategory of $\Lcoact(\OO)$.

More correctly, we should say that the category of quantum representations is rather $\Rep_{\FQA^\op}(\MM)$, where $\MM\in\Mon(\GrAlg^\op,\circ)$ is the monoid corresponding to the comonoid $\OO$. This is the category opposite to $\Corep_\FQA(\OO)$. The morphisms of quantum representations are reversed morphisms of corepresentations. However, it is more convenient for us to work in terms of category $\Corep_\FQA(\OO)$.

\bb{Multiplicative Manin matrices.} \label{bbMultMM}
An $m\times m$ matrix $M=(M^i_j)$ with entries in a bialgebra $\gR=(\gR,\Delta_\gR,\varepsilon_\gR)$ (or, more generally in a coalgebra) is called {\it multiplicative} if
\begin{align} \label{MultM}
 &\Delta_\gR(M^i_j)=\sum_{k=1}^mM^i_k\otimes M^k_j, &
 &\varepsilon_\gR(M^i_j)=\delta^i_j, &&i,j=1,\ldots,m
\end{align}
(see~\cite[\S~2.6]{Manin88}, \cite[\S~4.1.1]{ManinBook91}).
For example, the formulae~\eqref{dvMMM} means exactly that the universal $B$-Manin matrix $\M=(\M^i_j)$ is a multiplicative $m\times m$ matrix over the algebra $\bfcoend\big(\gX_B(\KK)\big)$. Recall that $\M^i_j$ are generators of this algebra.

Let us establish a relationship between the notions of multiplicative Manin matrix and of quantum representation. Consider a comonoid $\OO=(\A,\Delta,\varepsilon)$ in $(\GrAlg,\circ)$.

\begin{Th} \label{ThCorepMMM}
 Let $B\in\End(\KK^m\otimes\KK^m)$ be a matrix idempotent. Any corepresentation $\omega\colon\bfcoend(\B)\to\OO$ of the comonoid $\OO$ on the quadratic algebra $\B=\gX_B(\KK)$ acts on the generators by the formula
\begin{align} \label{omegaM}
 \omega(\M^i_j)=M^i_j,
\end{align}
where $M=(M^i_j)$ is a multiplicative first order $B$-Manin matrix over $\A$, that is an $m\times m$ matrix with entries $M^i_j\in\A_1$ satisfying~\eqref{MultM} and
\begin{align} \label{BMMB}
 B M^{(1)}M^{(2)}(1-B)=0.
\end{align}
Conversely, any such matrix $M$ defines a corepresentation of the comonoid $\OO$ on $\B$ by the formula~\eqref{omegaM}.
\end{Th}

\noindent{\bf Proof.} Any graded homomorphism $\omega\colon\bfcoend(\B)\to\A$ is uniquely determined by its values $M^i_j\in\A_1$ on the generators $\M^i_j$. The values~\eqref{omegaM} give a graded homomorphism $\omega$ iff they belong to the component $\A_1$ and satisfy the same commutation relations that $\M^i_j$ do, i.e.the relations~\eqref{BMMB}. The commutativity of the corresponding diagrams~\eqref{diagMonoidMor} is equivalent to the condition~\eqref{MultM}. \qed

\noindent{\bf Warning.} The same multiplicative first order matrix $M$ defines different corepresentations (different objects of the category $\Corep_\FQA(\OO)$) on different quadratic algebras $\B=\gX_B(\KK)$ and $\B'=\gX_{B'}(\KK)$ if $\B_1=\B'_1$ and $M$ is simultaneously $B$- and $B'$-Manin matrix. In other words, a matrix $M$ is not enough to fix a corepresentation, one also needs to know the quadratic algebra $\B$ (and also the comonoid $\OO$, of cause).

By a quantum representation of a quantum monoid we understand a corepresentation of the comonoid, obtained by embedding of the bialgebra $(\gR,\Delta_\gR,\varepsilon_\gR)$ via $\Alg\hookrightarrow\QAsc$. This is a corepresentation of the comonoid $\OO_\gR=(\gR\otimes\KK[u],\Delta_\gR\otimes\id_{\KK[u]},\varepsilon_\gR\otimes\id_{\KK[u]})$ defined by~\eqref{CCgR}. Theorem~\ref{ThCorepMMM} implies in this case a one-to-one correspondence between corepresentations of $\OO_\gR$ on $\B=\gX_B(\KK)$ and multiplicative $B$-Manin matrices over $\gR$.

A trivial but important example of quantum representation is the identical morphism $\bfcoend(\B)\xrightarrow{\id}\bfcoend(\B)$. It is a corepresentation of the comonoid $\OO=\bfcoend(\B)$ on the object $\B\in\FQA$. It is given by the universal $B$-Manin matrix $\M$. The corresponding coaction coincides with the coevaluation~\eqref{coevAB} for $\A=\B$.

\begin{Rem} \normalfont
 In~\cite[Def. 20]{FRT89} the following generalisation of representation of a group was introduced. Representation of a quantum group $\gR\in\Comon(\Alg,\otimes)$ on a vector space $W\in\FVect$ was defined as a coalgebra homomorphism $\rho\colon\bfend(W)^*\to\gR$. In the basis $(w_i)_{i=1}^m$ it is given by an arbitrary multiplicative first order $m\times m$ matrix $M=(M^i_j)$ over $\gR$.
 On one hand, $\rho$ can be considered as a particular case of quantum representation: the corresponding corepresentation $\omega\colon\bfcoend(TW^*)=T\bfend(W)^*\to\OO_\gR$ is given by the matrix $M$ or, equivalently, by the condition $\omega_1=\rho$. In that sense the representations $\rho$ of the quantum group $\gR$ form the subcategory $\Corep_\FTA(\OO_\gR)$ in $\Corep_\FQA(\OO_\gR)$, where $\FTA\subset\FQA$ is the subcategory of finitely generated tensor algebras (it is coreflective in the sense~\cite[\S~4.3]{Mcl} with coreflection $\FQA\to\FTA$, $\B\mapsto T\B_1$).
On the other hand, any quantum representation of a comonoid $\OO\in\Comon(\GrAlg,\circ)$ on $\B\in\FQA$ is described by a multiplicative first order matrix, which, in particular, gives corepresentation of the comonoid $\OO$ on $T\B_1\in\FTA$. We get a faithful functor $\Corep_\FQA(\OO)\to\Corep_\FTA(\OO)$ (it~is also a coreflection).
\end{Rem}

\bbo{Morphisms of quantum representations.}
Consider two quantum representations $\omega\colon\bfcoend(\B)\to\OO$ and $\omega'\colon\bfcoend(\B')\to\OO$ of the same $\OO=(\A,\Delta,\varepsilon)\in(\GrAlg,\circ)$ on quadratic algebras $\B,\B'\in\FQA$. Suppose that these algebras are defined by matrix idempotents $B\in\End(\KK^m\otimes\KK^m)$, $B'\in\End(\KK^{m'}\otimes\KK^{m'})$, i.e. $\B=\gX_B(\KK)$, $\B'=\gX_{B'}(\KK)$. Due to Theorem~\ref{ThCorepMMM} these corepresentations are given by multiplicative first order $B$- and $B'$-Manin matrices $M$ and $M'$ respectively: $\omega(\M)=M$, $\omega(\M')=M'$, where $\M$ and $\M'$ are the universal $B$- and $B'$-Manin matrices.

\begin{Prop} \label{PropMorCorep}
 Any morphism $f\colon(\B,\omega)\to(\B',\omega')$ in the category $\Corep_\FQA(\OO)$ has the form $f=f_K\colon\gX_B(\KK)\to\gX_{B'}(\KK)$ for an $m\times m'$ matrix $K\in\Hom(\KK^{m'},\KK^m)$ satisfying $BK^{(1)}K^{(2)}(1-B')=0$ and
\begin{align} \label{KMMK}
 KM'=MK.
\end{align}
This gives a bijection between the morphisms $f\colon(\B,\omega)\to(\B',\omega')$ and $(B,B')$-Manin matrices $K$ over $\KK$ satisfying~\eqref{KMMK}.
\end{Prop}

\noindent{\bf Proof.} The graded homomorphisms $f\colon\gX_B(\KK)\to\gX_{B'}(\KK)$ is in one-to-one correspondence with the $(B,B')$-Manin matrices over $\KK$ (see p.~\ref{bbManinMatr}). By taking the entries of the universal $(B,B')$-Manin matrix in the left upper corner of the diagram~\eqref{omegafC} we obtain exactly the condition~\eqref{KMMK}. Since these entries generate the algebra $\bfcohom(\B',\B)$ the commutativity of the diagram~\eqref{omegafC} is equivalent to~\eqref{KMMK}. \qed

\begin{Cor} \label{CorMorCorep}
Let $\B_1=\B'_1$. Suppose that this identification induces a morphism $f\colon\B\to\B'$ in $\FQA$, that is $f_1=\id_{\B_1}$ (this is not always true). Then $f$ is a morphism $(\B,\omega)\to(\B',\omega')$ in $\Corep_\FQA(\OO)$ iff the corepresentations $\omega$ and $\omega'$ are given by the same matrix $M=M'$.
\end{Cor}

\noindent{\bf Proof.} Since $f_1=\id_{\B_1}$ we have $f=f_K$ for the identity matrix $K=1$, so the condition~\eqref{KMMK} reduces to $M=M'$ (the equation $B(1-B')=0$ is equivalent to the requirement that the matrix $K=1$ defines a morphism $f_K\colon\gX_B(\KK)\to\gX_{B'}(\KK)$ in the category $\FQA$). \qed

\bb{Coactions in $(\GrAlg,\circ)$.}
Let $\OO$ be a comonoid in $(\GrAlg,\circ)$ corresponding to a bialgebra $(\A,\Delta_\A,\varepsilon_\A)$ in the sense of p.~\ref{bbComonGrAlg}. Consider a coaction $\delta\colon\B\to\A\circ\B$ of $\OO$ on $\B\in\GrAlg$. By composing it with the embedding $\A\circ\B\hookrightarrow\A\otimes\B$ we obtain a coaction $\delta_\A\colon\B\to\A\otimes\B$ of the bialgebra $(\A,\Delta_\A,\varepsilon_\A)$ in the monoidal category $(\Alg,\otimes)$. Conversely, a coaction $\delta_\A\colon\B\to\A\otimes\B$ of the bialgebra $(\A,\Delta_\A,\varepsilon_\A)$ gives a coaction $\delta\colon\B\to\A\circ\B$ of $\OO$ in $(\GrAlg,\circ)$ iff
\begin{align} \label{deltaAk}
 &\delta_\A(\B_k)\subset\A_k\otimes\B_k, &&k\in\NN_0.
\end{align}
If $\B\in\QA$ or $\B\in\QAsc$ then it is enough to require this condition for 
$k=1$ only or for $k=0,1$ only respectively.

\begin{Prop} \label{PropCoa}
 Each coaction $\delta\colon\gX_B(\KK)\to\A\circ\gX_B(\KK)$ has the form
\begin{align} \label{deltaM}
 \delta(w^i)=\sum_j M^i_j\otimes w^j,
\end{align}
where $M=(M^i_j)$ is a multiplicative first order $B$-Manin matrix and $(w^j)$ is the standard basis of $(\KK^m)^*$. Any such matrix $M$ define a coaction of $\OO$ on $\gX_B(\KK)$. If $\omega$ is a corepresentation of $\OO$ on $\gX_B(\KK)$ given by the matrix $M$, then the diagram
\begin{align} \label{deltaomegaetaDiag}
\xymatrix{
\gX_B(\KK)\ar[rr]^{\eta\qquad}\ar[drr]_{\delta}&&\bfcoend\big(\gX_B(\KK)\big)\circ\gX_B(\KK)\ar[d]^{\omega\circ\id} \\
 && \A\circ\gX_B(\KK)
}
\end{align}
commute. In particular, if $\A\in\QAsc$, then $\vartheta(\omega)=\delta$.
\end{Prop}

\noindent{\bf Proof.} Any graded homomorphism $\delta\colon\gX_B(\KK)\to\A\circ\gX_B(\KK)$ has the form~\eqref{deltaM} for a $B$-Manin matrix $M$ (see p.~\ref{bbManinMatr}). By taking $w^i$ in the left upper corners of the diagrams~\eqref{deltaDiag} we obtain the conditions~\eqref{MultM}. The commutativity of the diagram~\eqref{deltaomegaetaDiag} is checked on the generators $w^i$ in a straightforward way. \qed

We established the bijection between corepresentations $\omega\colon\bfcohom(\B)\to\OO$ and coactions $\delta\colon\B\to\A\circ\B$ for an arbitrary comonoid $\OO=(\A,\Delta,\varepsilon)\in\Comon(\GrAlg,\circ)$. This generalises the case $\bfC=\QAsc$ of the bijection~\eqref{thetaP}, so we denote it by the same letter: $\vartheta\colon\omega\leftrightarrow\delta$.

\begin{Prop} \label{PropcoactionMor}
 Let $\delta$ and $\delta'$ be coactions of $\OO$ on $\gX_B(\KK)$ and $\gX_{B'}(\KK)$. Let $M$ and $M'$ be the corresponding multiplicative first order $B$- and $B'$-Manin matrices. Then morphisms $f\colon\big(\gX_B(\KK),\delta\big)\to\big(\gX_{B'}(\KK),\delta'\big)$ in $\Lcoact(\OO)$ are homomorphisms $f=f_K$ for $(B,B')$-Manin matrices $K$ over $\KK$ satisfying~\eqref{KMMK}. Hence $\Corep_\FQA(\OO)$ is a subcategory of $\Lcoact(\OO)$.
\end{Prop}

\noindent{\bf Proof.} By substituting $w^i$ to the left upper corner of the left diagram~\eqref{LactMor} we obtain the condition~\eqref{KMMK}. \qed

\begin{Rem} \normalfont
 Non-linear versions of the quantum representations are the comodule algebras~\cite[Def.~3.7.1]{Kass}.
 By means of the notion of coaction in $(\Alg,\otimes)$ one can define a {\it comodule algebra} over a bialgebra $\gR=(\gR,\Delta_\gR,\varepsilon_\gR)\in\Comon(\Alg,\otimes)$ as an object of the category $\Lcoact(\gR)$. The strict monoidal embedding $(\Alg,\otimes)\hookrightarrow(\GrAlg,\circ)$, $\gR\mapsto\gR\otimes\KK[u]$, induces (not fully) faithful functor $\Lcoact(\gR)\to\Lcoact(\OO_\gR)$, which embeds the comodule algebras over the bialgebra $\gR$ to the category of coactions of comonoid $\OO_\gR$ defined by the formula~\eqref{CCgR}.
On the other hand, the forgetful functor $U\colon\GrAlg\to\Alg$ has a strict monoidal structure $(\GrAlg,\otimes)\to(\Alg,\otimes)$ and hence induces functors between comonoids and between coactions. Thus we get the functor $\Lcoact(\OO)\to\Lcoact(\A,\Delta_\A,\varepsilon_\A)\to\Lcoact(\BB)$, which translates coactions $\delta$ of a comonoid $\OO=(\A,\Delta,\varepsilon)\in\Comon(\GrAlg,\circ)$ on an object $\B\in\GrAlg$ to the comodule algebra $U\B\in\Alg$ over the bialgebra $\BB=(U\A,U\Delta_\A,U\varepsilon_\A)$, where the coaction is given by the homomorphism $U\delta_\A$. In particular, we obtain a forgetful functor from the category of quantum representations of the comonoid $\OO$ to the category of comodule algebras over $\BB$, which generalises the classical situation where the linear representations are regarded as actions in $(\AlgSet,\times)$.
\end{Rem}

\subsection{Operations with quantum representations}

We have tree operations with vector spaces: direct sum, tensor product and duality. They induce operations with representations (in the case of tensor product and duality we need an additional structure on the represented object). Similarly, the operations with quadratic algebras give operations with quantum representations. We construct these operations in terms of coactions and Manin matrices.

\bb{Direct sum of Manin matrices.}
According to the Manins' ideas~\cite{Manin88} the quantum analogue of the direct sum of vector spaces (in the `even' interpretation) is the (even) tensor product of quadratic algebras (see p.~\ref{bbSQVS} and p.~\ref{bbLaxMonS}). Therefore, we describe first the tensor product $\B\otimes\C$ for quadratic algebras $\B,\C\in\FQA$.

For arbitrary matrix idempotents $B\in\End(\KK^m\otimes\KK^m)$ and $C\in\End(\KK^n\otimes\KK^n)$ define the idempotent $D=\DiS(B,C)\in\End(\KK^{m+n}\otimes\KK^{m+n})$ with the only non-zero entries
\begin{align}
 &D^{ij}_{kl}=B^{ij}_{kl}, && D^{\bar a,\bar b}_{\bar c,\bar d}=C^{ab}_{cd} &
 &D^{i,\bar a}_{i,\bar a}=D^{\bar a,i}_{\bar a,i}=\frac12, &
 &D^{i,\bar a}_{\bar a,i}=D^{i,\bar a}_{\bar a,i}=-\frac12,
\end{align}
where $i,j,k,l=1,\ldots,m$, \quad $a,b,c,d=1,\ldots,n$,\quad $\bar a=m+a$.
Then $\gX_B(\KK)\otimes\gX_C(\KK)=\gX_D(\KK)$.

Let $B'\in\End(\KK^{m'}\otimes\KK^{m'})$ and $C'\in\End(\KK^{n'}\otimes\KK^{n'})$ be idempotents as well. Denote $D'=\DiS(B',C')\in\End(\KK^{m'+n'}\otimes\KK^{m'+n'})$. Graded operators $f\colon\gX_B(\KK)\to\gX_{B'}(\gR)$ and $g\colon\gX_{C}(\KK)\to\gX_{C'}(\gR)$ give one more graded operator by the composition
\begin{multline}
 \gX_D(\KK)=\gX_B(\KK)\otimes\gX_C(\KK)\xrightarrow{f\otimes g}
 \gR\otimes\gX_{B'}(\KK)\otimes\gR\otimes\gX_{C'}(\KK)
\xrightarrow{\sigma^{(23)}} \\
\gR\otimes\gR\otimes\gX_{B'}(\KK)\otimes\gX_{C'}(\KK)
\xrightarrow{\mu_\gR\otimes\id\otimes\id}
\gR\otimes\gX_{B'}(\KK)\otimes\gX_{C'}(\KK)=\gX_{D'}(\gR). \label{grOp}
\end{multline}
The condition that this operator is an algebra homomorphism can be written in terms of Manin matrices.

\begin{Prop} \label{PropLMN}
 Let $M$ and $N$ be $m\times m'$ and $n\times n'$ matrices over $\gR$. Consider their direct sum
\begin{align} \label{MoplusN}
 L=M\oplus N=
\begin{pmatrix}
 M & 0 \\
 0 & N
\end{pmatrix}.
\end{align}
This is an $(m+n)\times(m'+n')$ matrix over $\gR$ with the only non-zero entries $L^i_j=M^i_j$, $L^{\bar a}_{\bar b}=N^a_b$, where $i,j=1,\ldots,m$, $a,b=1,\ldots,n$. The matrix $L$ is a $(D,D')$-Manin matrix iff $M$ is $(B,B')$-Manin matrix, $N$ is $(C,C')$-Manin matrix and
\begin{align} \label{MNComm}
 &M^i_jN^a_b=N^a_bM^i_j, &&\forall\,i,j=1,\ldots,m, \quad a,b=1,\ldots,n. 
\end{align}
In other words, the graded operator~\eqref{grOp} is a graded homomorphisms iff $f$, $g$ are algebra homomorphisms and the corresponding Manin matrices entry-wise commute.
\end{Prop}

\noindent{\bf Proof.} By writing the condition $DL^{(1)}L^{(2)}(1-D')=0$ in entries one sees that it is equivalent to the conditions $BM^{(1)}M^{(2)}(1-B')=0$, $CN^{(1)}N^{(2)}(1-C')=0$ and~\eqref{MNComm}. \qed

\bb{Direct sum of quantum representations.} \label{bbDS}
Consider coactions $\delta\colon\B\to\A\circ\B$ and $\gamma\colon\C\to\A\circ\C$ of the comonoid $\OO=(\A,\Delta,\varepsilon)\in\Comon(\GrAlg,\circ)$ on the objects $\B,\C\in\GrAlg$. Let $(\A,\Delta_\A,\varepsilon_\A)\in\Comon(\Alg,\otimes)$ be the bialgebra corresponding to the comonoid $\OO$ and $\delta_\A\colon\B\to\A\otimes\B$, $\gamma_\A\colon\C\to\A\otimes\C$ be the corresponding coactions of the bialgebra $(\A,\Delta_\A,\varepsilon_\A)$. The letter one is the bimonoid $(\A,\mu_\A,\varepsilon_\A,\Delta_\A,\varepsilon_\A)$ in $(\Vect,\otimes)$, where $\mu_\A\colon\A\otimes\A\to\A$ and $\eta_\A\colon\KK\to\A$ are the multiplication and unity map in the algebra $\A$, i.e. $\mu_\A(a\otimes b)=ab$, $\eta_\A(1)=1_\A\in\A_0$.

We define an analogue of direct sum for the coactions $\delta$ and $\gamma$ as a coaction of $\OO$ on $\B\otimes\C$ by using the monoidal category $(\Alg,\otimes)$. According to p.~\ref{bbActBimon} we need to introduce the structure of a bimonoid on the comonoid $(\A,\Delta_\A,\varepsilon_\A)\in\Comon(\Alg,\otimes)$, this is a pair $(\mu'_\A,\eta'_\A)$ such that $(\A,\mu'_\A,\eta'_\A,\Delta_\A,\varepsilon_\A)$ is a bimonoid in $(\Alg,\otimes)$, where $\A$ is understood as the algebra $(\A,\mu_\A,\eta_\A)\in\Alg$. But $(\Alg,\otimes)=\Mon(\Vect,\otimes)$, so due to Prop.~\ref{PropMonMon} we obtain
\begin{multline*}
 \Bimon(\Alg,\otimes)=\Comon\big(\Mon(\Alg,\otimes)\big)= \Comon\Big(\Mon\big(\Mon(\Vect,\otimes)\big)\Big)= \\
 \Comon\big(\cMon(\Vect,\otimes)\big)=\Comon(\CommAlg,\otimes).
\end{multline*}
This implies $(\A,\mu'_\A,\eta'_\A)=(\A,\mu_\A,\eta_\A)\in\CommAlg$, so we can define the direct sum of coactions by the formula~\eqref{deltagamma} iff the algebra $\A$ is commutative. We need to use in this case the multiplication $\mu_\A$ in~\eqref{deltagamma}. In the non-commutative case we define the direct sum in an immediate way, but not for every pair $(\delta,\gamma)$.

Consider the linear map
\begin{align} \label{deltaotimes}
 \B\otimes\C\xrightarrow{\delta_\A\otimes\gamma_\A}
 \A\otimes\B\otimes\A\otimes\C\xrightarrow{\sigma^{(23)}}
 \A\otimes\A\otimes\B\otimes\C\xrightarrow{\mu_\A\otimes\id\otimes\id}
 \A\otimes\B\otimes\C.
\end{align}
The conditions~\eqref{deltaAk} and $\mu_\A(\A_k\otimes\A_l)\subset\A_{k+l}$ imply that the image of $\B_k\otimes\C_l$ under the map~\eqref{deltaotimes} lies in $\A_{k+l}\otimes\B_k\otimes\C_l$, so this map also satisfies~\eqref{deltaAk}. Hence it gives the graded map $\B\otimes\C\to\A\circ(\B\otimes\C)$, which we denote by $\delta\dotplus\gamma$.

\begin{Def} \normalfont
We say that the direct sum of the coactions $\delta,\gamma$ {\it exists} if $\delta\dotplus\gamma$ is an algebra homomorphism. In this case $\delta\dotplus\gamma$ is a coaction of $\OO$ on $\B\otimes\C$ called {\it direct sum of the coactions $\delta$ and $\gamma$}.
Suppose $\B,\C\in\FQA$, so we can consider the corresponding corepresentations $\omega=\vartheta^{-1}(\delta)\colon \bfcoend(\B)\to\OO$ and $\nu=\vartheta^{-1}(\gamma)\colon\bfcoend(\C)\to\OO$.
If the direct sum of $\delta,\gamma$ exists, then
\begin{align*}
 \omega\dotplus\nu:=\vartheta^{-1} (\delta\dotplus\gamma)\colon\bfcoend(\B\otimes\C)\to\OO
\end{align*}
is a corepresentation of the comonoid $\OO$ on the object $\B\otimes\C$ called {\it direct product of corepresentations $\omega$ and $\nu$}. In this case we say that the direct sum of corepresentations $\omega,\nu$ {\it exists}.
\end{Def}

\begin{Prop} \label{PropDiS}
Let $\B=\gX_B(\KK)$ and $\C=\gX_C(\KK)$, $D=\DiS(B,C)$. Let $\delta$ and $\gamma$ be coactions of $\OO$ given by multiplicative first order $B$- and $C$-Manin matrices $M$ and $N$ respectively.
The direct sum of the coactions $\delta,\gamma$ exists iff $M$ and $N$ entry-wise commute in the sense of~\eqref{MNComm}. Thus the direct sum of the corepresentations $\omega,\nu$ exists iff the corresponding Manin matrices $M$ and $N$ satisfy this condition. In this case the coaction $\delta\dotplus\gamma\colon\gX_D(\KK)\to\A\circ\gX_D(\KK)$ and corepresentation $\omega\dotplus\nu\colon\bfcoend\big(\gX_D(\KK)\big)\to\OO$ are given by the multiplicative first order $D$-Manin matrix~\eqref{MoplusN}.
\end{Prop}

\noindent{\bf Proof.} Since $\B_0=\C_0=\KK$ we have $(\B\otimes\C)_1=\B_1\oplus\C_1$, so the first order component of $\delta\dotplus\gamma$ coincides with $\delta_1\oplus\gamma_1$. Hence the graded map $\delta\dotplus\gamma\colon\B\otimes\C\to\A\circ(\B\otimes\C)$ is given by the matrix $L=M\oplus N$. Since $M$ and $N$ are multiplicative first order matrices, so is $L$. The rest follows from Prop.~\ref{PropCoa} and~\ref{PropLMN}. \qed

In general, the operation of direct sum does not give a monoidal structure on $\Lcoact(\OO)$ nor on $\Corep_\FQA(\OO)$. However if $\delta,\gamma,\delta',\gamma'$ are coactions of $\OO$ on $\B$, $\C$, $\B'$, $\C'$ such that the direct sums of the coactions $\delta,\gamma$ and of $\delta',\gamma'$ exist, then due to Prop.~\ref{PropDiS} and~\ref{PropcoactionMor} any morphisms $f\colon(\B,\delta)\to(\B',\delta')$ and $g\colon(\C,\gamma)\to(\C',\gamma')$ in $\Lcoact(\OO)$ give a morphism
 $$f\dotplus g\colon(\B\otimes\C,\delta\dotplus\gamma)\to(\B'\otimes\C',\delta'\dotplus\gamma').$$
Thus we have the functor from a full subcategory of $\Lcoact(\OO)\times\Lcoact(\OO)$ to $\Lcoact(\OO)$. The objects of this subcategory are the pairs of coactions for which the direct sum exists. Note that the direct sum of the coactions $\delta,\gamma$ exists iff it exists for $\gamma,\delta$; in this case the isomorphism $\sigma_{\B,\C}\colon\B\otimes\C\cong\C\otimes\B$ gives the isomorphism $(\B\otimes\C,\delta\dotplus\gamma)\cong(\C\otimes\B,\gamma\dotplus\delta)$ in $\Lcoact(\OO)$. In particular, if the algebra $\A$ is commutative, then the direct sum functor $(\B,\delta)\dotplus(\C,\gamma)=(\B\otimes\C,\delta\dotplus\gamma)$ defines a symmetric monoidal structure on the category $\Lcoact(\OO)$ and on its subcategory $\Rep_\FQA(\OO)$.

\bb{Categorical product of quantum representations.} In p.~\ref{bbCop} we constructed the coproduct of quadratic algebras, which corresponds to the product of quantum linear spaces in the category $\FQA^\op$. Let us modify the definition of direct sum of corepresentations for the case of coproduct $\B\amalg\C$.

Let $\gX_B(\KK)$ and $\gX_C(\KK)$ be quadratic algebras for idempotents $B\in\End(\KK^m\otimes\KK^m)$ and $C\in\End(\KK^n\otimes\KK^n)$. Define the idempotent $E=\CoP(B,C)\in\End(\KK^{m+n}\otimes\KK^{m+n})$ with the only non-zero entries
\begin{align}
 &E^{ij}_{kl}=B^{ij}_{kl}, && E^{\bar a,\bar b}_{\bar c,\bar d}=C^{ab}_{cd},
\end{align}
where $i,j,k,l=1,\ldots,m$,\quad $a,b,c,d=1,\ldots,n$, \quad $\bar a=m+a$ (the entries
 $E^{i,\bar a}_{i,\bar a}$, $E^{\bar a,i}_{\bar a,i}$,
$E^{i,\bar a}_{\bar a,i}$, $E^{i,\bar a}_{\bar a,i}$ are also zero in this case).
Then $\gX_B(\KK)\amalg\gX_C(\KK)=\gX_E(\KK)$.

\begin{Prop} \label{PropCoP}
 Consider idempotents $B'\in\End(\KK^{m'}\otimes\KK^{m'})$, $C'\in\End(\KK^{n'}\otimes\KK^{n'})$ and $E'=\CoP(B',C')\in\End(\KK^{m'+n'}\otimes\KK^{m'+n'})$.
 Let $M$ and $N$ be $m\times m'$ and $n\times n'$ matrices over $\gR$. Let $L=M\oplus N$ be defined as in Prop.~\ref{PropLMN}. Then $L$ is a $(E,E')$-Manin matrix iff $M$ is $(B,B')$-Manin matrix and $N$ is $(C,C')$-Manin matrix.
\end{Prop}

\noindent{\bf Proof.} The same as for Prop.~\ref{PropLMN}. \qed

Note that in this case the commutativity condition~\eqref{MNComm} is not needed. Therefore we can define corepresentation of $\OO$ on $\B\amalg\C$ for any corepresentations $\omega,\nu$ on $\B,\C\in\FQA$.

\begin{Def} \normalfont
 Let $\B=\gX_B(\KK)$ and $\C=\gX_C(\KK)$. Then we have $\B\amalg\C=\gX_E(\KK)$ for $E=\CoP(B,C)$. Consider corepresentations $\omega\colon\bfcoend(\B)\to\OO$ and $\nu\colon\bfcoend(\C)\to\OO$ given by multiplicative first order $B$- and $C$-Manin matrices $M$ and $N$ over $\A$. By taking into account Prop.~\ref{PropCoP} we see that $L=M\oplus N$ is a multiplicative first order $E$-Manin matrix over $\A$. We call the corepresentation $\omega\dotsqcup\nu\colon\bfcoend(\B\amalg\C)\to\OO$ given by the matrix $L$ {\it coproduct of corepresentations} $\omega$ and $\nu$. Let $\delta=\vartheta(\omega)$, $\gamma=\vartheta(\nu)$ be coactions of $\OO$ on $\B$ and $\C$ given by the matrices $M$ and $N$. By {\it coproduct of coactions} $\delta$ and $\gamma$ we mean the coaction $\delta\dotsqcup\gamma=\vartheta(\omega\dotsqcup\nu)$ on $\B\amalg\C$ given by the matrix $L$.
\end{Def}

Note that the coproduct of corepresentations (coactions) is given by the same matrix $L$ as the their direct sum, but this is a corepresentation (coaction) on another quadratic algebra (see the warning in p.~\ref{bbMultMM}).

By using the results of p.~\ref{bbCop} and Prop.~\ref{PropcoactionMor} one can check that $(\B\amalg\C,\omega\dotsqcup\nu)$ is a coproduct of the objects $(\B,\omega)$ and $(\C,\nu)$ in the category $\Corep_\FQA(\OO)$. Hence the latter is a category with finite coproducts (its initial object is the algebra $\KK$ with the trivial corepresentation). If we denote this coproduct by $(\B,\omega)\dotsqcup(\C,\nu):=(\B\amalg\C,\omega\dotsqcup\nu)$, then we obtain a bifunctor $-\dotsqcup-\colon\Corep_\FQA(\OO)\times\Corep_\FQA(\OO)\to\Corep_\FQA(\OO)$. The coproduct of corepresentations corresponds to the categorical product of quantum representations in $\Corep_\FQA(\OO)^\op$.

\bb{Tensor product of Manin matrices.} The quantum analogue of the tensor product is the Manin product (see p.~\ref{bbSQVS} and p.~\ref{bbLaxMonS}). Let us describe the Manin product of the connected finitely generated quadratic algebras $\B=\gX_B(\KK)$ and $\C=\gX_C(\KK)$, where $B\in\End(\KK^m\otimes\KK^m)$ and $C\in\End(\KK^n\otimes\KK^n)$ are matrix idempotents.

\begin{Prop}
 The operator $F=\TeP(B,C)\in\End(\KK^m\otimes\KK^n\otimes\KK^m\otimes\KK^n)$ defined as
\begin{align} \label{FAB}
 &F=\sigma^{(23)}(B\otimes1+1\otimes C-B\otimes C)\sigma^{(23)},
\end{align}
is an idempotent.
We have $\gX_B(\KK)\circ\gX_C(\KK)=\gX_F(\KK)$.
\end{Prop}

\noindent{\bf Proof.} The idempotentness of~\eqref{FAB} is equivalent to the idempotentness of the operator
\begin{align} \label{wtF}
 \wt F=\sigma^{(23)}F\sigma^{(23)}=B\otimes1+1\otimes C-B\otimes C\in\End(\KK^m\otimes\KK^m\otimes\KK^n\otimes\KK^n),
\end{align}
which, in turn, directly follows from $B^2=B$ and $C^2=C$. Let $W=(\KK^m)^*$ and $V=(\KK^n)^*$, then we have $\gX_B(\KK)=TW/(\Im B^*)$ and $\gX_C(\KK)=TV/(\Im C^*)$, where $B^*\in\bfend(W^{\otimes2})$ and $C^*\in\bfend(V^{\otimes2})$ are transposed idempotents (see~p.~\ref{bbIdem}). The algebra $\gX_B(\KK)\circ\gX_C(\KK)$ is the quotient $T(W\otimes V\big)/(R_{\mathrm w})$, where $R_{\mathrm w}=\sigma^{(23)}(\Im B^*\otimes V^{\otimes2}+W^{\otimes2}\otimes\Im C^*)$. Hence it is enough to check that the subspace $\Im B^*\otimes V^{\otimes2}+W^{\otimes2}\otimes\Im C^*=\Im(B^*\otimes1)+\Im(1\otimes C^*)\subset W^{\otimes2}\otimes V^{\otimes2}$ coincides with $\sigma^{(23)}\Im F^*=\Im\wt F^*$. We have $\Im\wt F^*=\Im(B^*\otimes1+1\otimes C^*-B^*\otimes C^*)\subset\Im(B^*\otimes1)+\Im(1\otimes C^*)$. On the other hand, the formulae $B^*\otimes1=\wt F(B^*\otimes1)$ and $1\otimes C^*=\wt F(1\otimes C^*)$ imply $\Im(B^*\otimes1)\subset\Im(\wt F^*)$ and $\Im(1\otimes C^*)\subset\Im(\wt F^*)$. \qed

Consider two more idempotents $B'\in\End(\KK^{m'}\otimes\KK^{m'})$, $C'\in\End(\KK^{n'}\otimes\KK^{n'})$. Define $F'=\DiS(B',C')\in\End(\KK^{m'+n'}\otimes\KK^{m'+n'})$. We have $F'=\sigma^{(23)}\wt F'\sigma^{(23)}$ with
\begin{align} \label{wtFp}
 \wt F'=B'\otimes1+1\otimes C'-B'\otimes C'.
\end{align}
For graded homomorphisms $f\colon\gX_B(\KK)\to\gX_{B'}(\gR)$ and $g\colon\gX_{C}(\KK)\to\gX_{C'}(\gR)$ consider a graded operator
\begin{align}
 \gX_F(\KK)=\gX_B(\KK)\circ\gX_C(\KK)\xrightarrow{f\circ g}
 \gX_{B'}(\gR)\circ\gX_{C'}(\gR)\xrightarrow{(\mu_\gR)} \gX_{F'}(\gR), \label{grOpTeP}
\end{align}
where $\gX_{B'}(\gR)\circ\gX_{C'}(\gR)\xrightarrow{(\mu_\gR)}\gX_{D'}(\gR)$ is the multiplication in the tensor factors $\gR$, that is a graded operator with the components
\begin{multline*}
 \gX_{B'}(\gR)_k\otimes\gX_{C'}(\gR)_k=
\gR\otimes\gX_{B'}(\KK)_k\otimes\gR\otimes\gX_{C'}(\KK)_k
\xrightarrow{\sigma^{(23)}} \\
\gR\otimes\gR\otimes\gX_{B'}(\KK)_k\otimes\gX_{C'}(\KK)_k
\xrightarrow{\mu_\gR\otimes\id\otimes\id}
\gR\otimes\gX_{B'}(\KK)_k\otimes\gX_{C'}(\KK)_k=
\gX_{F'}(\gR)_k.
\end{multline*}

The homomorphisms $f$ and $g$ are given by some $(B,B')$- and $(C,C')$-Manin matrices $M$ and $N$. These are $m\times m'$ and $n\times n'$ matrices over $\gR$. The graded operator~\eqref{grOpTeP} is given by the matrix $P=M^{(1)}N^{(2)}$, this is an $(mn)\times(m'n')$ matrix over $\gR$ with the entries
\begin{align}
 P^{ia}_{jb}=M^i_jN^a_b.
\end{align}
In the notations of p.~\ref{bbQSLS} this is the operator $P=M\dototimes N \colon\KK^{m'}\otimes\KK^{n'}\to\gR\otimes\KK^m\otimes\KK^n$.

\begin{Prop} \label{PropPMN}
Suppose $M$ and $N$ are $(B,B')$- and $(C,C')$-Manin matrices. The matrix $P=M^{(1)}N^{(2)}=M\dototimes N$ is an $(F,F')$-Manin matrix iff
\begin{align} \label{MNFComm}
 \wt FM^{(1)}(M^{(2)}N^{(3)}-N^{(3)}M^{(2)})N^{(4)}(1-\wt F')=0,
\end{align}
where $\wt F$ and $\wt F'$ are given by the formulae~\eqref{wtF}, \eqref{wtFp}.
Let $f=f_M$ and $g=f_N$, then~\eqref{MNFComm} is the condition for the graded operator~\eqref{grOpTeP} to be an algebra homomorphism.
\end{Prop}

\noindent{\bf Proof.} The matrix $P$ is an $(F,F')$-Manin matrix iff
 $FM^{(1)}N^{(2)}M^{(3)}N^{(4)}(1-F')=0$. By means of the conjugation by $\sigma^{(23)}$ we obtain the equivalent condition
\begin{align} \label{F32}
 \wt FM^{(1)}N^{(3)}M^{(2)}N^{(4)}(1-\wt F')=0.
\end{align}
Due to $BM^{(1)}M^{(2)}=BM^{(1)}M^{(2)}B'$ and $CN^{(1)}N^{(2)}=CN^{(1)}N^{(2)}C'$ we derive
\begin{align*}
 (B\otimes 1)M^{(1)}M^{(2)}N^{(3)}N^{(4)}&=(B\otimes 1)M^{(1)}M^{(2)}N^{(3)}N^{(4)}(B'\otimes1), \\
\big((1-B)\otimes C\big)M^{(1)}M^{(2)}N^{(3)}N^{(4)}&=\big((1-B)\otimes C\big)M^{(1)}M^{(2)}N^{(3)}N^{(4)}(1\otimes C').
\end{align*}
Addition and multiplication by $1-\wt F'=(1-B')\otimes(1-C')$ on the right gives
\begin{align} \label{F23}
 \wt FM^{(1)}M^{(2)}N^{(3)}N^{(4)}(1-\wt F')=0.
\end{align}
By subtracting~\eqref{F32} from \eqref{F23} we get~\eqref{MNFComm}. \qed

Note that the condition~\eqref{MNFComm} is valid for entry-wise commuting matrices $M$ and $N$ in the sense of Prop.~\ref{PropLMN}. However, the condition~\eqref{MNComm} is not necessary. For example, if $B,B',C,C'$ are $0$, then any matrices $M$ and $N$ satisfy~\eqref{MNFComm}.

\bb{Tensor product of quantum representations.}
Let $\OO=(\A,\Delta,\varepsilon)$ be a comonoid in $(\GrAlg,\circ)$, where $\A=(\A,\mu_\A,\eta_\A)\in\GrAlg$. By tensor product of its coactions or corepresentations on $\B$ and $\B'$ we mean the corresponding coaction or, respectively, corepresentation on $\B\circ\B'$. One can define it as a monoidal product in the category $\Lcoact(\OO)$ introduced in p.~\ref{bbActBimon} if $\OO$ has additionally a structure of bimonoid $\BB=(\A,\mu,\eta,\Delta,\varepsilon)$, where $\mu\colon\A\circ\A\to\A$ and $\eta\colon\KK[u]\to\A$ are morphisms of graded algebras compatible with the comonoid structure of $\OO$ in the sense of p.~\ref{bbBimon} (the homomorphisms $\mu$, $\eta$ are not to be confused with the maps $\mu_\A\colon\A\otimes\A\to\A$ and $\eta_\A\colon\KK\to\A$). This means the following conditions on the components $\mu_k\colon\A_k\otimes\A_k\to\A_k$ and $\eta_k\colon\KK\to\A_k$. First, they must give an associative multiplication on each component $\A_k$ with the unities $1_k=\eta_k(1)=\eta(u^k)\in\A_k$ such that $(\A_k,\mu_k,\eta_k,\Delta_k,\varepsilon_k)\in\Bimon(\Vect,\otimes)$, where $\Delta_k\colon\A_k\to\A_k\otimes\A_k$ and $\varepsilon_k\colon\A_k\to\KK$ are the graded components of $\Delta$ and $\varepsilon$. Second, they must be compatible with the multiplication $\mu_\A(a\otimes b)=ab$ and with the unity $1_\A=\eta_\A(1)$, that is $\mu_{k+l}(a_kb_l\otimes a'_kb'_l)=\mu_k(a_k\otimes a'_k)\mu_l(b_l\otimes b'_l)$ for any $a_k,a'_k\in\A_k$, $b_l,b'_l\in\A_l$ and $1_{k+l}=1_k1_l$, $1_0=1_\A$. In particular, the vector space $\A_0$ has two multiplications $(\mu_A)_0$ and $\mu_0$ with the same unity. By Eckmann--Hilton Principle (Prop.~\ref{PropMonMon}) the compatibility condition $\mu_0(ab\otimes cd)=\mu_0(a\otimes c)\mu_0(d\otimes d)$ implies that $\mu_0=(\mu_\A)_0$ and this is a commutative multiplication on $\A_0$. Moreover, one can show that a monoid structure $(\mu,\eta)$ on a graded algebra $\gR\otimes S(V)$ exists only if $\dim V\le 1$.

Consider a bialgebra $(\gR,\Delta_\gR,\varepsilon_\gR)\in\Comon(\Alg,\otimes)$. The embedding $\Alg\hookrightarrow\QAsc$ gives the comonoid $\OO_\gR=(\gR\otimes\KK[u],\Delta,\varepsilon)$ defined by the formula~\eqref{CCgR}.

If $\gR\in\CommAlg$, then the comonoid $\OO_\gR$ has a unique bimonoid structure:
\begin{align} \label{BBgR}
 &\BB_\gR=(\gR\otimes\KK[u],\mu,\eta,\Delta,\varepsilon), &&&
 &\mu=\mu_\gR\otimes\id_{\KK[u]}, &
 &\eta=\eta_\gR\otimes\id_{\KK[u]}.
\end{align}
In this case we get a tensor product of coactions as the monoidal product in $\Lcoact(\BB_\gR)$.

Here we have a situation similar to one from p.~\ref{bbDS}: we can define the tensor product without commutativity condition on $\mu_\gR$, but this product not always exists.
Namely, consider an arbitrary $\gR\in\Alg$. Note that $(\gR\otimes\KK[u])\circ\gX_B(\KK)=\gX_B(\gR)$. Let $\delta\colon\gX_B(\KK)\to\gX_B(\gR)$, $\gamma\colon\gX_C(\KK)\to\gX_C(\gR)$ be coactions of $\OO_\gR$ on $\B=\gX_B(\KK)$ and $\C=\gX_C(\KK)$ given by multiplicative $B$- and $C$-Manin matrices $M$ and $N$ over $\gR$. By substituting the homomorphisms $f=f_M=\delta$ and $g=f_N=\gamma$ to~\eqref{grOpTeP} we obtain the graded operator
\begin{align} \label{Fdelta}
 \gX_F(\KK)=\gX_B(\KK)\circ\gX_C(\KK)\xrightarrow{\delta\circ\gamma}
 \gX_{B}(\gR)\circ\gX_{C}(\gR)\xrightarrow{(\mu_\gR)} \gX_{F}(\gR),
\end{align}
where $F=\TeP(B,C)=\sigma^{(23)}(B\otimes1+1\otimes C-B\otimes C)\sigma^{(23)}$. Denote the composition~\eqref{Fdelta} by $\delta\dotcirc\gamma$. Due to Prop.~\ref{PropPMN} this operator is a coaction of $\OO_\gR$ on $\B\circ\C$ iff the matrices $M$ and $N$ satisfy~\eqref{MNFComm}.

\begin{Def} \normalfont
We say that the tensor product of the coactions $\delta,\gamma$ {\it exists} if $\delta\dotcirc\gamma$ is an algebra homomorphism, i.e. if $$\delta\dotcirc\gamma\colon\B\circ\C\to\A\circ\B\circ\C=\gR\otimes\B\circ\C$$
is a coaction of $\OO_\gR$ on $\B\circ\C$, we call it {\it tensor product of the coactions $\delta$ and $\gamma$}.
Let $\omega=\vartheta^{-1}(\delta)\colon\bfcoend(\B)\to\OO_\gR$, $\nu=\vartheta^{-1}(\gamma)\colon\bfcoend(\C)\to\OO_\gR$ be the corresponding corepresentations of $\OO_\gR$. If $\delta\dotcirc\gamma$ is a coaction, then $$\omega\dotcirc\nu:=\vartheta^{-1}(\delta\dotcirc\gamma)\colon\bfcoend(\B\circ\C)\to\OO_\gR$$
is a corepresentation of $\OO_\gR$ on $\B\circ\C$ called {\it tensor product of the corepresentations $\omega$ and $\nu$}. In this case we say that the tensor product of corepresentations $\omega,\nu$ {\it exists}.
\end{Def}

It follows from Prop.~\ref{PropMorCorep} again, that the tensor product of corepresentations is a functor $-\dotcirc-$ from a full subcategory of $\Corep_\FQA(\OO_\gR)\times\Corep_\FQA(\OO_\gR)$ to $\Corep_\FQA(\OO_\gR)$. It has the form $(\B,\omega)\dotcirc(\C,\nu)=(\B\circ\C,\omega\dotcirc\nu)$. If $\gR\in\CommAlg$, this functor is defined on the whole $\Corep_\FQA(\OO_\gR)\times\Corep_\FQA(\OO_\gR)$ and, consequently, gives a structure of monoidal category $\big(\Corep_\FQA(\OO_\gR),\dotcirc\big)$. It is  symmetric due to commutativity of $\mu_\gR$.

Similarly we can define a tensor product of the quantum representations $\omega$ and $\nu$ on the black Manin product $\B\bullet\C=\gX_G(\KK)$, where $G=\sigma^{(23)}(B\otimes C)\sigma^{(23)}$ is the corresponding idempotent. Since $(\B\bullet\C)_1=(\B\circ\C)_1$, a corepresentation on the black Manin product can be defined by the same matrix $P=M\dot\otimes N=M^{(1)}N^{(2)}$, which is multiplicative with respect to $\Delta_\gR$. It is a $G$-Manin matrix iff
\begin{align} \label{BCMMNN}
 (B\otimes C)M^{(1)}(M^{(2)}N^{(3)}-N^{(3)}M^{(2)})N^{(4)}(1-B\otimes C)=0.
\end{align}
Under this condition we obtain {\it black tensor product}
\begin{align} \label{bpCorep}
 \omega\dotbullet\nu
\colon\bfcoend(\B\bullet\C)\to\OO_\gR
\end{align}
defined by the matrix $M\dot\otimes N$.
This is a corepresentation of $\OO_\gR$ on the quadratic algebra $\B\bullet\C$.
We get one more partially defined binary operation with quantum representations. In the case $\gR\in\CommAlg$, it is a bifunctor $-\dotbullet-\colon\Corep_\FQA(\OO_\gR)\times\Corep_\FQA(\OO_\gR)\to\Corep_\FQA(\OO_\gR)$.

\bb{Opposite and coopposite quantum representations.}
For an arbitrary comonoid $\OO=(\A,\Delta,\varepsilon)\in\Comon(\GrAlg,\circ)$ denote $\OO^\op=(\A^\op,\Delta,\varepsilon)$, where $\A^\op$ is the opposite graded algebra (the opposite monoid in $(\GrVect,\otimes)$ in the sense of p.~\ref{bbOpMon}). The corresponding bialgebra $(\A^\op,\Delta_\A,\varepsilon_\A)$ is the bimonoid in $(\Vect,\otimes)$ opposite to the bimonoid $(\A,\Delta_\A,\varepsilon_\A)$. This bialgebra satisfies~\eqref{DeltaAk}, so the triple $\OO^\op=(\A^\op,\Delta,\varepsilon)$ is also a comonoid in $(\GrAlg,\circ)$.

The opposite graded algebra to a quadratic algebra $\A=TV/(R)\in\QA$ has the form $\A^\op=TV/(\sigma_{V,V}R)$ (see~\cite{Manin88}). Let $\B=\gX_B(\KK)\in\FQA$ for some $B\in\bfend(\KK^m\otimes\KK^m)$ and denote $B^{(21)}=\sigma\cdot B\cdot\sigma\in\bfend(\KK^m\otimes\KK^m)$, where $\sigma=\sigma_{\KK^m,\KK^m}$. This is an $m^2\times m^2$ matrix over $\KK$ with entries $(B^{(21)})^{ij}_{kl}=B^{ji}_{lk}$. The corresponding quadratic algebra is $\B^\op=\gX_{B^{(21)}}(\KK)$. Let $\M=(\M^i_j)$ be the universal $B$-Manin matrix. The algebra $\bfcohom(\B^\op,\B^\op)$ is generated by the same $\M^i_j$, but with the relations $B^{(21)}\M^{(1)}\M^{(2)}(1-B^{(21)})=0$. In other words, the matrix $\M$ considered as a matrix over $\bfcohom(\B^\op,\B^\op)$ is the universal $B^{(21)}$-Manin matrix. On the other hand the algebra $\bfhom(\B,\B)^\op$ is generated by $\M^i_j$ with the relations $B\M^{(2)}\M^{(1)}(1-B)=0$, which is equivalent to the previous relations. Thus we obtain $\bfcoend(\B^\op)=\bfcoend(\B)^\op$ for any $\B\in\FQA$.

Let $\omega\colon\bfcoend(\B)\to\OO$ be a corepresentation. The same linear map is a comonoid morphism $\bfcoend(\B^\op)=\bfcoend(\B)^\op\to\OO^\op$, so this is a corepresentation of $\OO^\op$ on $\B^\op$. Let us denote it by $\omega^\op$ and call {\it the corepresentation opposite to $\omega$.} Note that the corepresentations $\omega$ and $\omega^\op$ is given by the same multiplicative $B$-Manin matrix $M$, which is considered as $B$- or as $B^{(21)}$-Manin matrix over $\A$ or over $\A^\op$ respectively.

We have the algebra isomorphism $\bfcohom(\B,\B)=\B^!\bullet\B=(\B^!)^!\bullet\B^!=\bfcohom(\B^!,\B^!)$, so the formulae~\eqref{dBcop} implies the comonoid isomorphism $\bfcoend(\B^!)=\bfcoend(\B)^\cop$ (see~\cite[\S~5.10]{Manin88}).
 Any corepresentation $\omega\colon\bfcoend(\B)\to\OO$ can be regarded as a corepresentation $\omega^\cop\colon\bfcoend(\B^!)\to\OO^\cop$ of the coopposite comonoid $\OO^\cop$ on the Koszul dual quadratic algebra $\B^!=\Xi_B(\KK)$, we call it {\it coopposite corepresentation} to the corepresentation $\omega$. It corresponds to the $(1-B^\top)$-Manin matrix $M^\top$ over $\A$, which is multiplicative with respect to $\Delta^\cop$.

By applying these two operations (in any order) to the corepresentation $\omega$ we obtain opposite coopposite corepresentation $\omega^{\op,\cop}\colon\bfcoend\big((\B^!)^\op\big)\to\OO^{\op,\cop}$. This is a corepresentation of the comonoid $\OO^{\op,\cop}=(\A^\op,\Delta^\cop,\varepsilon)$, which corresponds to the opposite coopposite bialgebra $(\A^\op,\Delta_\A^\cop,\varepsilon_\A)$, on the quadratic algebra $(\B^!)^\op=(\B^\op)^!=\Xi_{B^{(21)}}(\KK)$.

\bb{Dual and Koszul dual quantum representations.} \label{bbDualQR}
Assume that the bialgebra $\gR=(\gR,\mu_\gR,\eta_\gR,\Delta_\gR,\varepsilon_\gR)$ is a Hopf algebra and let $\zeta_\gR\colon\gR\to\gR$ be its antipode. This is an algebra anti-homomorphism as well as a coalgebra anti-homomorphism in the sense
\begin{align}
 &\zeta_\gR\cdot\mu_\gR=\mu_\gR^\op\cdot(\zeta_\gR\otimes\zeta_\gR), && \zeta_\gR\cdot\eta_\gR=\eta_\gR, \\
 &(\zeta_\gR\otimes\zeta_\gR)\cdot\Delta_\gR=\Delta_\gR^\cop\cdot\zeta_\gR, && \varepsilon_\gR\cdot\zeta_\gR=\varepsilon_\gR
\end{align}
(see~\cite[Th.~III.3.4]{Kass}). This means that $\zeta_\gR$ is a morphism $\gR^{\op,\cop}\to\gR$ in $\Bimon(\Vect,\otimes)$, where $\gR^{\op,\cop}=(\gR,\mu_\gR^\op,\eta_\gR,\Delta_\gR^\cop,\varepsilon_\gR)$. By tensoring with $\KK[u]$ we obtain the morphism $\zeta=\zeta_\gR\otimes\id_{\KK[u]}\colon\OO_\gR^{\op,\cop}\to\OO_\gR$, where $\OO_\gR=(\gR\otimes\KK[u],\Delta=\Delta_\gR\otimes\id_{\KK[u]},\varepsilon=\varepsilon_\gR\otimes\id_{\KK[u]})$ is the bialgebra $\gR$ embedded to $\Comon(\GrAlg,\circ)$ and $\OO_\gR^{\op,\cop}=(\gR^\op\otimes\KK[u],\Delta^\cop,\varepsilon)$.

Let $\omega\colon\bfcoend(\B)\to\OO_\gR$ be a corepresentation of the comonoid $\OO_\gR$ on a quadratic algebra $\B=\gX_B(\KK)\in\FQA$ given by a multiplicative $B$-Manin matrix $M=(M^i_j)$ over $\gR$. Then the composition
\begin{align} \label{Dualcorep}
 \omega^D\colon\bfcoend\big((\B^!)^\op\big)\xrightarrow{\omega^{\op,\cop}} \OO_\gR^{\op,\cop}\xrightarrow{\zeta}\OO_\gR
\end{align}
is a corepresentation of the same comonoid $\OO_\gR$ on the opposite Koszul dual quadratic algebra $(\B^!)^\op=\Xi_{B^{(21)}}(\KK)$. Let us call it {\it dual corepresentation} to the corepresentation $\omega$. This is a quantum analogue of dual (contragredient) representation on a dual vector space.

Since $M$ is a multiplicative matrix over a Hopf algebra, it is invertible. The entries of the inverse matrix $M^{-1}$ are $(M^{-1})^i_j=\zeta_\gR(M^i_j)$. Due to the fact that $\zeta_\gR$ is an anti-endomorphism of $\gR$ the matrix $M^{-1}$ is a $B^{(21)}$-Manin matrix. The dual corepresentation~\eqref{Dualcorep} is given by the $(1-B^{(21)})^\top$-Manin matrix $(M^\top)^{-1}=(M^{-1})^\top$ over $\gR$, which is multiplicative with respect to the comultiplication $\Delta_\gR$.

Consider the corepresentation opposite to~\eqref{Dualcorep}. It has the form
\begin{align} \label{rhoKDual}
 \omega^{KD}\colon\bfcoend(\B^!)\xrightarrow{\omega^{\cop}} \OO_\gR^{\cop}\xrightarrow{\zeta}\OO_\gR^{\op}.
\end{align}
This is a corepresentation of $\OO_\gR^\op$ on the Koszul dual quadratic algebra $\B^!$. Let us call it {\it Koszul dual corepresentation} to the corepresentation $\omega\colon\bfcoend(\B)\to\OO_\gR$. It corresponds to the $(1-B^\top)$-Manin matrix $(M^\top)^{-1}=(M^{-1})^\top$ over $\gR^\op$.

Let $\nu\colon\bfcoend(\C)\to\OO_\gR$ be a corepresentation of $\OO_\gR$ on $\C\in\FQA$ and $N$ be the corresponding multiplicative matrix. By using~\eqref{KoszulDualExch} we get $\B\bullet\C=(\B^!\circ\C^!)^!$. One can show that the tensor product of $\omega^{KD},\nu^{KD}$ exists iff the black tensor product of $\omega,\nu$ exists (the condition~\eqref{BCMMNN} fulfils). In this case the Koszul dual to $\omega^{KD}\dotcirc\nu^{KD}$ is the corepresentation $\big(\omega^{KD}\dotcirc\nu^{KD}\big)^{KD}$ of $\OO_\gR$ on the black Manin product $\B\bullet\C$, given by the multiplicative matrix $M^{(1)}N^{(2)}$. It coincides with the corepresentation~\eqref{bpCorep}.

Suppose the Hopf algebra $\gR$ is commutative as an algebra. Then the Koszul dual to a corepresentation $\omega\colon\bfcoend(\B)\to\OO_\gR$ is a corepresentation $\omega^{KD}$ of the same comonoid $\OO_\gR^\op=\OO_\gR$ on the object $\B^!$. The comonoid $\OO_\gR$ can be considered as a commutative bimonoid $\BB_\gR=(\gR\otimes\KK[u],\mu=\mu_\gR\otimes\id_{\KK[u]},\eta=\eta_\gR\otimes\id_{\KK[u]},\Delta,\varepsilon)$, it is a Hopf monoid with the antipode~$\zeta$.
By virtue of the formulae~\eqref{KoszulDualExch}, \eqref{cohomBA} we derive $\bfcohom(\B,\C)=\B^!\bullet\C=(\B\circ\C^!)^!$, so the homomorphism
\begin{align} \label{cohomCorep}
 \la\omega,\nu\ra:=\omega^{KD}\dotbullet\nu= \big(\omega\dotcirc\nu^{KD}\big)^{KD}\colon\bfcoend\big(\bfcohom(\B,\C)\big)\to\OO_\gR
\end{align}
is a corepresentation of $\OO_\gR$ on the internal cohom of the corresponding quadratic algebras. It is given by the matrix $\big((M^{-1})^\top\big)^{(1)}N^{(2)}$.

\begin{Prop}
 If the Hopf algebra $\gR$ is commutative, then the symmetric monoidal category $\Corep_\FQA(\BB_\gR)$ is coclosed and $\bfcohom\big((\B,\omega),(\C,\nu)\big)$ is the pair consisting of the object $\bfcohom(\B,\C)$ and the corepresentation~\eqref{cohomCorep}.
\end{Prop}

\noindent{\bf Proof.} Let $\wt\omega\colon\bfend(\wt\B)\to\OO_\gR$ be a corepresentation of the comonoid $\OO_\gR$ on a quadratic algebra $\wt\B\in\FQA$ given by a multiplicative matrix $\wt M$. Due to the definition of internal cohom we have the adjunction
\begin{align} \label{HomBBB}
 \Hom(\bfcohom(\B,\C),\wt\B)\cong\Hom(\C,\wt\B\circ\B)
\end{align}
in the category $\FQA$. Let $(w_i)$ and $(w^i)$ be dual bases of the vector spaces $W=(\B_1)^*$ and $W^*=\B_1$. Let $(z^l)$ and $(\wt w^a)$ be bases of $Z^*=\C_1$ and $\wt W^*=\wt\B_1$.
The graded homomorphisms $h\colon\bfcohom(\B,\C)\to\wt\B$ and $f\colon\C\to\wt\B\circ\B$ related to each other by the bijection~\eqref{HomBBB} have the form
\begin{align}
 &h\colon w_i\otimes z^k\mapsto\sum_a K^k_{ia}\wt w^a, &
 &f\colon z^k\mapsto\sum_{i,a} K^k_{ia}\wt w^a\otimes w^i
\end{align}
for the same coefficients $K^k_{ia}\in\KK$. The corepresentation $\la\omega,\nu\ra$ is given by the matrix $\big((M^{-1})^\top\big)^{(1)}N^{(2)}$, hence $h$ is a morphism from corepresentation $\la\omega,\nu\ra$ to the corepresentation $\wt \omega$ iff $\sum_a K^k_{ia}\wt M^a_b=\sum_{j,l}(M^{-1})^j_iN^k_lK^l_{jb}$.
The condition that the homomorphism $f$ is a morphism from $\nu$ to $\wt\omega\circ\omega$ has the form  $\sum_{i,a} K^k_{ia}\wt M^a_bM^i_j=\sum_{l} N^k_lK^l_{jb}$. Due to the equivalence of these conditions the adjunction~\eqref{HomBBB} gives a natural bijection between the morphisms of the corresponding corepresentations. \qed

\subsection{$S$-embedding of classical representations}

Now we `embed' the classical Representation Theory to the Quantum Representation Theory. Let us begin with representations of finite-dimensional algebra on finite-dimensional vector spaces, i.e. with representations in the closed category $(\FVect,\otimes)$. We embed them to the category of representations in $(\FQA^\op,\circ)$ via the functor $S^*\colon\FVect\hookrightarrow\FQA^\op$. Then we generalise this embedding for the monoidal category $(\SLAffSch,\otimes)$, which contains $(\AffSch,\times)$ and $(\AlgSet,\times)$ as (not full) monoidal subcategories.

\bb{Representations of a finite-dimensional algebra as quantum representations.} \label{bbRepFAlg}
According to p.~\ref{bbMonF} the lax monoidal functor $S^*\colon(\FVect,\otimes)\to(\FQA^\op,\circ)$ induces the contravariant functor $\Mon(S^*)\colon\FAlg\to\Comon(\FQA,\circ)$, where we used $\Mon(\FVect,\otimes)=\FAlg$, $\Mon(\FQA^\op,\circ)=\Comon(\FQA,\circ)^\op$ (see~p.~\ref{bbTPMon}, p.~\ref{bbMorMonComon}). Let $\gA=(V,\mu_V,\eta_V)\in\FAlg$. As it is described in p.~\ref{bbFact} and p.~\ref{bbFrep} the actions and representations of the algebra $\gA$ (objects of $\Lact(\gA)=\Rep_\FVect(\gA)$) are embedded contravariantly to $\Lcoact(\OO)=\Corep_\FQA(\OO)$, where $\OO=\Mon(S^*)\gA\in\Comon(\FQA,\circ)$. We call it {\it $S$-embedding} of finite-dimensional representations of $\gA$ into the category of quantum representations. Let us describe this embedding in details.

First we need to examine how the algebra $\gA$ is translated to the quantum level. Recall that $S^*$ is a composition of the contravariant and covariant functors $(-)^*$ and $S$. Since the duality functor $(-)^*\colon(\FVect,\otimes)\to(\FVect^\op,\otimes)$ is a strong monoidal, the dual vector space $V^*\in\FVect$ is equipped with a structure of coalgebra by the maps $\Delta_{V^*}\colon V^*\xrightarrow{\mu_V^*}(V\otimes V)^*\cong V^*\otimes V^*$ and $\varepsilon_{V^*}=\eta_V^*\colon V^*\to\KK$.
 Further, the colax monoidal functor $S\colon(\FVect,\otimes)\to(\FQA,\circ)$ translates the coalgebra $\gA^*=(V^*,\Delta_{V^*},\varepsilon_{V^*})$ to the comonoid $\OO=\Comon(S)\gA^*=(SV^*,\Delta,\varepsilon)$ with the morphisms
 $\Delta\colon SV^*\xrightarrow{S\mu_V^*}
S^*(V\otimes V)\cong S(V^*\otimes V^*)\xrightarrow{\phi_{V^*,V^*}} SV^*\circ SV^*$ and $\varepsilon=S\eta_V^*\colon SV^*\to S\KK=\KK[u]$, where $\phi$ is the colax monoidal structure of $S$.

Let $vv'=\mu_V(v\otimes v')$ be the multiplication in the algebra $\gA$ and $1_\gA=\eta_V(1)$ be its unity. The quadratic algebras $SV^*$ and $SV^*\otimes SV^*\cong S^*(V\oplus V)$ are identified with the algebras of the polynomials on the spaces $V$ and $V\times V=V\oplus V$ considered as affine spaces (their linear structure gives the grading). Then the structure of the bialgebra $(SV^*,\Delta_{SV^*},\varepsilon_{SV^*})$ corresponding to the comonoid $\OO$ has the form
\begin{align}
 &\Delta_{SV^*}\colon SV^*\to SV^*\otimes SV^*, &&(\Delta_{SV^*}p)(v,v')=p(vv'), \label{Deltap} \\
 &\varepsilon_{SV^*}\colon SV^*\to\KK, &&\varepsilon_{SV^*}(p)=p(1_\gA), \label{varepsilonp}
\end{align}
where $p\in SV^*$ is a polynomial function on $V$ and $v,v'\in V$. In terms of dual bases $(v_i)_{i=1}^n$ and $(v^i)_{i=1}^n$ of $V$ and $V^*$ we have
 $\Delta(v^k)=\sum_{i,j=1}^n c_{ij}^k v^i\otimes v^j$, $\varepsilon(v^k)=d^k$ 
for coefficients $c_{ij}^k,d^k\in\KK$, determined from the formulae
 $v_iv_j=\sum_{k=1}^n c_{ij}^kv_k$, $1_\gA=\sum_{k=1}^n d^kv_k$.
In particular, this implies the condition~\eqref{Delta01}.

Consider a representation $\rho\colon\gA\to\bfend(W)$ of the algebra $\gA$ on a finite-dimensional vector space $W\in\FVect$. It corresponds to an action $a\colon V\otimes W\to W$. Due to the natural isomorphism $(V\otimes W)^*\cong V^*\otimes W^*$ the functor $(-)^*\colon\FVect\to\FVect^\op$ gives the coaction $a^*\colon W^*\to V^*\otimes W^*$ of the coalgebra $\gA^*\in\Comon(\FVect)$ on the vector space $W^*\in\FVect$.
Then, by applying $S$ we get the coaction $\delta\colon SW^*\xrightarrow{Sa^*}S(V^*\otimes W^*)\xrightarrow{\phi_{V^*,W^*}} SV^*\circ SW^*$ of the comonoid $\OO$ on the quadratic algebra $SW^*\in\FQA$.

Thus, the lax monoidal functor $S^*$ induces the functor $\Rep_\FVect(\gA)\to\Corep_{\FQA}(\OO)$, which maps $(W,\rho)$ to $(SW^*,\omega)$, where $\omega\colon\bfcoend(SW^*)\to\OO$ is a corepresentation corresponding to the coaction $\delta$. This corepresentation has the form
\begin{align}
 \omega\colon\bfcoend(SW^*)\xrightarrow{\Phi_{W,W}}S\big(\bfend(W)\big)^*\xrightarrow{S\rho^*}\OO, \label{omegarho}
\end{align}
where $\Phi_{W,Z}\colon\bfcohom(SW^*,SZ^*)\to S\big(\bfhom(W,Z)\big)^*$ is the graded homomorphism whose first order component is the isomorphism $\bfhom(W^*,Z^*)\cong W\otimes Z^*\cong(W^*\otimes Z)^*\cong\bfhom(W,Z)^*$ natural in $W\in\FVect$, $Z\in\Vect$.
Thus, the $S$-embedding of representations is the functor
\begin{align} \label{Sembedding}
 &\Rep_\FVect(\gA)\to\Corep_\FQA(\OO)^\op, &&(W,\rho)\mapsto(TW^*,\omega).
\end{align}

Let $(w_i)_{i=1}^m$ and $(w^i)_{i=1}^m$ be dual bases of the vector spaces $W$ and $W^*$ respectively. The latter basis gives the isomorphism $SW^*=\gX_{A_m}(\KK)$. Any representation $\rho\colon\gA\to\bfend(W)$ has the form $\rho(v)w_j=\rho^i_j(v)w_i$ where $\rho^i_j\colon V\to\KK$ are linear functions such that
\begin{align} \label{rhoij}
 &\rho^i_j(vv')=\sum_{k=1}^m\rho^i_k(v)\rho^k_j(v'), &&\rho^i_j(1_\gA)=\delta^i_j.
\end{align}
These functions form an $m\times m$ matrix $M$ over $SV^*$ with the entries $M^i_j=\rho^i_j\in V^*=(SV^*)_1$. Due to the formulae~\eqref{Deltap}, \eqref{varepsilonp} the conditions~\eqref{rhoij} means exactly that $M$ is a multiplicative matrix. Since $SV^*$ is a commutative algebra this is a usual Manin matrix. Thus $M=(\rho^i_j)$ is the multiplicative first order $A_m$-Manin matrix defining the corepresentation~\eqref{omegarho}. Then, by using Prop.~\ref{PropMorCorep} one can show that the functor~\eqref{Sembedding} is fully faithful.

\begin{Rem} \label{RemInfAlg} \normalfont
 The actions $a\colon\gA\otimes W\to W$ of an infinite-dimensional algebra $\gA=(V,\mu_V,\eta_V)$ on a vector space $W\in\FVect$ can be also translated to the graded homomorphisms $\delta\colon SW^*\to S(V^*\otimes W^*)\to SV^*\circ SW^*$ due to the isomorphism $(V\otimes W)^*\cong V^*\otimes W^*$ natural in $V\in\Vect$ and $W\in\FVect$. However, the multiplication $\mu_V\colon V\otimes V\to V$ does not give a comultiplication $V^*\to V^*\otimes V^*$ in general, since the natural embedding $V^*\otimes W^*\hookrightarrow(V\otimes W)^*$ is not bijective for some infinite-dimensional vector spaces, even if $V=W$. To consider representations of this algebra $\gA$ as quantum representations we need to introduce a topology on $\gA$ and to replace $\Vect$ with an appropriate category of topological $\KK$-modules $W$. This could extend Quantum Representation Theory for infinitely generated quadratic algebras and allow to consider the classical (including infinite-dimensional) representations of an arbitrary algebra at the quantum level.
\end{Rem}

\begin{Rem} \label{RemDeq} \normalfont
 The dequantisation functor $(-)_1^*\colon\QA^\op\to\Vect$ with the strong monoidal structure $(\QA^\op,\circ)\to(\Vect,\otimes)$ translates a corepresentation $\omega\colon\bfcoend(\B)\to\OO$
 of a co\-mo\-noid $\OO=(\A,\Delta,\varepsilon)$ on $\B=\gX_A(\KK)\in\FQA$ to its ``classical limit'' -- the representation $\rho=\omega_1^*\colon\gA\to\bfend(W)$ of the algebra $\gA=(\A_1^*,\Delta_1^*,\varepsilon_1^*)$ on the vector space $W=\B_1^*$.
 The entries of the operators $\rho(v)\in\bfend(W)$ have the form $\rho^i_j(v)=v(M^i_j)$, where $M=(M^i_j)$ is the multiplicative first order $B$-Manin matrix over $\A$ defining the corepresentation~$\omega$. In particular, this functor recovers a classical finite-dimensional representation of a finite-dimensional algebra after its lifting to the quantum level by the $S$-embedding~\eqref{Sembedding}. For a finite-dimensional representation $\rho$ of an infinite-dimensional algebra there can also exist a quantum representation with the classical limit $\rho$, but it is not guaranteed.
\end{Rem}

\bbo{Representations of a semi-linear affine monoid scheme as quantum representations.} \label{bbRepSLM}
In p.~\ref{bbSLAS} we extended the functor $S^*$ to all the semi-linear affine schemes. It has a structure of a lax monoidal functor $S^*\colon(\SLAffSch,\otimes)\to(\FQAsc^\op,\circ)$, so it translates actions and representations of a monoid $\SSS\in\Mon(\SLAffSch,\otimes)$ to actions and representations of the corresponding comonoid $\OO=\Mon(S^*)\SSS\in\Comon(\QAsc,\circ)$.

Recall that any semi-linear affine monoid scheme $\SSS\in\Mon(\SLAffSch,\otimes)$ is given by a bialgebra $(\gR,\alpha,\beta)\in\Comon(\CommAlg,\otimes)$ and a vector space $V\in\FVect$ with linear maps $t\colon V^*\to\gR\otimes\gR\otimes V^*\otimes V^*$, $\varepsilon_{V^*}\colon V^*\to\KK$ such that the diagrams~\eqref{DiagAlphat} commute (see p.~\ref{bbSLAM}). Hence the monoid $\SSS$ can be regarded as the coalgebra $\gC=(\gR\otimes V^*,\Delta_\gC,\varepsilon_\gC)$ with a comultiplication $\Delta_\gC=\sigma^{(23)}\cdot(\alpha,t)\colon\gR\otimes V^*\to\gR\otimes V^*\otimes\gR\otimes V^*$ and a counit $\varepsilon_\gC=\beta\otimes\varepsilon_{V^*}\colon\gR\otimes V^*\to\KK$. The comonoid $\OO=\Mon(S^*)\SSS$ is the semi-connected quadratic algebra $\A=\gR\otimes SV^*$ with morphisms $\Delta\colon\A\to\A\circ\A$ and $\varepsilon\colon\A\to\KK[u]$ whose first order components are $\Delta_\gC$ and $\varepsilon_\gC$.

Let us write more detailed formulae for the case $\SSS\in\Mon(\SLAlgSet,\otimes)$, i.e. for a semi-linear algebraic monoid. The role of the affine scheme $\Spec\gR$ is played by an algebraic set $X$ with the algebra of functions $A(X)=\gR$. It is equipped with a structure of algebraic monoid $(X,\mu_X,\eta_X)$ such that $\alpha=\mu_X^*$ and $\beta=\eta_X^*$; the maps $t$ and $\varepsilon_{V^*}$ are given by $f\colon X\times X\times(V\otimes V)\to V$ and $1_V\in V$ (see~p.~\ref{bbSLAM}).
Let $(w_i)$ and $(w^i)$ be dual bases of vector spaces $W$ and $W^*$.
A representation $\rho\colon X\times V\to\bfend(W)$ and the corresponding action $a\colon X\times(V\otimes W)\to W$ have the form $\rho(x,v)w_j=a(x,v\otimes w_j)=\sum_{i=1}^m \rho^i_j(x,v)w_i$, where $\rho^i_j\in A(X)\otimes V^*$ are semi-linear functions on $X\times V$. 
The commutativity of the diagrams~\eqref{aDiag} is equivalent to $\rho^i_j\big(xy,f_{x,y}(v\otimes v')\big)=\sum_{k=1}^m\rho^i_k(x,v)\rho^k_j(y,v')$ and $\rho^i_j(e,1_V)=\delta^i_j$. This means, in turn, that the matrix $M=(\rho^i_j)$ is multiplicative: $\Delta(\rho^i_j)=\sum_{k=1}^m\rho^i_k\otimes\rho^k_j$, $\varepsilon(\rho^i_j)=\delta^i_j$.
In particular, for $X=\{0\}=\Spec\KK$ we have the case of p.~\ref{bbRepFAlg} with $\rho^i_j(v)=\rho^i_j(0,v)$.

For general $\gR\in\CommAlg$ denote $X=\Spec\gR$. An action $a\colon X\times(V\otimes W)\to W$ is given by a linear map $a^*\colon W^*\to\gR\otimes V^*\otimes W^*$ that should be a coaction of the coalgebra $\gC$ on the space $W^*$. It has the form $a^*(w^i)=\sum_{j=1}^m \rho^i_j\otimes w^j$ for some $\rho^i_j\in\gR\otimes V^*$ such that $M=(\rho^i_j)$ is a multiplicative matrix. The corresponding representation $\rho=\theta^{-1}(a)\colon X\times V\to\bfend(W)$ is given by the coalgebra morphism $\rho^*\colon\bfend(W)^*\to\gC$, $w^i\otimes w_j\mapsto\rho^i_j$.
The lax monoidal functor $S^*\colon(\SLAffSch,\otimes)\to(\QAsc^\op,\circ)$ translates the representation $\rho$ to the corepresentation $\omega\colon\bfcohom(SW^*)\to\OO$. It also has the form~\eqref{omegarho}. In terms of bases it is given by the $A_m$-Manin matrix $M=(\rho^i_j)$. We obtain an extension of the $S$-embedding~\eqref{Sembedding} for the representations of $\SSS$. This is the fully faithful functor
\begin{align} \label{SembeddingS}
 &\Rep_\FVect(\SSS)\hookrightarrow\Corep_\FQA(\OO)^\op, &&(W,\rho)\mapsto(SW^*,\omega),
\end{align}
where $\OO=\Mon(S^*)(\SSS)$.

\bb{Representations of an algebraic monoid or group as quantum representations.} \label{bbRepMon}
Consider an algebraic monoid $\MM=(X,\mu_X,\eta_X)$, this is an algebraic set $X$ equipped with the multiplication $\mu_X\colon X\times X\to X$ and the unity $\eta_X\colon\{0\}\to X$ which give a structure of monoid in $\AlgSet$. We use the notations $xy=\mu_X(x,y)$, $e=\eta_X(0)$. The (not full) embedding $\AlgSet\hookrightarrow\SLAlgSet$ gives a semi-linear algebraic monoid $\SSS_\MM:=(X\times\KK,\mu_X\times\id_\KK,\eta_X\times\id_\KK)$, so we can apply the results of p.~\ref{bbRepSLM}. The function algebra $\gR=A(X)$ is a commutative bialgebra with the comultiplication $(\Delta_\gR f)(x,y)=f(xy)$ and counit $\varepsilon_\gR(f)=f(e)$. Let $\rho\colon\MM\to\bfend(W)$ be a representation of the algebraic monoid $\MM$ on $W\in\FVect$. In the basis $(w_i)$ it has the form $\rho(x)w_j=\sum_i\rho^i_j(x)w_i$, where $\rho^i_j\in A(X)$. The corresponding action $(x,w)\mapsto x.w=\rho(x)w$ satisfies $(xy).w=x.(y.w)$ and $e.w=w$, which is equivalent to the condition of multiplicativity of the matrix $M=(\rho^i_j)$ is multiplicative. This matrix gives the corresponding corepresentation $\omega\colon\bfcoend(SW^*)\to\OO_\gR$, where $\OO_\gR=\Mon(S^*)\SSS_\MM=(\gR\otimes\KK[u],\Delta,\varepsilon)$ is defined by~\eqref{CCgR}.

If $\MM$ is an algebraic group (a group in $\AlgSet$), then $\gR=A(X)$ is a commutative Hopf algebra with the antipode $\zeta_\gR\colon\gR\to\gR$, $(\zeta_\gR f)(x)=f(x^{-1})$, hence the matrix $M=(\rho^i_j)$ is invertible. The inverse matrix $M^{-1}$ has entries $\zeta_\gR(\rho^i_j)$. In this case we have the dual representation $\rho^D\colon\MM\to\bfend(W^*)$. It is defined by the formula $\rho^D(x)\xi=\xi_*\big(\rho(x^{-1})\big)$. In the dual basis $(w^i)$ it reads $\rho^D(x)w^i=\sum_j\rho^i_j(x^{-1})w^j$, so it is given by the matrix $(M^{-1})^\top$. At the quantum level this matrix gives the Koszul dual corepresentation $\omega^{KD}\colon\bfcoend(\Lambda W)\to\OO_\gR$, which we defined in p.~\ref{bbDualQR}, since $(SW^*)^!=\Lambda W$. Note that it does not coincide with the quantum representation corresponding to the classical representation $\rho^D$, because this quantum representation is a corepresentation on $SW$, not on $\Lambda W$, despite they are defined by the same matrix $(M^{-1})^\top$ (see the warning in p.~\ref{bbMultMM}).

Finite monoids and groups are particular cases of algebraic monoids and groups respectively. We can embed a finite monoid $X\in\Mon(\Set,\times)$ into $\AA^1$ as a set of isolated points given by the equation $\prod_{i=1}^N(x-a_i)=0$, where $N$ is the number of elements in $X$ and $a_1,\ldots,a_N$ are arbitrary pair-wise different elements of the infinite field $\KK$. Thus, the description given for general algebraic monoids and groups is valid for the finite case.

\bb{Translation of binary operations with representations by the $S$-embedding.} In the classical representation theory there are two important binary operations with representations: direct sum and tensor product. Consider the question: how the $S$-embedding translates them to the quantum level?

Since the functor $S^*$ has a strong monoidal structure~\eqref{strongSs}, it translates direct sum of classical representations on the vector spaces $W,Z\in\FVect$ to a corepresentation on  the quadratic algebra $S^*(W\oplus Z)=S(W^*\oplus Z^*)=SW^*\otimes SZ^*$. The corresponding quantum representations are corepresentations on $\B=SW^*$ and $\C=SZ^*$, their direct sum defined in p.~\ref{bbDS} is a corepresentation on the same quadratic algebra $\B\otimes\C=SW^*\otimes SZ^*$. Let us describe them in details. Consider a semi-linear monoid $\SSS$ as in p.~\ref{bbRepSLM}. Let $\rho\colon\SSS\to\bfend(W)$ and $\pi\colon\SSS\to\bfend(Z)$ be its representations.
In terms of dual bases $(w_i)$, $(w^i)$ of $W$, $W^*$  and $(z_k)$, $(z^k)$ of $Z$, $Z^*$ these representations are given by the elements $\rho^i_j=\rho^*(w^i\otimes w_j)\in\gR\otimes V^*$ and $\pi^k_l=\pi^*(z^k\otimes z_l)\in\gR\otimes V^*$. The corresponding corepresentations $\omega\colon\bfcoend(SW^*)\to\OO$ and $\nu\colon\bfcoend(SZ^*)\to\OO$ of the comonoid $\OO=\Mon(S^*)\SSS$ are defined by the multiplicative matrices $M=(\rho^i_j)$ and $N=(\pi^k_l)$. The direct sum $\rho\oplus\pi\colon\SSS\to\bfend(W\oplus Z)$ is described via the multiplicative matrix $L=M\oplus N$ as
\begin{align*}
 &(\rho\oplus\pi)^*\colon &
 &w^i\otimes w_j\mapsto\rho^i_j, &
 &w^i\otimes z_l\mapsto0, &
 &z^k\otimes w_j\mapsto0, &
 &z^k\otimes z_l\mapsto\pi^k_l
\end{align*}
(cf. the definition of $L=M\oplus N$ in Prop.~\ref{PropLMN}). Since the entries of $M$ and $N$ belongs to the commutative algebra $\A=\gR\otimes SV^*$, the direct sum of $\omega,\nu$ exists. This is a corepresentation $\omega\dotplus\nu\colon\bfcoend(SW^*\otimes SZ^*)\to\OO$ given by the multiplicative matrix $L=M\oplus N$. On the other hand, the $S$-embedding of the direct sum $\rho\otimes\pi\colon\SSS\to\bfend(W\oplus Z)$ is a corepresentation $\bfcoend(SW^*\otimes SZ^*)\to\OO$ defined by the same matrix $L$, and hence it coincides with the direct sum $\omega\dotplus\nu$. Thus, the $S$-embedding translates a direct sum of classical representations to the direct sum of the corresponding quantum representations.
Since the algebra $\A$ is commutative, we have the symmetric monoidal category $(\Corep_\FQA(\OO),\dotplus)$ and we obtain the following statement.

\begin{Prop}
 The $S$-embedding~\eqref{SembeddingS} has a structure of symmetric strict monoidal functor $\big(\Rep_\FVect(\SSS),\oplus\big)\to\big(\Corep_\FQA(\OO)^\op,\dotplus\big)$.
\end{Prop}

The tensor product of the representations $\rho$ and $\pi$ is defined if $\SSS=(X\times V,\mu_{X\times V},\eta_{X\times V})$ is equipped with a structure of a bimonoid $(X\times V,\mu_{X\times V},\eta_{X\times V},\Delta_{X\times V},\varepsilon_{X\times V})$ in $(\SLAffSch,\otimes)$. The comultiplication $\Delta_{X\times V}$ and the counit $\varepsilon_{X\times V}$ give the monoidal product and unit object in $\Lact(\SSS)$ as it is described in p.~\ref{bbActBimon}. In order to translate the bimonoid structure to the quantum level we need a strong monoidal embedding functor, but unfortunately the monoidal structure~\eqref{laxSs} of the functor $S^*$ is only lax monoidal, so we can not translate the comultiplication $\Delta_{X\times V}$ to the quantum level in general.

Suppose that the monoid $\SSS\in\Mon(\SLAffSch,\otimes)$ is obtained by the category embedding $\AffSch\hookrightarrow\SLAffSch$ from an affine monoid scheme $\MM=(X,\mu_X,\eta_X)\in\Mon(\AffSch,\times)$, where the morphisms $\mu_X$ and $\eta_X$ are given by algebra homomorphisms $\alpha\colon\gR\to\gR\otimes\gR$ and $\beta\colon\gR\to\KK$, so that $\SSS=\SSS_\MM=(X\times\KK,\mu_X\times\id_\KK,\eta_X\times\id_\KK)$. Note that an action $X\times(\KK\otimes W)\to W$ of $\SSS_\MM$ on $W$ is the same as an action $X\times W\to W$ of $\MM$ on $W$, hence we can identify the representations of $\MM$ on $W$ with the representations of $\SSS_\MM$ of $W$. In particular, representations $\rho\colon\MM\to\bfend(W)$ and $\pi\colon\MM\to\bfend(Z)$ are considered as representations of $\SSS=\SSS_\MM$. According to p.~\ref{bbMonCatFP} the monoid $\MM=(X,\mu_X,\eta_X)$ has the unique structure of bimonoid $(X,\mu_X,\eta_X,\Delta_X,\varepsilon_X)$, where the diagonal morphism $\Delta_X\colon X\to X\times X$ and the unique morphism $\varepsilon_X\colon X\to\{0\}$ correspond to the multiplication $\mu_\gR\colon\gR\otimes\gR\to\gR$ and the unity $\eta_\gR\colon\KK\to\gR$. In particular, $\Delta_X$ defines the tensor product $\rho\otimes\pi\colon\MM\to\bfend(W\otimes Z)$. In the case $X\in\AlgSet$ the diagonal morphism has the form $\Delta_X(x)=(x,x)$, so we obtain $(\rho\otimes\pi)(x)=\rho(x)\otimes\pi(x)$. For general $X\in\AffSch$ the representation $\rho\otimes\pi$ is a morphism corresponding to the homomorphism $(\rho\otimes\pi)^*\colon S^*\big(\bfcoend(W\otimes Z)\big)\to\gR$, $\wh w\otimes\wh z\mapsto\rho^*(\wh w)\pi^*(\wh z)$, where $\wh w\in\bfend(W)^*$, $\wh z\in\bfend(Z)^*$, $\wh w\otimes\wh z\in\bfend(W)^*\otimes\bfend(Z)^*=\bfend(W\otimes Z)^*\subset S^*\big(\bfcoend(W\otimes Z)\big)$. In terms of bases one yields $(\rho\otimes\pi)^*(w^i\otimes w_j\otimes z^k\otimes z_l)=\rho^i_j\pi^k_l$.

The $\Bimon$-functor of the embedding $\Alg\hookrightarrow\QAsc$, $\gR\mapsto\gR\otimes\KK[u]$, translates the bimonoid $(\gR,\mu_\gR,\eta_\gR,\Delta_\gR=\alpha,\varepsilon_\gR=\beta)$ to the bimonoid $\BB_\gR=(\gR\circ\KK[u],\mu,\eta,\Delta,\varepsilon)$ defined by the formulae~\eqref{CCgR}, \eqref{BBgR}. The commutative multiplication $\mu=\mu_\gR\otimes\id_{\KK[u]}$ gives the tensor product of corepresentations, which always exists in this case. Thereby, we get the tensor product $\omega\dotcirc\nu\colon\bfcoend(SW^*\circ SZ^*)\to\BB_\gR$. It corresponds to the multiplicative usual Manin matrix $P=M\dototimes N=M^{(1)}N^{(2)}$ with the entries $P^{ik}_{jl}=\rho^i_j\pi^k_l$. On the other hand, $\rho\otimes\pi$ is translated to the corepresentation $\lambda\colon\bfcoend\big(S(W^*\otimes Z^*)\big)\to\BB_\gR$ given by the same matrix $P=M\dototimes N$. Due to Corollary~\ref{CorMorCorep} this means that the graded homomorphism $\phi_{W,Z}\colon S(W^*\otimes Z^*)\to SW^*\circ SZ^*$ is a morphism between the corepresentations $\lambda$ and $\omega\dotcirc\nu$. We can formulate it in the following form.

\begin{Prop}
 In the case $\SSS=\SSS_\MM$ the $S$-embedding~\eqref{SembeddingS} has a structure of symmetric lax monoidal functor $\big(\Rep_\FVect(\SSS_\MM),\otimes\big)\to\big(\Corep_\FQA(\BB_\gR)^\op,\dotcirc\big)$.
\end{Prop}

\begin{Rem} \normalfont
 The $S$-embedding of $\rho\oplus\pi$ differs from the coproduct of $\omega$ and $\nu$ in general, since $S^*(W\oplus Z)\ne S^*(W)\amalg S^*(Z)$.
\end{Rem}

\bbo{Dualisation of a finite-dimensional bialgebra.} \label{bbDuali}
For representations of a finite-dimensional algebra equipped with a bialgebra structure there is another way to translate them to the quantum level such that the representation spaces are translated via the same functor $S^*$.

Let $\gB=(V,\mu_V,\eta_V,\Delta_V,\varepsilon_V)\in\Bimon(\FVect,\otimes)$ be a finite-dimensional bialgebra. This is an algebra $\gA=(V,\mu_V,\eta_V)\in\FAlg$ with a comultiplication $\Delta_\gA=\Delta_V\colon V\to V\otimes V$ and a counit $\varepsilon_\gA=\varepsilon_V\colon V\to\KK$. Application of the contravariant strong monoidal functor $(-)^*\colon(\FVect,\otimes)\to(\FVect,\otimes)$ gives the dual bialgebra
\begin{align} \label{Bdual}
 \gB^*=\Bimon\big((-)^*\big)\gB=(V^*,\mu_{V^*}=\Delta_V^*,\eta_{V^*}=\varepsilon_V^*,\Delta_{V^*}=\mu_V^*,\varepsilon_{V^*}=\eta_V^*)
\end{align}
(see e.g.~\cite[\S~3.2, Ex.~1]{Kass}). Consider the latter one as a comonoid in $(\Alg,\otimes)$. The embedding $\Alg\hookrightarrow\QAsc$ gives the comonoid $\OO_{\gB^*}=(\gB^*\otimes\KK[u],\Delta,\varepsilon)\in\Comon(\FQAsc,\circ)$ defined by~\eqref{CCgR}. It differs from the comonoid $\OO=\Mon(S^*)\gA=(SV^*,\Delta_{SV^*},\varepsilon_{SV^*})$ used in p.~\ref{bbRepFAlg}, but sometimes we can translate a representation of $\gA$ to a corepresentations of $\OO_{\gB^*}$.

Suppose first that the bialgebra $\gB$ is cocommutative. Then the dual bialgebra~\eqref{Bdual} is commutative, so by virtue of Prop.~\ref{PropMonMon} the comonoid $\OO_{\gB^*}$ has the structure of commutative bimonoid $\BB_{\gB^*}=(\gB^*\otimes\KK[u],\mu,\eta,\Delta,\varepsilon)\in\Bimon(\FQAsc,\circ)$ defined by~\eqref{BBgR}. In this case any representation $\rho\colon\gA\to\bfend(W)$ on a vector space $W\in\FVect$ is translated to a corepresentations of $\BB_{\gB^*}$ on $SW^*\in\FQA$. Indeed, due to the commutativity of $\mu_{V^*}$ the linear map $(SV^*)_1=V^*\isoright(\BB_{\gB^*})_1$ induces the graded homomorphism $SV^*\to\gB^*\otimes\KK[u]$. In this way we obtain a comonoid morphism $\OO\to\BB_{\gB^*}$. By composing it with~\eqref{omegarho} we get the corepresentation $\omega'$ of the bimonoid $\BB_{\gB^*}$ on the quadratic algebra $SW^*$:
\begin{align}
 \omega'\colon\bfcoend(SW^*)\xrightarrow{\Phi_{W,W}}S\big(\bfend(W)\big)^*\xrightarrow{S\rho^*}\OO\to\BB_{\gB^*}. \label{omegarhoBB}
\end{align}
Let $(w_i)_{i=1}^m$ be a basis of $W$ and $\rho^i_j\in V^*$ be such that $\rho(v)w_j=\sum_i\rho^i_j(v)w_i$, then $\rho$ is translated to the corepresentation $\omega'\colon\bfcoend(SW^*)\to\BB_{\gB^*}$ given by the usual Manin matrix $M=(\rho^i_j)$. In this way we obtain the fully faithful functor
\begin{align} \label{RepBB}
 &\Rep_\FVect(\gA)\hookrightarrow\Corep_\FQA(\BB_{\gB^*})^\op, &&(W,\rho)\mapsto(SW^*,\omega').
\end{align}

Let $\pi\colon\gA\to\bfend(Z)$ be a representation of $\gA$ on a vector space $Z\in\FVect$. The functor~\eqref{RepBB} translates $\pi$ to a corepresentation $\nu'\colon\bfcoend(SZ^*)\to\BB_{\gB^*}$. The direct sum $\rho\oplus\pi\colon\gA\to\bfend(W\oplus Z)$ is translated to $\omega'\dotplus\nu'\colon\bfcoend(SW^*\otimes SZ^*)\to\BB_{\gB^*}$, so we can regard~\eqref{RepBB} as a strict monoidal functor $\big(\Rep_\FVect(\gA),\oplus\big)\to\big(\Corep_\FVect(\BB_{\gB^*})^\op,\dotplus\big)$.
The translation of the tensor product by the functor~\eqref{RepBB} gives a lax monoidal functor $\big(\Rep_\FVect(\gA),\otimes\big)\to\big(\Corep_\FVect(\BB_{\gB^*})^\op,\dotcirc\big)$. If $\gB$ is a Hopf algebra with an antipode $\zeta_V\colon V\to V$, then~\eqref{Bdual} is also a Hopf algebra, its antipode is $\zeta_V^*\colon V^*\to V^*$ (see~\cite[\S~3.3, Prop.~3.3.3]{Kass}), so $\BB_{\gB^*}$ is a Hopf monoid and, hence, one can define the Koszul dual corepresentation $(\omega')^{KD}\colon\bfcoend(\Lambda W)\to\BB_{\gB^*}$.

Let $\MM=(X,\mu_X,\eta_X)\in\Mon(\Set,\times)$ be a finite monoid. The algebra $\gB=\KK[\MM]$ defined in p.~\ref{bbRepAlg} is a cocommutative bialgebra with the coalgebra structure
\begin{align}
 &\Delta_V\Big(\sum_{x\in X}\alpha_x x\Big)=\sum_{x\in X}\alpha_x (x\otimes x), &
 &\varepsilon_V\Big(\sum_{x\in X}\alpha_x x\Big)=\sum_{x\in X}\alpha_x.
\end{align}
It is a Hopf algebra iff $\MM$ is a group. The dual bialgebra $\gB^*$ coincides with the bialgebra $\gR=A(X)$ constructed in p.~\ref{bbRepMon}. A representation $\MM\to\bfend(W)$ is uniquely extended to a representation $\KK[\MM]\to\bfend(W)$. If we translate these representations to the corepresentations $\omega\colon\bfcohom(SW^*)\to\OO_{\gR}$ and $\omega'\colon\bfcohom(SW^*)\to\OO_{\gB^*}$ respectively, then the results coincide: $\omega=\omega'$, so in this case the functor~\eqref{RepBB} is just another description of the representation $S$-embedding obtained in p.~\ref{bbRepMon}.

For a non-cocommutative bialgebra $\gB\in\Bimon(\FVect,\otimes)$ we can not translate an arbitrary representation $\rho\colon\gA\to\bfend(W)$ to a corepresentation of $\OO_{\gB^*}$, but if the corresponding matrix $M=(\rho^i_j)$ is an $A_m$-Manin matrix as a matrix over $\gB^*$, then we have a corepresentation $\omega'\colon\bfcoend(SW^*)\to\OO_{\gB^*}$ given by this matrix, since this matrix is multiplicative with respect to $\Delta_{V^*}$ and $\varepsilon_{V^*}$ anyway.

\subsection{$T$-embedding of classical representations}

An alternative way to translate a classical representation to the quantum level is to use the functor $T^*\colon\FVect\to\FQA^\op$ defined in p.~\ref{bbT} instead of $S^*$.

\bb{$T$-embedding of a finite-dimensional algebra and its representations.}
Let $\gA=(V,\mu_V,\eta_V)\in\FAlg$. The strong monoidal functor $T^*\colon(\FVect,\otimes)\to(\FQA^\op,\circ)$ induces the functor $\Mon(T^*)$, which maps the algebra $\gA\in\Mon(\FVect,\otimes)$ to the comonoid
 $\OO=\Mon(T^*)\gA=\Comon(T^\op)(V^*,\mu_V^*,\eta_V^*)=(TV^*,\Delta,\varepsilon)$ in the category $(\FQA,\circ),$
 where  $\Delta\colon TV^*\xrightarrow{T\mu_V^*}
T^*(V\otimes V)\cong T(V^*\otimes V^*)\cong TV^*\circ TV^*$, $\varepsilon=T\eta_V^*\colon TV^*\to T\KK=\KK[u]$.

Due to Prop.~\ref{PropQA} we obtain $\bfcohom(TW^*,TZ^*)=T\big(\bfhom(W^*,Z^*)\big)\cong T\big(\bfhom(W,Z)\big)^*$, so the morphisms~\eqref{PhiDef} for the monoidal functor $F=T^*\colon(\FVect,\otimes)\to(\FQA^\op,\circ)$ is the isomorphisms $\Phi_{W,Z}\colon\bfcohom(TW^*,TZ^*)\isoright T\big(\bfhom(W,Z)\big)^*$ given by the same identification $\bfhom(W^*,Z^*)\cong\bfhom(W,Z)^*$ as for $S^*$.

Let $\rho\colon\gA\to\bfend(W)$ be a representation of the algebra $\gA$ on $W\in\FVect$. In a basis $(w_i)_{i=1}^m$ of $W$ it is given by the linear functions $\rho^i_j\in V^*$ satisfying~\eqref{rhoij}. By using the isomorphism $\Phi_{W,W}$ we obtain the corresponding quantum representation
\begin{align}
 \omega\colon\bfcoend(TW^*)\cong T\big(\bfend(W)\big)^*\xrightarrow{T\rho^*}\OO. \label{omegarhoT}
\end{align}
This is a corepresentation of $\OO=(TV^*,\Delta,\varepsilon)$ on the quadratic algebra $TW^*=\gX_{0_m}(\KK)$, where $0_m\in\bfend(\KK^m\otimes\KK^m)$ is the zero idempotent and the identification $W=\KK^m$ is fixed by the basis $(w_i)$. Recall that any $m\times m$ matrix is a $0_m$-Manin matrix. The corepresentation~\eqref{omegarhoT} is defined by the multiplicative first order matrix $M=(\rho^i_j)$.
Thus, the {\it $T$-embedding} of the representations is the fully faithful functor
\begin{align} \label{Tembedding}
 &\Rep_\FVect(\gA)\hookrightarrow\Corep_\FQA(\OO)^\op, &&(W,\rho)\mapsto(TW^*,\omega).
\end{align}

\bbo{$T$-embedding of a semi-linear affine monoid scheme and its representations.}
The functor $T^*$ can be extended to the functor $T^*\colon\SLAffSch\to\QAsc^\op$ by the formula $T^*\colon(\Spec\gR)\times V\mapsto\gR\otimes TV^*$ with the obvious mapping of the morphisms. This is a strong monoidal functor $(\SLAffSch,\otimes)\to(\QAsc^\op,\circ)$.

Let $\SSS\in\Mon(\SLAffSch,\otimes)$ be a semi-linear affine monoid scheme as in p.~\ref{bbRepSLM}. Then the comonoid $\OO=\Mon(T^*)\SSS$ is the semi-connected quadratic algebra $\A=\gR\otimes T^*V$ with the graded homomorphisms $\Delta\colon\A\to\A\circ\A$ and $\varepsilon\colon\A\to\KK[u]$ such that their first order components $\Delta_1$ and $\varepsilon_1$ coincide with $\Delta_\gC\colon\gR\otimes V^*\to\gR\otimes V^*\otimes\gR\otimes V^*$ and $\varepsilon_\gC\colon\gR\otimes V^*\to\KK$, which were described in p.~\ref{bbRepSLM}.

Let $\rho\colon\SSS\to\bfend(W)$ be a representation of $\SSS$ on a space $W\in\FVect$ with a basis $(w_i)_{i=1}^m$. It is given by the coalgebra morphism $\rho^*\colon\bfend(W)^*\to\gC$, $w^i\otimes w_j\to\rho^i_j$, where $\gC=(\gR\otimes V^*,\Delta_\gC,\varepsilon_\gC)$ and $M=(\rho^i_j)$ is a multiplicative first order matrix over the algebra $\gR\otimes TV^*$ with respect to $\Delta$ and $\varepsilon$. Since the functor $T^*\colon(\SLAffSch,\otimes)\to(\QAsc^\op,\circ)$ is strong monoidal, it translates $\rho$ to the corepresentation $\omega\colon\bfcoend(TW^*)\to\OO$ corresponding to the $0_m$-Manin matrix $M=(\rho^i_j)$.
This extends $T$-embedding functor~\eqref{Tembedding} to the fully faithful functor
\begin{align} \label{TembeddingS}
 \Rep_\FVect(\SSS)\hookrightarrow\Corep_\FQA(\OO)^\op, &&(W,\rho)\mapsto(TW^*,\omega),
\end{align}
where $\OO=\Mon(T^*)\SSS$.

\bb{Translation of binary operations under the $T$-embedding.}
Let representations $\rho\colon\SSS\to\bfend(W)$, $\pi\colon\SSS\to\bfend(Z)$ be given by the elements $\rho^i_j,\pi^k_l\in\gR\otimes V^*$ in bases $(w_i)$ and $(z_k)$. Denote the corresponding corepresentations $\bfcoend(TW^*)\to\OO$ and $\bfcoend(TZ^*)\to\OO$ given by the matrices $M=(\rho^i_j)$ and $N=(\pi^k_l)$ by $\omega$ and $\nu$.
The functor $T^*\colon\FVect\to\FQA^\op$ has a structure of strong monoidal functor~\eqref{TstFVectplus}. By using Corollary~\ref{CorMorCorep} again we see that the $T$-embedding translates the direct sum representation $\rho\oplus\pi\colon\SSS\to\bfend(W\oplus Z)$ to the coproduct $\omega\dotsqcup\nu\colon\bfcoend(TW^*\amalg TZ^*)\to\OO$ given by the matrix $L=M\oplus N$.
This means that $T$-embedding~\eqref{TembeddingS} preserves the finite products. In other words, the functor~\eqref{TembeddingS} has a structure of symmetric strong monoidal functor $\big(\Rep_\FVect(\SSS),\oplus\big) \to \big(\Corep_\FQA(\OO)^\op,\dotsqcup\big)$.
 
By taking into account the lax monoidal structure~\eqref{laxTstFVectplus} of the category embedding $T^*\colon\FVect\hookrightarrow\FQA^\op$ we see that the $T$-embedding gives a symmetric lax monoidal functor $\big(\Rep_\FVect(\SSS),\oplus\big) \to \big(\Corep_\FQA(\OO)^\op,\dotplus\big)$.

Since $T^*\colon\FVect\hookrightarrow\FQA^\op$ has a structure of symmetric strong monoidal functor~\eqref{TstFVectotimes}, we can apply the results of p.~\ref{bbTrTP}. In this way we obtain a strong monoidal functor $\big(\Rep_\FVect(\SSS),\otimes\big) \to \big(\Corep_\FQA(\OO)^\op,\dotcirc\big)$ given by the $T$-embedding~\eqref{TembeddingS}.

\bb{Representations of a finite-dimensional bialgebra.}
Let $\gB=(V,\mu_V,\eta_V,\Delta_V,\varepsilon_V)$ be a bialgebra as in p.~\ref{bbDuali}. Then we can translate any representation of the algebra $\gA=(V,\mu_V,\eta_V)$ on $W\in\FVect$ to a quantum representation of the comonoid $\OO_{\gB^*}$ on $TW^*$. We obtain the fully faithful functor
\begin{align} \label{TEmbBiAlg}
 &\Rep_\FVect(\gA)\hookrightarrow\Corep_{\FQA}(\BB_{\gB^*})^\op, &&(W,\rho) \mapsto(TW^*,\omega'),
\end{align}
where $\rho\colon\gA\to\bfend(W)$, $\rho(v)w_j=\sum_i\rho^i_j(v)w_i$, $\rho^i_j\in V^*$ and $\omega'\colon\bfcoend(TW^*)\to\OO_{\gB^*}$ is the corepresentation given by the multiplicative matrix $M=(\rho^i_j)$. Note that in this case we do not need to require the condition of cocommutativity of $\gB$ nor any Manin condition on $M$.

Let $\pi\colon\gA\to\bfend(Z)$ be representation and $\nu'\colon\bfcoend(TZ^*)\to\OO_{\gB^*}$ be the corresponding corepresentation. The embedding~\eqref{TEmbBiAlg} translates direct sum $\rho\oplus\pi\colon\gA\to\bfend(W\oplus Z)$ to the coproduct $\omega'\dotsqcup\nu'\colon\bfcoend(TW^*\amalg TZ^*)\to\OO_{\gB^*}$. It translates the tensor product of representations $\rho\otimes\pi\colon\gA\to\bfend(W\otimes Z)$ defined by means of $\Delta_V$ to the tensor product $\omega'\dotcirc\nu'\colon\bfcoend(TW^*\circ TZ^*)\to\OO_{\gB^*}$.

\section{Examples}
\label{sec6}

\subsection{Corepresentations of $M_q(m)$}

The relationship of some Manin matrices with the Lax operators of $U_q(\mathfrak{gl}_n)$ type was described in details in~\cite{S}. It can be interpreted in terms of Quantum Representation Theory.

\bb{Matrix algebra.} \label{bbMatAlg}
Consider the algebra $\gA=\Mat_m(\KK)$. It consists of $m\times m$ matrices over $\KK$. The multiplication in $\gA$ is the usual matrix multiplication and the unity is the identity matrix. The quadratic algebra $M(m)=M(m,\KK)=S\gA^*$ is the algebra of regular functions on $\Mat_m(\KK)$ (as on an algebraic set). It is generated by the functions $a^i_j\colon M\mapsto M^i_j$, $M\in\Mat_m(\KK)$. Thus $M(m)$ is the graded algebra generated by $m^2$ order 1 elements $a^i_j\in\gA^*$. The multiplication and unity of $\Mat_m(\KK)$ gives the structure of comonoid on $M(m)$. We obtain the comonoid $\OO=M(m)=(M(m),\Delta,\varepsilon)\in\Comon(\FQA,\circ)$, where the graded homomorphisms $\Delta\colon M(m)\to M(m)\circ M(m)$, $\varepsilon\colon M(m)\to\KK$ have the form
\begin{align} \label{Deltaa}
 &\Delta(a^i_j)=\sum_{k=1}^m a^i_k\otimes a^k_j, &&\varepsilon(a^i_j)=\delta^i_j.
\end{align}
The algebra $M(m)$ with this structure can be considered as a bialgebra, but this is not a Hopf algebra.

The matrices from $\Mat_m(\KK)$ act on elements of $W=\KK^m$ as on column vectors. The corresponding representation $\rho\colon\Mat_m(\KK)\to\bfend(\KK^m)$ is the isomorphism identifying operators on $\KK^m$ with the matrices. By the $S$-embedding described in p.~\ref{bbRepFAlg} we get the corepresentation $\omega\colon\bfcoend\big(\gX_{A_m}(\KK)\big)\to M(m)$ given by the matrix $M=(a^i_j)$.

The right action of $\Mat_m(\KK)$ on the space $\KK^m$ considered as the space of row-vectors gives the representation $\rho^\top\colon\Mat_m(\KK)^\op\to\bfend(\KK^m)$. The corresponding corepresentation $\omega^\top\colon\bfcoend\big(\gX_{A_m}(\KK)\big)\to M(m)^\cop$ is defined by the transposed matrix $M^\top$. Note that the commutation relations $a^i_ja^k_l=a^k_la^i_j$, which define the quadratic algebra $M(m)$, are equivalent to the requirement that both $M$ and $M^\top$ are $A_m$-Manin matrices (see e.g.~\cite[Prop.~3.1]{S}).

\bb{A $q$-deformation of $\KK^m$.}
Consider the $\KK$-plane $\KK^2$. At the quantum level it is described by the quadratic algebra $\gX_{A_2}(\KK)=\KK[x,y]$. Let $q\in\KK\backslash\{0\}$. Replacement of the commutation relation $yx=xy$ by $yx=qxy$ gives a $q$-deformation of $\KK[x,y]$ called {\it quantum plane}. It can be interpreted as a $q$-deformation of the vector space $\KK^2$ in the category of quantum linear spaces (where the vector spaces considered as quadratic algebras via the embedding $S^*\colon\FVect\hookrightarrow\FQA^\op$).

For a general dimension $m$ a quantum $q$-deformation of $\KK^m$ is constructed as the quadratic algebra generated by $x^1,\ldots,x^m$ with the commutation relations $x^jx^i=qx^ix^j$, $i<j$. This is the algebra $\gX_{A^q_m}(\KK)$ for the idempotent $A_m^q\in\bfend(\KK^m\otimes\KK^m)$ defined by the formulae%
\footnote{More generally, one can construct a multi-parametric deformation~\cite{Manin89} by imposing the conditions $q_{ij}=q_{ji}^{-1}$, $q_{ii}=1$ only (see~\cite[\S~3.3]{S}).
}
\begin{align}
 &A^q_m=\frac12(1-P^q_m), && (P^q_m)^{ij}_{kl}=q_{ji}\delta^i_l\delta^j_k,
\end{align}
$q_{ij}=q_{ji}^{-1}=q$ for $i<j$, $q_{ii}=1$.

\bb{A $q$-deformation of the matrix algebra.}
Let $M_q(m)=M_q(m,\KK)$ be the bialgebra generated by $a^i_j$, $i,j=1,\ldots,m$, with the commutation relations
\begin{gather}
\begin{split} \label{aaq12}
 &a^j_ka^i_k=qa^i_ka^j_k, \qquad\qquad a^i_la^j_k=a^j_ka^i_l, \\
 &a^i_la^i_k=qa^i_ka^i_l, \qquad\qquad
 a^i_ka^j_l-a^j_la^i_k=(q^{-1}-q)a^j_ka^i_l,
\end{split}
\end{gather}
where $i<j$, $k<l$, the comultiplication and counit have the same form~\eqref{Deltaa} (see e.g. \cite[\S~4.10]{Kass}, \cite{FRT89}). Since these commutation relations are quadratic and the maps~\eqref{Deltaa} satisfy~\eqref{DeltaAk} (for $k=1$) the algebra $M_q(m)$ can be considered as a comonoid in $(\FQA,\circ)$.

In~\cite{Manin87,Manin88} Manin interpreted $M_q(2)$ in terms of the quantum plane and gave a description of related quantum groups via the quantum linear spaces. In the case of the bialgebra $M_q(m)$ this approach can be formulated as the following fact, which was discussed in the works~\cite[\S~3.3, Prop.~2]{CFR}, \cite[\S~2.4, Prop.~2.4]{qManin}, \cite[\S~4.1, Th.~4.4]{S}. 

\begin{Prop}\label{PropMq}
 The relations~\eqref{aaq12} are equivalent to the requirement that both $M=(a^i_j)$ and its transposed $M^\top$ are $A^q_m$-Manin matrices.
\end{Prop}

It follows that the $A^q_m$-Manin matrices $M=(a^i_j)$ and $M^\top$ over $M_q(m)$ define the corepresentations
\begin{align} \label{omegaq}
 &\omega_q\colon\bfcoend\big(\gX_{A^q_m}(\KK)\big)\to M_q(m), &
 &\omega_q^\top\colon\bfcoend\big(\gX_{A^q_m}(\KK)\big)\to M_q(m)^\cop
\end{align}
respectively.
These corepresentations are $q$-deformations of the standard matrix representations $\rho$ and $\rho^\top$ defined in p.~\ref{bbMatAlg}. Namely, we apply the embedding $S^*$ to the matrix algebra $\Mat_m(\KK)$, to the representation space $\KK^m$, to the representations $\rho$, $\rho^\top$ and then we simultaneously deform them all.

If $m=2$ we obtain the quantum representation of $M_q(2)$ on the quantum plane. This case is described in details in~\cite{Manin87,Manin88,Kass} in terms of bialgebras and coactions in the monoidal category $(\Alg,\otimes)$.

The comonoid $M_q(m)$ and its quantum representations~\eqref{omegaq} can be interpreted in the following way. Consider the graded homomorphisms $f\colon T\big(\Mat_m(\KK)\big)^*\to\bfcoend\big(\gX_{A^q_m}(\KK)\big)$ and $f^\top\colon T\big(\Mat_m(\KK)\big)^*\to\bfcoend\big(\gX_{A^q_m}(\KK)\big)^\cop$ defined as $f(a^i_j)=\M^i_j$ and $f^\top(a^i_j)=\M^j_i$, where $a^i_j\in\big(\Mat_m(\KK)\big)^*$ are linear functions on the vector space $\Mat_m(\KK)$ (see p.~\ref{bbMatAlg}) and $\M^i_j$ are entries of the universal $A^q_m$-Manin matrix. These homomorphisms are morphisms in the category $\Comon(\FQA,\circ)$, since they preserve the comultiplications and counits.
Prop.~\ref{PropMq} means exactly that the graded algebra $M_q(m)=M_q(m,\KK)$ is the pushout
\begin{align}
 \xymatrix{
T\big(\Mat_m(\KK)\big)^*\ar[r]^{f^\top}\ar[d]^f & \bfcoend\big(\gX_{A^q_m}(\KK)\big)^\cop\ar[d]^{\omega_q^\top} \\
\bfcoend\big(\gX_{A^q_m}(\KK)\big)\ar[r]^{\omega_q} & M_q(m)
} 
\end{align}
in the category $\GrAlg$ and, consequently, in $\FQA$ (see~\cite[\S~3.3]{Mcl} for a definition of pushout). Moreover, the comonoid $M_q(m)$ is a pushout of $f$ and $f^\top$ in the category $\Comon(\FQA,\circ)$.

\subsection{An extension of the Yangian $Y(\mathfrak{gl}_m)$ and its corepresentations}
\label{secYangA}

In~\cite{CF} the authors presented an example of a usual Manin matrix defined as a product of the Lax operator for the $\mathfrak{gl}_m$ Yangian with a shift operator. We interpret this extension as a corepresentation of some extension of this Yangian on the quadratic algebra $\gX_{A_m}(\KK)$.

\bb{Manin matrix from the $Y(\mathfrak{gl}_m)$ Lax operator.} \label{bbYgln}
Denote by $P_m$ the operator acting by permutation of the tensor factors in $\KK^m\otimes\KK^m$, that is $P_m=\sigma_{\KK^m,\KK^m}\in\bfend(\KK^m\otimes\KK^m)$. It has entries $(P_m)^{ij}_{kl}=\delta^i_l\delta^j_k$ and is related with the anti-symmetrizer as $A_m=\dfrac{1-P_m}2$.

Consider the rational $R$-matrix $R(z)=z-P_m$. The Yangian $Y(\mathfrak{gl}_m)$ is the algebra generated (over $\KK$) by $t_{ij}^r$, $i,j=1,\ldots,m$, $r\in\ZZ_{\ge1}$, with the relations
\begin{align}
 R(z-v)T^{(1)}(z)T^{(2)}(w)=T^{(2)}(w)T^{(1)}(z)R(z-w), \label{RTT}
\end{align}
where $T(z)$ is the $m\times m$ matrix over $Y(\mathfrak{gl}_m)[[z^{-1}]]$ with the entries $T(z)^i_j=\delta^i_j+\sum_{r\ge1}t_{ij}^rz^{-r}$. The Yangian is equipped with a structure of Hopf algebra:
\begin{align} \label{DeltaY}
 &\Delta\big(T(z)^i_j\big)=\sum_{k=1}^m T(z)^i_k\otimes T(z)^k_j, &
 &\varepsilon\big(T(z)^i_j\big)=\delta^i_j, &
 &\zeta\big(T(z)\big)=T(z)^{-1}.
\end{align}

Consider the shift operator $e^{-\frac{\partial}{\partial z}}$. The relation~\eqref{RTT} implies that the matrix $T(z)e^{-\frac{\partial}{\partial z}}$ is an $A_m$-Manin matrix over the algebra $Y(\mathfrak{gl}_m)[[z^{-1}]][e^{-\frac{\partial}{\partial z}}]$ (see details in~\cite{CF}, \cite{S}). This fact is equivalent to the relation
\begin{align} \label{ATTtau}
 A_mT(z)T(z-1)(1-A_m)=0.
\end{align}

\bb{Extended Yangian and its corepresentation.}
Let us add one more generator $\tau$ to the Yangian and postulate the commutation relation
\begin{align} \label{Ttau}
 \tau T(z)=T(z-1)\tau,
\end{align}
see~\cite{Kh}; via generators it can be written in the form $\tau t_{ij}^r=\sum_{k=1}^r{r-1\choose k-1}t_{ij}^k\tau$. We obtain a bialgebra, where the comultiplication and counit of the new generator are defined by the formulae $\Delta(\tau)=\tau\otimes\tau$, $\varepsilon(\tau)=1$. This is an extension of $Y(\mathfrak{gl}_m)$ as a bialgebra, denote it by $Y(\mathfrak{gl}_m)[\tau]$. It has the form of the tensor product $Y(\mathfrak{gl}_m)\otimes\KK[\tau]$ as a vector space, but not as an algebra, since $\tau$ does not commute with the generators $t_{ij}^r$. 

 The relations~\eqref{ATTtau}, \eqref{Ttau} imply that the matrix $M=T(z)\tau$ is an $A_m$-Manin matrix over the algebra $\gR=Y(\mathfrak{gl}_m)[\tau]((z^{-1}))$, where we extend the basic field to the field of formal Laurent series
\begin{align} \label{KKz}
 \KK((z^{-1}))=\Big\{\sum_{k=-\infty}^{N}\alpha_k z^k\mid N\in\ZZ,\alpha_k\in\KK\Big\}.
\end{align}
The $A_m$-Manin matrix $M=T(z)\tau$ is multiplicative with respect to the comultiplication $\Delta_\gR\colon\gR\to\gR\otimes_{\KK((z^{-1}))}\gR$ and counit $\varepsilon_\gR\colon\gR\to\KK((z^{-1}))$ obtained by the field extension, so it gives a corepresentation
\begin{align}
 \omega\colon\bfend\Big(\gX_{A_m}\big(\KK((z^{-1}))\big)\Big)\to\OO_\gR
\end{align}
of the comonoid $\OO_\gR=\big(\gR\circ_{\KK((z^{-1}))}\KK((z^{-1}))[u],\Delta_\gR\otimes_{\KK((z^{-1}))}\id_{\KK((z^{-1}))[u]},\varepsilon_\gR\otimes_{\KK((z^{-1}))}\id_{\KK((z^{-1}))[u]}\big)$, where the operation $\circ_{\KK((z^{-1}))}$ is defined over the extended field~\eqref{KKz}, we have the bialgebra isomorphism $\gR\circ_{\KK((z^{-1}))}\KK((z^{-1})))[u]=\gR\circ_{\KK}\KK[u]$ over $\KK$. The quantum representation $\omega$ has the classical representation space $\KK^m$ (with respect to the $S$-embedding), but the quantum algebra we represent is not classical.

In this construction we consider the elements of $\gR=Y(\mathfrak{gl}_m)[\tau]((z^{-1}))$ as elements of order $0$, the graded algebra is $\gR[u]$. Alternatively one can set $\deg\tau=1$. Then the bialgebra $Y(\mathfrak{gl}_m)[\tau]((z^{-1}))$ itself becomes a comonoid in $(\GrAlg,\circ)$, and entries of $M=T(z)\tau$ become elements of order $1$. This is an affinely generated quadratic algebra over $Y(\mathfrak{gl}_m)((z^{-1}))$, but it is not semi-connected.

\begin{Rem} \normalfont
 The localisation of $Y(\mathfrak{gl}_m)[\tau]$ by $\tau^{-1}$ gives us a Hopf algebra $Y(\mathfrak{gl}_m)[\tau^{\pm1}]$ with the antipode extended as $\zeta(\tau)=\tau^{-1}$. However, the formula $\deg\tau=1$ and its consequence $\deg\tau^{-1}=-1$ define a $\ZZ$-grading rather than an $\NN_0$-grading on the algebra $Y(\mathfrak{gl}_m)[\tau^{\pm1}]$.
\end{Rem}

\subsection{Corepresentations of an extended $Y(\mathfrak{so}_m)$}

As we see in~\cite[\S~7.1]{S}, the Lax operator of the Yangians $Y(\mathfrak{so}_m)$ and $Y(\mathfrak{sp}_m)$ also give Manin matrices for some idempotents, so we can construct quantum representations for them by analogy with Subsection~\ref{secYangA}. For simplicity we consider the orthogonal case only. The symplectic version is completely analogues.

Here we suppose that $\chara\KK=0$.

\bb{The algebras $X(\mathfrak{so}_m)$ and $Y(\mathfrak{so}_m)$.}
Consider again the collection of generators $t_{ij}^r$, where $i,j=1,\ldots,m$, $r\in\ZZ_{\ge1}$. Let $T(u)$ be the $m\times m$ matrix defined by the same formula as in p.~\ref{bbYgln}. Denote by $X(\mathfrak{so}_m)$ the algebra with these generators and commutation relations~\eqref{RTT} but for another $R$-matrix $R(z)=R_{\mathfrak{so}_m}(z)=1-\dfrac{P_m}{z}+\dfrac{Q_m}{z-m/2+1}$, where $Q_m$ is the operator on $\KK^m\otimes\KK^m$ with the entries $(Q_m)^{ij}_{kl}=\delta^{i+j}_{m+1}\delta^{m+1}_{k+l}$. The Yangian $Y(\mathfrak{so}_m)$ can be defined as a quotient of $X(\mathfrak{so}_m)$ by the relations $\sum_{k=1}^mT(z)^i_kT(z+\frac{m}{2}-1)^j_{m+1-k}=\delta^{i+j}_{m+1}$. Both $X(\mathfrak{so}_m)$ and $Y(\mathfrak{so}_m)$ are Hopf algebras with the comultiplication, counit and antipode of the same form~\eqref{DeltaY} (see details in~\cite[\S~2, \S~3, Cor.~3.2]{AMR}, \cite[\S~11.1]{MolevSO} and \cite[\S~7.1]{S}).
\footnote{The algebra $X(\mathfrak{so}_m)$ is usually called `extended Yangian', but we use this term for another algebra. 
}

\bb{An algebra $\gX_{B_m}(\KK)$ and a $B_m$-Manin matrix.}
Consider the idempotent
\begin{align} \label{Bn}
 &B_m=\frac{1-P_m}2+\frac{Q_m}m=A_m+\frac{Q_m}m\in\bfend(\KK^m\otimes\KK^m)
\end{align}
(introduced in~\cite[\S~7.1]{S}). It gives the commutative quadratic algebra $\gX_{B_m}(\KK)$ isomorphic to the quotient of $\gX_{A_m}(\KK)=\KK[x^1,\ldots,x^m]$ by the quadratic relations
\begin{align}
 \sum_{i=1}^m x^ix^{m+1-i}=0.
\end{align}

It is proved in~\cite[\S~7.1]{S} that the matrix $T(z)e^{-\frac{\partial}{\partial z}}$ is a $B_m$-Manin matrix over the algebra $Y(\mathfrak{so}_m)[[z^{-1}]][e^{-\frac{\partial}{\partial z}}]$. This is equivalent to the relation
\begin{align} \label{BTTtau}
 B_mT(z)T(z-1)(1-B_m)=0.
\end{align}

\bbo{A corepresentation on $\gX_{B_m}(\KK)$.}
Let us extend the Yangian $Y(\mathfrak{so}_m)$ by a new generator $\tau$ and impose the relations~\eqref{Ttau}. This is a bialgebra over $\KK$ (the comultiplication and counit has the same form as for the $\mathfrak{gl}_m$ case). It is isomorphic to $Y(\mathfrak{so}_m)\otimes\KK[\tau]$ as a vector space, so we analogously denote it by $Y(\mathfrak{so}_m)[\tau]$. The field extension gives the bialgebra $\gR=Y(\mathfrak{so}_m)[\tau]((z^{-1}))$ over~\eqref{KKz}.

Due to~\eqref{BTTtau} the matrix $M=T(z)\tau$ is a $B_m$-Manin matrix over $\gR$. It is multiplicative with respect to the $\KK((z^{-1}))$-linear comultiplication and counit of the bialgebra $\gR$, so it defines a corepresentation
\begin{align}
 \omega\colon\bfend\Big(\gX_{B_m}\big(\KK((z^{-1}))\big)\Big)\to\OO_\gR
\end{align}
of the comonoid $\OO_\gR=(\gR\circ\KK((z^{-1})))[u],\Delta_\gR\otimes_{\KK((z^{-1}))}\id_{\KK((z^{-1}))[u]},\varepsilon_\gR\otimes_{\KK((z^{-1}))}\id_{\KK((z^{-1}))[u]})$.

\begin{Rem} \normalfont
 Since the relation~\eqref{BTTtau} follows from~\eqref{RTT} only the matrix $T(z)\tau$ can be considered as a $B_m$-Manin matrix over the bialgebra $X(\mathfrak{so}_m)[\tau]((z^{-1}))$. By using the localisation we can regard it as a $B_m$-Manin matrix over the Hopf algebra $Y(\mathfrak{so}_m)[\tau^{\pm1}]((z^{-1}))$ or $X(\mathfrak{so}_m)[\tau^{\pm1}]((z^{-1}))$.
\end{Rem}

\section*{Conclusion and further directions}
\addcontentsline{toc}{section}{Conclusion and further directions}

\noindent{\bf Quantum linear spaces and quantum representations.}
 The Manin product `$\circ$' allows us to construct Quantum Representation Theory on quantum linear spaces $\A\in\FQA^\op$. Namely, we applied the general approach described in Section~\ref{sec3} to the symmetric monoidal category $\bfC=(\QAsc^\op,\circ)$ and its subcategory $\bfP=\FQA^\op$. By using the category embedding $S^*\colon\FVect\hookrightarrow\FQA^\op$ we can consider the finite-dimensional vector spaces and as quantum linear spaces. Its extension $S^*\colon\SLAffSch\hookrightarrow\QAsc^\op\hookrightarrow\GrAlg^\op$ helps us to interpret representations of finite-dimensional algebras and algebraic groups as quantum representations. The binary operations with vector spaces and representations agree with their quantum versions (modulo some natural transformations $\phi$).

\vspace*{2mm}\noindent{\bf Semi-linear spaces.}
 Introduction of semi-linear spaces gives the unified theory of finite-dimensional representations of finite-dimensional algebras and algebraic groups: one needs to apply the approach of Section~\ref{sec3} to $\bfC=(\SLAffSch,\otimes)$ and $\bfP=\FVect$. We have the same situation at the quantum level, where the role of semi-linear spaces is played by semi-connected quadratic algebras: $\bfC=(\QAsc^\op,\circ)$, $\bfP=\FQA^\op$.

\vspace*{2mm}\noindent{\bf Multiplicative Manin matrices.}
 It is shown that the quantum representations of a quantum monoid on $\gX_B(\KK)$ is in one-to-one correspondence with the multiplicative $B$-Manin matrices over the corresponding bialgebra. In the case of more general quantum algebra an additional condition on the Manin matrix is imposed: its entries should belong to the first order graded component.

\vspace*{2mm}\noindent{\bf $(A,B)$-Manin matrices.}
 The multiplicative $B$-Manin matrices correspond to the comonoid morphisms from the graded algebra $\bfcoend\big(\gX_B(\KK)\big)$ equipped with the comultiplication and counit~\eqref{dvMMM}. A natural question arises: how to generalise the theory described here in order to include more general $(A,B)$-Manin matrices and cocomposition~\eqref{dKMN}? In this case we probably need to consider a family of Manin matrices compatible with the cocompositions and counits in some way. We plan to construct such a generalisation in future works.

\vspace*{2mm}\noindent{\bf Another subcategories of graded algebras.}
 The quantum linear spaces $\A\in\FQA^\op$ can be generalised by considering more wide subcategory of $\GrAlg$. The main condition we should keep is the connectivity: $\A_0=\KK$. For example, one can take all the connected affinely generated graded algebras -- this will allow to connect the (generalised) quantum linear spaces with a quantum analogue of the projective spaces. We hope to investigate this case in the future.

\vspace*{2mm}\noindent{\bf Super-case.}
 If one considers a vector space $V\in\FVect$ as a purely odd super-vector space, the embedding $\Lambda^*\colon\FVect\to\FQA^\op$ should be applied (see p.~\ref{bbSQVS}). By combining the even and odd embeddings one can embed an arbitrary finite-dimensional super-vector space $V=V_{\bar0}\oplus V_{\bar 1}$ into the category of quadratic super-algebras as $S^*(V_{\bar0})\otimes\Lambda^*(V_{\bar 1})$. The general approach of Section~\ref{sec3} can be applied to the super-case as well. The super-version of Quantum Representation Theory is also a subject for future publications.

\vspace*{2mm}\noindent{\bf $S$-embedding vs $T$-embedding.}
 Both functors $S^*$ and $T^*$ lift the vector spaces and representations to the quantum level. These two ways to embed the vector spaces have advantages and disadvantages with respect to each other. The $S$-embedding has more intuitive interpretation and leads to the usual Manin matrices, which have a lot of applications. It has a strong monoidal structure $(\FVect,\oplus)\to(\FQA^\op,\otimes)$, but it is not strong monoidal as a monoidal functor $(\FVect,\otimes)\to(\FQA^\op,\circ)$. The $T$-embedding has a structure of strong monoidal functor $(\FVect,\otimes)\to(\FQA^\op,\circ)$. As a monoidal functor $(\FVect,\oplus)\to(\FQA^\op,\otimes)$ it is not strong monoidal, however it has a strong monoidal structure $(\FVect,\oplus)\to(\FQA^\op,\amalg)$, so it is more natural to consider coproduct $\amalg$ as an analogue of direct sum of vector spaces for the quantum case, when we use the $T$-embedding (note that $\amalg$ is exactly the categorical product in $\FQA^\op$).

\vspace*{2mm}\noindent{\bf Tensor product.}
 If we translate a finite-dimensional algebra or an algebraic monoid/group to the quantum level by either functor $S^*$ or $T^*$, then the tensor product of quantum representations of the obtained quantum algebra always exists due to the `commutativity' of the classical algebra. For a more general comonoid $\OO\in\Comon(\GrAlg,\circ)$ this existence is not guaranteed even if $\OO\in\Comon(\FQA,\circ)$ or $\OO=\OO_\gR$ for a quantum monoid/group $\gR=(\gR,\Delta_\gR,\varepsilon_\gR)\in\Mon(\Alg,\otimes)$. It means that the category of quantum representations loses the tensor product under a quantum deformation. We conjecture that the tensor product can be restored in some important cases by means of deformation of the symmetric structure.
 
\vspace*{2mm}\noindent{\bf Infinite-dimensional case.}
 By considering the purely algebraic theory of quantum representations we are restricted by the finite-dimensional case by three reasons. First, the functor $(-)^*\colon\Vect\to\Vect^\op$ is not full (it is only faithful). Second, it is not lax monoidal with respect to $\otimes$ (only colax monoidal). Finally, the monoidal category $(\QA,\circ)$ is not coclosed (the proof of Prop.~\ref{PropQA} essentially uses the finiteness of the basis $(w_i)$). It seems that a generalisation of the results given here to the infinite-dimensional case is possible, if we introduce some kind of topology on vector spaces and algebras.

\end{document}